\documentclass[a4paper,11pt]{article}
\usepackage{amsmath}
\usepackage{amsthm}
\let\newemptytheorem\newtheorem
\usepackage{amssymb}
\usepackage{bm}
\usepackage[margin=0.5in]{caption}
\usepackage{dsfont}
\usepackage{enumitem}
\usepackage{graphicx}
\usepackage{indentfirst}
\usepackage{mathrsfs}
\usepackage{mathtools}
\usepackage{nameref}
\usepackage[new]{old-arrows}
\usepackage{stmaryrd}
\usepackage{svg}
\usepackage[x11names]{xcolor}
\usepackage{xparse}
\usepackage{mathabx} 
\usepackage{blindtext}
\usepackage{titlesec}
\usepackage{titletoc}
\usepackage{mdframed}
\usepackage{letltxmacro}
\LetLtxMacro\amsproof\proof
\LetLtxMacro\amsendproof\endproof
\usepackage{thmbox}

\AtBeginDocument{
  \LetLtxMacro\proof\amsproof
  \LetLtxMacro\endproof\amsendproof
}

\usepackage{sansmathfonts}
\usepackage[sfdefault]{atkinson}
\usepackage[T1]{fontenc}

\usepackage[colorlinks=true,
    citecolor=blue,
    linkcolor=blue,
    urlcolor=blue,
    linktoc=page,
    hyperindex=true,
    pdfcreator={}]{hyperref} 
    
    \hypersetup{
        pdfauthor={Naomi Andrew and Yassine Guerch and Camille Horbez and Richard D. Wade},
        pdftitle={Finiteness properties of centralisers in Out(Fn) and automorphisms of free-by-cyclic groups},
        pdfkeywords={Out(Fn), Centralisers, Free-by-cyclic groups, Automorphisms, Finiteness properties, Group actions on trees, Limiting trees, Trees of cylinders, Transverse coverings, Twists},
}

\usepackage{geometry}
\usepackage{setspace} 
\usepackage{enumitem}
\setlist{nosep}

\usepackage{graphics}
\normalsize
\usepackage{tcolorbox}

\newcounter{theo}[section]
\renewcommand{\thetheo}{\arabic{section}.\arabic{theo}}
\newenvironment{theo}[2][]{
\refstepcounter{theo}
\ifstrempty{#1}
{\mdfsetup{
frametitle={
\tikz[baseline=(current bounding box.east),outer sep=0pt]
\node[anchor=east,rectangle,fill=blue!15]
{\strut Theorem~\thetheo};}}
}
{\mdfsetup{
frametitle={
\tikz[baseline=(current bounding box.east),outer sep=0pt]
\node[anchor=east,rectangle,fill=blue!15]
{\strut Theorem~\thetheo:~#1};}}
}
\mdfsetup{innertopmargin=2pt,linecolor=black,
linewidth=1pt,topline=true,
frametitleaboveskip=\dimexpr-\ht\strutbox\relax
}
\begin{mdframed}[]\relax
\label{#2}}{\end{mdframed}}

\global\mdfdefinestyle{exampledefault}{topline=false,bottomline=false,rightline=false}

\theoremstyle{definition}
\newtheorem[S,bodystyle=\normalfont\noindent,titlestyle=\normalfont{~(#1)}]{thm}[theo]{Theorem}
\newtheorem[S,bodystyle=\normalfont\noindent,titlestyle=\normalfont{~(#1)}]{lem}[theo]{Lemma}
\newtheorem[S,bodystyle=\normalfont\noindent,titlestyle=\normalfont{~(#1)}]{lemma}[theo]{Lemma}

\newtheorem[S,bodystyle=\normalfont\noindent,titlestyle=\normalfont{~(#1)}]{coro}[theo]{Corollary}
\newtheorem[S,bodystyle=\normalfont\noindent,titlestyle=\normalfont{~(#1)}]{cor}[theo]{Corollary}
\newtheorem[S,bodystyle=\normalfont\noindent,titlestyle=\normalfont{~(#1)}]{prop}[theo]{Proposition}
\newtheorem[S,bodystyle=\normalfont\noindent,titlestyle=\normalfont{~(#1)}]{facts}[theo]{Facts}

\theoremstyle{remark}
\newemptytheorem{rk}[theo]{Remark}
\newemptytheorem{notation}[theo]{Notation}

\newemptytheorem{claim}[theo]{Claim}
\newemptytheorem{exer}[theo]{Exercise}
\newemptytheorem{q}[theo]{Question}
\theoremstyle{definition}
\newtheorem[S,bodystyle=\normalfont\noindent,titlestyle=\normalfont{~(#1)}]{defn}[theo]{Definition}
\newemptytheorem{ex}[theo]{Example}

\newtheorem[S,bodystyle=\normalfont\noindent,titlestyle=\normalfont{~(###1)}]{manualtheoreminner}{Theorem}

\renewcommand{\phi}{\varphi}

\renewcommand{\subset}{\subseteq}

\renewcommand{\subsetneq}{\varsubsetneq}

\newcommand {\cala} {{\mathcal {A}}}

\newcommand {\cale} {{\mathcal {E}}}   
\newcommand {\calf} {{\mathcal {F}}}   
   
\newcommand {\calh} {{\mathcal {H}}}

\newcommand {\calk} {{\mathcal {K}}}   
   
\newcommand {\calm} {{\mathcal {M}}}

\newcommand {\calp} {{\mathcal {P}}}

\newcommand {\caly} {{\mathcal {Y}}}

\newcommand{\ZZ}{\mathbb{Z}}
\newcommand{\NN}{\mathbb{N}}
\newcommand{\res}{\mathrm{res}}
\newcommand{\vpc}{\mathrm{vpc}}
\newcommand{\slow}{\mathrm{slow}}

 \usepackage{bbm}

\newcommand{\Isom} {{\mathrm{Isom}}}

\newcommand{\rank}{{\mathrm{rank}}}
\newcommand{\Stab} {{\mathrm{Stab}}}
\newcommand{\SetStab}{{\mathrm{SetStab}_M}}
\newcommand{\PtStab} {{\mathrm{PtStab}_M}}
\newcommand{\Slow} {{\mathrm{Slow}}}

\newcommand{\Fix}{\mathrm{Fix}}

\newcommand{\Out} {{\mathrm{Out}}}
\newcommand{\Hom} {{\mathrm{Hom}}}
\newcommand{\Aut} {{\mathrm{Aut}}}

\newcommand{\Inn} {{\mathrm{Inn}}}

\newcommand{\id} {\mathrm{id}}
\newcommand{\ad} {{\mathrm{ad}}}
\newcommand{\ab}{\mathrm{ab}}

\newcommand{\Comp}{\mathrm{Comp}}
\newcommand{\F}{\rm{F}}
\newcommand{\VF}{\rm{VF}}

\newcommand{\IA}{\mathrm{IA}_F(\mathbb{Z}/3\mathbb{Z})}

\newcommand{\Inc}{\mathrm{Inc}}

\newcommand{\edge}{\mathrm{edge}}
\newcommand{\vertex}{\mathrm{vertex}}

\newcommand{\gobble}[1]{} 

\newcommand{\per}{\mathrm{per}} 

\makeatletter
\edef\@tempa#1#2{\def#1{\mathaccent\string"\noexpand\accentclass@#2 }}
\@tempa\rond{017}
\makeatother

\usepackage{tikz-cd}

\title{Finiteness properties of centralisers in $\Out(F_N)$ and automorphisms of free-by-cyclic groups}
\author{Naomi Andrew \and Yassine Guerch \and Camille Horbez \and Richard D.\ Wade}
\date{\today}

\begin{document}

\maketitle

\begin{abstract}
    Let $F$ be a finitely generated free group. We prove that the centraliser of any element of $\Out(F)$ or $\Aut(F)$ is of type $\VF$. We also prove that the (outer) automorphism group of any free-by-cyclic group $M=F\rtimes_\varphi\mathbb{Z}$ is of type VF, and so are the McCool subgroups of $\Out(M)$. This relies on the construction of a new $\Out(M)$-invariant splitting of $M$, that enables us to reduce the proof of the above statements to the case where the outer class of $\varphi$ has finite order. 
\end{abstract}

\vspace{0.5cm}
\setcounter{tocdepth}{1}
\startcontents \printcontents{}{0}[1]{}

\vfill

\noindent \textcolor{black!70}{This work is openly licensed via \href{https://creativecommons.org/licenses/by/4.0/}{Creative Commons CC-BY 4.0}}
\medskip

\noindent \textcolor{black!70}{This work was completed without use of artificial intelligence. The authors do not consent for all or part of this work to be used for training of artificial intelligence models.}
\newpage

\section{Introduction}

Finite generation is at the heart of geometric group theory: when a group $G$ is finitely generated it has a well-defined word metric (up to quasi-isometry), and this word metric is the departure point for the exploration of a group as a geometric object. After finite generation comes finite presentability, and this is followed by a host of \emph{higher finiteness properties}, the strongest of which is for $G$ to have a classifying space which is a finite cell complex (in which case we say that $G$ is \emph{of finite type}). As many groups of interest have torsion and every group of finite type is torsion-free, in practice we instead ask for $G$ to be  \emph{of type $\VF$}, which means that $G$ has a finite-index subgroup of finite type.

The study of finiteness properties often appears hand-in-hand with the study of subgroup structure. Even when a group $G$ is of finite type, its subgroups need not be so well-behaved. A well-known manifestation of this phenomenon is the Bestvina--Brady subgroup of a right-angled Artin group: right-angled Artin groups are of finite type but their Bestvina--Brady subgroups can exhibit a wide range of finiteness properties \cite{Bestvina-Brady}. 

In this paper the main subgroups of concern are centralisers in (outer) automorphism groups of free groups. Centralisers are a naturally occurring family of subgroups and are intricately linked to the (famously open) conjugacy problem in $\Out(F)$ (where $F$ is a finite rank free group). Despite their importance, and many partial results (see below for a more detailed history), it was not known whether arbitrary centralisers in $\Out(F)$ were finitely generated. The first main goal of this paper is to show that they are not only finitely generated, they are of type $\VF$.

\begin{theo}{theo:centralisers-out}
Let $F$ be a finitely generated free group. 
\begin{enumerate}
\item\label{main-1} For every $\Phi\in\Out(F)$, the centraliser of $\Phi$ in $\Out(F)$ is of type $\mathrm{VF}$.
\item\label{main-2} For every $\varphi\in\Aut(F)$, the centraliser of $\varphi$ in $\Aut(F)$ is of type $\mathrm{VF}$. 
\end{enumerate}
\end{theo}

In the same way that passing to subgroups need not preserve finiteness properties, the same is true of passage to (outer) automorphism groups \cite{CollinsLevin1983}. When $F$ is a finitely generated free group, Nielsen proved that $\Out(F)$ was finitely generated in 1924 \cite{Nie}, however we did not know that $\Out(F)$ was of type $\VF$ until Culler and Vogtmann's introduction of Outer space, just over sixty years later \cite{CV86}. 

Given an automorphism $\varphi\in\Aut(F)$, one can form the semi-direct product $F\rtimes_\varphi\mathbb{Z}$, which only depends (up to isomorphism) on the outer class of $\varphi$. A \emph{free-by-cyclic group} is a group of the form $M=F\rtimes_\varphi\mathbb{Z}$, for some finitely generated free group $F$ and some $\varphi\in\Aut(F)$. Our second main theorem resolves Conjecture~7.1 of \cite{aimpl}.

\begin{theo}{theo:fbyz}
For every free-by-cyclic group $M$, the groups $\Out(M)$ and $\Aut(M)$ are of type $\mathrm{VF}$.
\end{theo}

\begin{rk}\label{rk:linton}
Although in this paper we adopt the convention that free-by-cyclic groups have finitely generated fibres, one can use work of Linton~\cite{Linton} to show that the above theorem holds in the more general case that $M$ (and not necessarily $F$) is finitely generated, see Section~\ref{sec:linton}. 
\end{rk} 

These two theorems are related, and in fact both follow from a more general statement, Theorem~\ref{theo:mccool} below, regarding McCool subgroups of automorphism groups of free-by-cyclic groups. We will explain how they are linked in Section~\ref{sec:main}. Before that, we first recall some of the history and further motivation behind these questions.

\subsection{Some history}

\paragraph*{1 -- Centralisers in $\Out(F)$.} Work on understanding the behaviour of subgroups of $\mathrm{Out}(F)$ has led not only to interesting results but also to the development of much of the theory underpinning modern study. This is embodied by the resolution of the Tits alternative \cite{BFHTits1,BFHTits2,BFHTits3}, which expanded Bestvina and Handel's train-track theory and highlighted the importance of laminations and limiting trees in $\Out(F)$ (we will see the latter play a crucial role in the sequel).

In terms of subgroup structure away from the Tits alternative, finite subgroups of $\Out(F)$ have an explicit description as stabilisers of points in Outer space \cite{Cul,Khr,Zim} and Feighn--Handel gave tools to describe  abelian subgroups in \cite{FH}. Centralisers are the next natural step in the quest to investigate commutation in $\Out(F)$: this has been achieved so far in some specific cases described below. As we briefly mentioned above, the behaviour of centralisers is closely related to the conjugacy problem. While the conjugacy problem is solved in a host of specific cases \cite{Sel95,Khr95,Los,CL99,KLV01,Lustig,dahmani2016suspensions,DFMT,DT24,FeiHan25,bartlett2026}, it remains a major open problem in general. 

To organise this overview of previous results on centralisers we use a standard growth dichotomy for outer automorphisms $\Phi\in\Out(F)$: either $\Phi$ has polynomial growth (all conjugacy classes grow polynomially fast under iteration of $\Phi$) or it has  exponential growth (some conjugacy class grows exponentially fast). This dichotomy is crucially related to the geometry of the extension $M=F\rtimes_\Phi\mathbb{Z}$, namely: $M$ is relatively hyperbolic if and only if $\Phi$ is exponentially growing \cite{gaulus,DahLi,Ghosh,Hagen}. A further subset of previous results concern low-rank phenomena, which we treat separately at the end of this section. 

\textbf{$\drsh$ Polynomial growth.} The identity automorphism is the simplest example with polynomial growth. In \cite{CV86}, Culler and Vogtmann proved that $\Out(F)$ itself is of type $\VF$, by constructing a contractible simplicial complex, namely the spine of \emph{Outer space}, on which $\Out(F)$ acts with finite stabilisers (since $\Out(F)$ is virtually torsion free~\cite{BT}, this shows that it is of type $\VF$). Later, Krsti\'c and Vogtmann showed that centralisers of finite-order elements of $\Out(F)$ are also of type $\VF$, by proving the contractibility of their fixed point sets in Outer space \cite{KV93}. In \cite{RW}, Rodenhausen and Wade proved that the centraliser of any Dehn twist automorphism (in the sense of Cohen and Lustig \cite{CL99}) is of type $\VF$. Their result was later extended by Rodenhausen \cite{Rod} to all unipotent polynomially-growing automorphisms.

However, replacing $\Phi$ with a root can change the commensurability class of the centraliser. One example of this phenomenon is given by finite-order irreducible elements of $\Out(F)$, whose centralisers are finite (see \cite[Corollary~3.10]{FMSmindisplaced}). 
An infinite-order example involving a root of a Dehn twist is given in  \cite[Section 2]{AMCentralisers}. 
In that direction, Andrew and Martino \cite{AMCentralisers} extended type $\mathrm{VF}$ to all linearly-growing automorphisms (allowing roots of Dehn twists). The extension to all polynomially-growing automorphisms is new in this paper.

\textbf{$\drsh$ Exponential growth.} On the other side of the growth dichotomy, it follows from Bestvina, Feighn and Handel's work on laminations for irreducible elements with irreducible powers \cite{BFHlaminations} that their centralisers are virtually cyclic, and some examples where the centraliser is isomorphic to $\mathbb{Z}$ were given in \cite{AKP}. For automorphisms that are atoroidal, or where there are finitely many conjugacy classes of periodic cyclic subgroups, the centraliser is virtually abelian \cite{FH} (see also \cite{GuerchRoots} and \cite{MutanguhaPretreesUniv}).  In this case, the relative hyperbolicity of the extension $F\rtimes_\Phi\mathbb{Z}$ provides useful tools including those coming from JSJ theory~\cite{RS,GLjsj}.

\textbf{$\drsh$ Low rank structure theorems.} When $F$ has rank $2$, the group $\Out(F)$ is virtually free, so centralisers of infinite-order elements are all virtually cyclic. When $F$ has rank $3$, finite generation of all centralisers in $\Out(F)$ is due to Francaviglia, Martino and Syrigos \cite[Theorem~8.2.1]{FMScocompact}, and type $\mathrm{VF}$ was proved by Guerch, Hughes and S\'anchez Salda\~na \cite{GHSS} for all centralisers of automorphisms in a certain finite-index subgroup of $\Out(F_3)$.

\paragraph*{2 -- Automorphisms of free-by-cyclic groups.}

Despite their tight relationships to centralisers, reviewed in the next section, automorphism groups of free-by-cyclic groups have received less attention so far. Let $\Phi\in\Out(F)$, and let $M$ be the associated free-by-cyclic group.

\textbf{$\drsh$ Polynomial growth.} If $\Phi$ has finite order in $\Out(F)$, then $\Out(M)$ is of type $\mathrm{VF}$ as a consequence of Levitt's work on automorphisms of generalised Baumslag--Solitar groups \cite{LevittGBS}. And if $\Phi$ is linearly growing, then Andrew--Martino proved that $\Out(M)$ is finitely generated \cite{AMFbyZ}.

\textbf{$\drsh$ Exponential growth.} In this case $M$ is relatively hyperbolic. If $\Phi$ is atoroidal, then $M$ is hyperbolic \cite{BFcombination,Bri}, and it follows from JSJ theory that $\Out(M)$ is virtually abelian, noting that the canonical JSJ decomposition of $M$ has no QH vertex~\cite[Proposition~2.14]{dahmani2016suspensions}.

\textbf{$\drsh$ Low rank structure theorems.} The outer automorphism group $\Out(M)$ was described by Bogopolski--Martino--Ventura when $F$ has rank $2$ \cite{BMV}, and was proven to be finitely generated by Andrew--Martino when $F$ has rank $3$ \cite{AMFbyZ}.

\subsection{From centralisers to free-by-cyclic groups: main theorem}\label{sec:main}

For both the conjugacy problem \cite{Sel95,dahmani2016suspensions} and the study of centralisers \cite{AMCentralisers}, it turns out to be useful to switch the viewpoint to that of free-by-cyclic groups. Indeed, assuming that $F$ is non-abelian, let $\Phi\in\Out(F)$, let $\varphi\in\Aut(F)$ be a representative of the outer class $\Phi$, and let $M=F\rtimes_{\varphi}\mathbb{Z}$ be the associated free-by-cyclic group. Let $\Aut(M,F)\subseteq\Aut(M)$ be the subgroup consisting of all automorphisms that preserve the fibre $F$, and let $\Aut^+(M,F)$ be the subgroup acting trivially on the quotient $M/F\simeq\mathbb{Z}$. The homomorphism $\Aut^+(M,F)\to\Aut(F)$ given by restricting an automorphism to the fibre $F$ induces an injective homomorphism $\Aut^+(M,F)/\Inn(F)\to\Out(F)$, whose image coincides with the centraliser $C_{\Out(F)}(\Phi)$. This observation, already exploited in \cite{AMCentralisers}, is recorded in Section~\ref{sec:autos} of the present paper. 

The above observation leads us to study fibre-preserving automorphisms of free-by-cyclic groups, and our techniques also enable us to study the full automorphism group. For induction purposes, it is also useful to study \emph{McCool subgroups} of these: given a finite set $\calh$ of conjugacy classes of subgroups of $M$, we let $\Out(M,\calh^{(t)})$ be the subgroup of $\Out(M)$ consisting of all outer automorphisms that act trivially (i.e.\ by conjugation by an element $g_H\in M$) on every subgroup $H$ whose conjugacy class belongs to $\calh$. Likewise, we denote by $\Out(M,F)$ the image of $\Aut(M,F)$ in $\Out(M)$, and by $\Out(M,F,\calh^{(t)})$ its intersection with $\Out(M,\calh^{(t)})$. The main theorem of the present paper is as follows.

\begin{theo}{theo:mccool}
Let $\varphi\in\Aut(F)$, let $M=F\rtimes_{\varphi}\mathbb{Z}$, and let $\calh$ be a finite set of conjugacy classes of finitely generated subgroups of $M$. 

Then $\Out(M,\calh^{(t)})$ and $\Out(M,F,\calh^{(t)})$ are of type $\VF$.
\end{theo}

Theorems~\ref{theo:centralisers-out} and~\ref{theo:fbyz} both follow from Theorem~\ref{theo:mccool}. For the former, this uses the above observation relating centralisers in $\Out(F)$ to automorphisms of free-by-cyclic groups (and its version for centralisers in $\Aut(F)$), see Section~\ref{sec:autos}.

\subsection{Canonical trees and proof strategy}

Our proof of Theorem~\ref{theo:mccool} relies on the construction of canonical (i.e.\ $\Aut$-invariant) actions of $M$ on trees. A central construction in the paper is an $\Aut(M)$-invariant action on a tree $T$ whose vertex groups either contain a finite-index subgroup isomorphic to $F'\times\mathbb{Z}$ for some finitely generated free group $F'$, or have few automorphisms, in a sense made precise in our next statement. 
Given an $\Aut(M)$-invariant $M$-tree $T$, we let $\Out^0_T(M)$ be the finite-index subgroup of $\Out(M)$ consisting of outer automorphisms that act trivially on the quotient graph $T/M$. For every vertex $v\in VT$, we denote by $M_v$ the stabiliser of $v$ under the action of $M$. Then there is a \emph{restriction homomorphism} $\mathrm{res}_v:\Out^0_T(M)\to\Out(M_v)$, defined as follows (see \cite{Lev}): every $\Phi\in\Out^0_T(M)$ has a representative $\varphi\in\Aut(M)$ that fixes $v$, and we let $\mathrm{res}_v(\Phi)$ be the outer class of $\varphi|_{M_v}$ (this is well defined because any two such representatives differ by an inner automorphism of $M_v$). 

\begin{theo}[see Theorem~\ref{theo:dévissage}]{theo:tree}
Let $M$ be a free-by-cyclic group. Then there exists an $\Aut(M)$-invariant $M$-tree $S$ with finitely generated edge stabilisers, such that 
\begin{itemize}
\item the $M$-action on $S$ is acylindrical;
\item for every vertex $v\in VS$, one of the following two holds:
\begin{itemize}
\item $M_v$ has a finite-index subgroup isomorphic to $F_N\times\mathbb{Z}$ for some finitely generated free group $F_N$; 
\item the image of the restriction homomorphism $\mathrm{res}_v:\Out^0_S(M)\to\Out(M_v)$ is virtually polycyclic.
\end{itemize}
\end{itemize}
\end{theo}

Theorem~\ref{theo:tree} essentially reduces the proof of Theorem~\ref{theo:mccool} to the case where $M$ is virtually free-times-cyclic. Along the way we show how to resolve the virtually free-times-cyclic case using earlier works of Krsti\'c--Vogtmann \cite{KV93} and Bestvina--Feighn--Handel \cite{BFH2023} on the structure of centralisers of finite-order outer automorphisms of $F$. Passing from Theorem~\ref{theo:tree} to Theorem~\ref{theo:mccool} also requires us to have a control on the subgroup of outer automorphisms of $M$ that preserve $T$ and act trivially on every vertex group. Using the acylindricity of the $M$-action on $T$, we show in Section~\ref{sec:twist} that this coincides with the \emph{group of twists} of the splitting $T$ in the sense of Levitt \cite{Lev}. And, provided that $M$ has trivial centre, we then prove that for any $M$-tree with finitely generated edge stabilisers, there exists a finite-index subgroup of $\Out(M)$ whose intersection with the group of twists is of type $\mathrm{F}$, see Corollary~\ref{cor:twist-type-F}.

\paragraph*{On the proof of Theorem~\ref{theo:tree}.} Again there are several regimes to consider, depending on the growth of the outer automorphism $\Phi$ that defines the free-by-cyclic group. 

\smallskip

If $\Phi$ is periodic, then $M$ is virtually free-times-$\mathbb{Z}$. In this case the theorem is vacuously true, with $S$ reduced to one vertex.

\smallskip

In the polynomially-growing case (but not virtually free-times-$\mathbb{Z}$), our construction of the tree $S$ in Theorem~\ref{theo:tree} has three main steps.

$\drsh$ \textbf{Polynomial case, Step 1} (see Section~\ref{sec:limiting-tree}). The first step is based on a limit argument, following works of Cohen--Lustig \cite{CL99} in the case of linear growth, and of Bestvina--Feighn--Handel \cite{BFHTits2} for quadratic growth and higher. By iterating $\Phi$ on a tree in Outer space and taking a limiting tree, one obtains a $\Phi$-invariant simplicial $F$-tree whose edge stabilisers are either infinite cyclic (for linear growth) or trivial (for quadratic growth and higher). A collapsed tree of cylinders construction gives us a first $\Aut(M)$-invariant tree $T$ in which an element of $F$ is elliptic if and only if its conjugacy class does not grow at the top speed under iteration of $\Phi$. We would now like to blow-up these slow-growing vertices into their own limiting trees, but for this we would need the incident edge groups to be elliptic in these trees, which may fail.

$\drsh$ \textbf{Polynomial case, Step 2} (see Section~\ref{sec:tree-of-cylinders}). To bypass this difficulty, we use a transverse covering construction. This construction is completely new, and we consider it to be one of the most novel and important ideas introduced in this paper. Given an $\Aut(M)$-invariant $M$-tree $T$ and a finite set $\calh$ of conjugacy classes of finitely generated subgroups of $M$, we build a new $\Aut(M,\calh)$-invariant tree $U_{T,\calh}$ in which all subgroups whose conjugacy class belongs to $\calh$ become elliptic. The idea is to consider the skeleton of the transverse covering of $T$ by the minimal subtrees $T_H$, where $H$ varies among all subgroups whose conjugacy class belongs to $\calh$ and which are not elliptic in $T$. Therefore, vertices of $U_{T,\calh}$ are either vertices of the original tree $T$, or else they encode the overlapping patterns of the subtrees $T_H$. In Section~\ref{sec:tree-of-cylinders}, we investigate the vertex and edge stabilisers of this skeleton, and the image of $\Out^0_T(M,\calh^{(t)})$ in the outer automorphism groups of the edge and vertex groups of $U_{T,\calh}$. In particular, we prove in Proposition~\ref{prop:tree-cyl-F-by-Z-case} that, if $v \in VU_{T,\calh}$ is not a vertex of the original tree $T$, then the image of $\Out^0_T(M,\calh^{(t)})$ in $\Out(M_v)$ is virtually polycyclic. 

$\drsh$ \textbf{Polynomial case, Step 3} (see Section~\ref{sec:dévissage}). Let $T$ be the tree obtained in Step~1, let $v\in VT$, and let $\calh$ be the set of all incident edge stabilisers. Let $T_v$ be a limiting tree for $M_v$ (obtained as in Step~1 but in the vertex group $M_v$, of smaller growth). By applying Step~2 to $T_v$ and $\calh$, we get a new splitting of $M_v$ where all subgroups whose conjugacy class belongs to $\calh$ are elliptic. We can therefore blow-up the vertex $v$ of $T$ using this $M_v$-tree. By repeating this procedure finitely many times until the vertex groups have growth rate $0$ (i.e.\ they are virtually free-times-$\mathbb{Z}$), we reach a tree $S$ as in the statement of Theorem~\ref{theo:tree}.

\smallskip

Finally, if $\Phi$ has exponential growth, then $M$ is relatively hyperbolic, relative to suspensions of maximal polynomially-growing subgroups, and our goal is to reduce to the polynomially-growing case. The group $M$ has a canonical JSJ tree $T$ relative to its peripheral subgroups. Some vertex groups of $T$ are polynomially-growing free-by-cyclic groups (the others being rigid with finite outer automorphism group relative to incident edge groups), and we want to blow-up each such vertex $v$ into its associated canonical tree $S_v$ (built above). Before this, as in Step~2 above, we need to pass to a skeleton of a transverse covering to make the incident edge groups elliptic. The details are carried in Section~\ref{sec:dévissage}.

\paragraph*{Organisation of the paper.} After a general background section (Section~\ref{sec:background}), we review the relationship between centralisers in $\Out(F)$ or $\Aut(F)$ and automorphisms of $F\rtimes\mathbb{Z}$ in Section~\ref{sec:autos}. In Section~\ref{sec:finite-order}, relying on earlier works of Krsti\'c--Vogtmann and Bestvina--Feighn--Handel, we prove our main theorems in the case where $M$ is virtually free-times-$\mathbb{Z}$ (equivalently, the monodromy has finite order in $\Out(F)$). Section~\ref{sec:arboriculture} gives a general review of group actions on trees. The next three sections aim at establishing Theorem~\ref{theo:tree}, constructing our canonical tree for free-by-cyclic groups, following the outline presented right above. Once Theorem~\ref{theo:tree} is proved, we still need to understand the group of twists of the tree we just built (Section~\ref{sec:twist}), and extend automorphisms of the vertex groups to automorphisms of $M$ (Section~\ref{sec:extend}). The proof of our main theorem (Theorem~\ref{theo:mccool}, from which Theorems~\ref{theo:centralisers-out} and~\ref{theo:fbyz} follow) is then completed in Section~\ref{sec:conclusion}. The final section uses work of Linton to extend our main results to finitely generated groups of the form $F\rtimes\mathbb{Z}$, where the free group $F$ is no longer assumed finitely generated, as explained in Remark~\ref{rk:linton}.

\paragraph*{Acknowledgments.} We give our deepest thanks to Monika Kudlinska for her active participation and numerous decisive discussions in early stages of this project, and Gilbert Levitt for a three-day discussion in Caen that led to the covering argument construction presented in Section~\ref{sec:tree-of-cylinders} of the present work. We thank Gabriel Bartlett, François Dahmani and Vincent Guirardel for their interest in our work and interesting discussions. Finally, we thank Martin Bridson, Armando Martino and Moritz Rodenhausen for many nurturing discussions and the ideas that emerged from previous collaborations on this topic.

\section{General background}\label{sec:background}

A \emph{free-by-cyclic group} is a group $M=F\rtimes_{\varphi}\ZZ$, where $F$ is a finitely generated free group and $\varphi\in\Aut(F)$ -- note that by convention we impose $F$ to be finitely generated, possibly trivial or cyclic. We say that $F$ is a \emph{fibre} of $M$, and that $\varphi$ (or its outer class $\Phi\in\Out(F)$) is a \emph{monodromy} of $M$. We mention that the isomorphism class of $M$ only depends on the outer automorphism $\Phi$. Therefore, we will sometimes allow ourselves to denote a free-by-cyclic group by $F\rtimes_\Phi\mathbb{Z}$.

In what follows, we fix a free-by-cyclic group $M$, a fibre $F$, a  monodromy $\varphi\in\Aut(F)$ and its outer class $\Phi$. We let $\chi:M\to M/F\simeq\mathbb{Z}$ be the associated character. We also choose a generator $s\in M$ of a section of $\ZZ$ in $M$, such that the conjugation action of $s$ on $F$ induces the automorphism $\varphi$, i.e.\ $sgs^{-1}=\varphi(g)$ for every $g\in F$. In particular $M=\langle F,s\rangle$.

\subsection{Growth of outer automorphisms of free groups}\label{sec:growth_prelim}

Fix a finite generating set $S$ of $F$, and denote by $|\cdot|_S$ the associated word length. Given an element $g\in F$, the \emph{conjugacy length} of $g$ with respect to $S$ is $||g||_S=\inf_{h\in F}|hgh^{-1}|_S$. Given $\Phi\in\Out(F)$, it follows from train-track theory \cite{BH} (see also \cite{Levittgrowth}) that the conjugacy length of every element $g\in F$ grows either at most polynomially fast, or at least exponentially fast, under iteration of $\Phi$ (in fact $||\Phi^n(g)||_S$ always grows, up to a multiplicative and additive error, like $n^d\lambda^n$ for some $d\in\mathbb{N}$ and $\lambda\ge 1$, and $d$ and $\lambda$ do not depend on the choice of the generating set $S$). We say that $\Phi$ is \emph{polynomially growing} if every conjugacy class grows polynomially fast under iteration of $\Phi$, and \emph{exponentially growing} otherwise. If $\Phi$ is polynomially growing, its \emph{degree} is the maximal $d \geq 0$ such that there exists $g \in F$ whose conjugacy length grows like $n^d$ under iteration of $\Phi$. The degree of $\Phi$ is bounded by the rank of $F$, see \cite[Theorem~1]{Levittgrowth}.

We say that the free-by-cyclic group $M$ \emph{has polynomially-growing monodromy} if the outer class $\Phi$ is polynomially growing, and that it \emph{has exponentially-growing monodromy} if $\Phi$ is exponentially growing. We say that $M$ \emph{has finite-order monodromy} if $\Phi$ is a finite order element of $\Out(F)$. These notions are independent of the choice of a monodromy, see e.g.\ \cite[Theorem~3.2]{MutanguhaPG}, and Proposition~\ref{prop:periodic_implies_centre} below. 

\begin{thm}[\cite{Macura}] 
\label{thm:Macura}If $M$ and $M'$ are quasi-isometric free-by-cyclic groups with polynomially-growing monodromies $\Phi$ and $\Phi'$ respectively, then the degrees of growth of $\Phi$ and $\Phi'$ are equal. 
\end{thm}
 
Let $\calp(\Phi)$ be the set of $F$-conjugacy classes of maximal non-trivial subgroups $P$ of $F$ for which there exist $k>0$ and $\alpha\in\Aut(F)$ representing the outer class $\Phi^k$, such that $\alpha(P)=P$ and $\alpha|_{P}$ is polynomially growing. The set $\calp(\Phi)$ is a $\Phi$-invariant finite set of $F$-conjugacy classes of finitely generated malnormal subgroups of $F$ (see \cite[Proposition~1.4]{Levittgrowth}). For every subgroup $P\subseteq F$ whose conjugacy class belongs to $\calp(\Phi)$, let $k_P >0$ be the size of the orbit of the conjugacy class of $P$ under $\Phi$, and let $\phi_P\in\Aut(G)$ representing the outer class $\Phi^{k_P}$. 
Let $\calp(M)$ be the set of all $M$-conjugacy classes of subgroups of the form $P \rtimes_{\phi_P} \ZZ$, as $P$ varies among all subgroups whose conjugacy class belongs to $\calp(\Phi)$. This is a finite set of conjugacy classes of non-cyclic subgroups of $M$. As the notation suggests, the set $\calp(M)$ does not depend on the monodromy $\Phi$ (see for instance~\cite[Proposition~3.6]{MutanguhaPG}). 

\begin{thm}[\cite{gaulus,DahLi,Ghosh}]
\label{thm:relhyp}Any free-by-cyclic group $M$ is hyperbolic relative to the set $\calp(M)$.
\end{thm}

The following lemma allows us to compare the growth of elements under iteration of $\Phi$ or under iteration of a restriction of $\Phi$ to a finitely generated subgroup of $F$.

\begin{lem}
\label{lem:growth_subgroups}Let $\Phi \in \Out(F)$. Let $K \subseteq F$ be a finitely generated subgroup, let $k \geq 1$ and let $\psi$ be a representative of $\Phi^k$ preserving $K$. Let $\Psi$ be the outer class of $\psi|_{K}$, an element of $\Out(K)$. 

For every $g \in K$, the element $g$ has polynomial growth $d$ under iteration of $\Phi$ if and only if $g$ has polynomial growth $d$ under iteration of $\Psi$. In particular, if $\Phi$ is polynomially growing, so is $\Psi$.
\end{lem}

\begin{proof}
Let $S$ be a finite generating set of $F$ containing a generating set $S_K$ of $K$. Let $g \in K$. Since $K \subseteq F$, for every $n \geq 0$, we have $\lVert \Phi^{kn}(g) \rVert_S \leq \lVert \Psi^{n}(g) \rVert_{S_K}$. Therefore, if $g$ has polynomial growth $d$ under iteration of $\Psi$, then it has polynomial growth at most $d$ under iteration of $\Phi$.

Conversely, since $K$ is finitely generated, it is a virtual retract of $F$. Therefore, there exists $C\geq 1$ such that, for every $n \geq 0$, we have $\lVert \Psi^{n}(g) \rVert_{S_K} \leq C \lVert \Phi^{kn}(g) \rVert_S+C$. In particular, if $g$ has polynomial growth $d$ under iteration of $\Phi$, then it has polynomial growth at most $d$ under iteration of $\Psi$.
\end{proof}

\subsection{Automorphisms of free-by-cyclic groups} 

We denote by $\Aut(M,F)$ the subgroup of $\Aut(M)$ consisting of automorphisms $\alpha$ that are \emph{fibre preserving}, i.e.\ satisfy $\alpha(F)=F$. We let $\Out(M,F)$ be its image in $\Out(M)$. Every fibre-preserving automorphism of $M$ induces an automorphism of $M/F\simeq\mathbb{Z}$. We let $\Aut^+(M,F)$ be the \emph{character-preserving subgroup} of $\Aut(M,F)$, consisting of automorphisms that induce the trivial automorphism of $M/F$; it has index at most $2$. We let $\Out^+(M,F)$ be the image of $\Aut^+(M,F)$ in $\Out(M)$.

We will frequently use the following lemma. 

\begin{lemma}
\label{lemma:non-abelian-fibre}Let $M$ be a free-by-cyclic group, let $H\subseteq M$ be a subgroup such that $H\cap F$ is non-abelian, and let $\alpha\in\Aut^+(M,F)$.

If $\alpha|_{H\cap F}=\mathrm{id}_{H\cap F}$, then $\alpha|_{H}=\mathrm{id}_{H}$.
\end{lemma}

\begin{proof}
Let $u\in H$ be such that $H=\langle H\cap F,u\rangle$. Since $u$ normalises $H\cap F$, for every $h\in H\cap F$, we have $uhu^{-1}=\alpha(uhu^{-1})=\alpha(u)h\alpha(u)^{-1}$, showing that $\alpha(u)^{-1}u$ centralises $H\cap F$. Since $\alpha$ preserves the character, we have $\alpha(u)^{-1}u\in F$. Since $H\cap F$ is a non-abelian free group, it follows that $\alpha(u)=u$, so $\alpha$ is the identity on $H$.
\end{proof}

In particular, this lemma (applied with $H=M$) shows that if $F$ is non-abelian then an automorphism of $M$ acting as the identity on $F$ is the identity. 

\paragraph*{Relative (outer) automorphism groups.}

Given a set $\calh$ of conjugacy classes of subgroups of $M$, we denote by $\Aut(M,\calh)$ the subgroup of $\Aut(M)$ preserving the set $\calh$. Let $\Aut(M,\calh^{(t)})$ be the subgroup of $\Aut(M,\calh)$ consisting of automorphisms that coincide with a conjugation (by an element $g_H\in M$) in restriction to every subgroup $H\subseteq M$ whose conjugacy class belongs to $\calh$. We set \[\Aut(M,F,\calh)=\Aut(M,F) \cap \Aut(M,\calh);\] \[ \Aut(M,F,\calh^{(t)})=\Aut(M,F) \cap \Aut(M,\calh^{(t)}).\] We denote by $\Out(M,\calh)$, $\Out(M,\calh^{(t)})$, $\Out(M,F,\calh)$ and $\Out(M,F,\calh^{(t)})$ their respective images in $\Out(M)$. When $\calh=\{[H]\}$ consists of a single conjugacy class $[H]$ (of some subgroup $H\subseteq M$), we will allow ourselves to write $\Aut(M,[H])$ instead of $\Aut(M,\{[H]\})$ and likewise for $\Out$. 

\subsection{Free-by-cyclic groups with finite order monodromy}

We will often use the following standard proposition -- this appears for instance in various parts of \cite{LevittGBS,LevittRankOfGBS}, and we include a proof for the convenience of the reader.

\begin{prop} 
\label{prop:periodic_implies_centre}Let $M$ be a free-by-cyclic group with monodromy $\varphi$ representing $\Phi \in \Out(F)$. The following are equivalent:
\begin{enumerate}
    \item\label{center-1} The outer automorphism $\Phi$ has finite order.
    \item\label{center-2} The centre $Z(M)$ is non-trivial.
    \item\label{center-3} $M$ has a finite-index subgroup containing $F$ and isomorphic to $F \times \mathbb{Z}$. 
    \item\label{center-4} $M$ is virtually free-times-cyclic.
    \item\label{center-5} $M/Z(M)$ is a virtually free group.
\end{enumerate} 
If $M$ satisfies the above conditions and is not isomorphic to $\ZZ^2$, then $Z(M) = \langle as^k \rangle$, where $k$ is the order of $\Phi$ and $a\in F$ is uniquely determined and fixed by $\varphi$.
\end{prop}

\begin{rk} \label{rem:periodic_implies_GBS}
These are also equivalent to $M$ being a generalised Baumslag--Solitar group with trivial modulus, see \cite[Proposition~4.1]{LevittRankOfGBS}.
\end{rk}

\begin{proof}
If $F=\{1\}$ or $F\simeq\mathbb{Z}$ then $M$ is isomorphic to $\mathbb{Z}, \mathbb{Z}^2$, or the Klein bottle group and the statements (\ref{center-1})--(\ref{center-5}) all hold true. If $M=\mathbb{Z}$ then $M=Z(M)=\langle s \rangle$. If $M$ is isomorphic to the fundamental group of a Klein bottle, it has a presentation of the form $\langle a,s \;|\; sas^{-1}=a^{-1} \rangle$, with $F=\langle a\rangle$. One then checks that $Z(M)=\langle s^2\rangle$, and the conclusion of the proposition follows.

We may therefore assume that $F$ is non-abelian. Then $Z(M) \cap F=\{1\}$, so $\chi\colon M \to \ZZ$ is injective when restricted to $Z(M)$. In particular $Z(M)$ is either trivial or isomorphic to $\ZZ$. 

\underline{(\ref{center-1})$\Rightarrow$(\ref{center-2})}: If $\Phi^l$ is trivial then there exists $c \in F$ such that $\ad_c \varphi^l=\mathrm{id}$. Hence conjugation by $cs^l$ is the identity on $F$, so is also the identity on $M$ by Lemma~\ref{lemma:non-abelian-fibre}. So $cs^l \in Z(M)$.

\underline{(\ref{center-2})$\Rightarrow$(\ref{center-3})}: Suppose $Z(M)$ is non-trivial. Let $a\in F$ and $k\ge 0$ be such that $as^k$ generates $Z(M)$. Then $k \neq 0$ as $F$ is non-abelian, so $\langle F, as^k \rangle$ is a finite-index subgroup of $M$ isomorphic to $F \times \mathbb{Z}$.

\underline{(\ref{center-3})$\Rightarrow$(\ref{center-4})}: This is clear.

\underline{(\ref{center-4})$\Rightarrow$(\ref{center-1})}: If $M'=F' \times \mathbb{Z}$ is a finite-index subgroup of $M$, then we can pick an element $as^k \in Z(M')$. Every element $b \in F$ has a power $b^l \in M'$, so that $as^k$ commutes with $b^l$ and therefore $\ad_a \phi^k$ fixes $b^l$, and hence fixes $b$ by the unique root property in $F$. Hence $\ad_a \phi^k$ is the identity automorphism and $\Phi^k$ is trivial.

\underline{(\ref{center-2})$\Rightarrow$(\ref{center-5})}: Since (\ref{center-2}) implies (\ref{center-4}), $M$ has a finite-index subgroup $M'$ isomorphic to $F' \times \mathbb{Z}$. Then $Z(M)\cap M'$ is an infinite subgroup of $M'$ contained in $\mathbb{Z}$, so $M'$ projects to a finite-index subgroup of $M/Z(M)$ isomorphic to $F' \times \mathbb{Z}/k\mathbb{Z}$. 

\underline{(\ref{center-5})$\Rightarrow$(\ref{center-2})}: If $M/Z(M)$ is virtually free, then $Z(M)$ is non-trivial as $M$ is not virtually free (for instance, $M$ has a finitely generated infinite index infinite normal subgroup). 

\medskip 

For the last assertion of the proposition, we first prove that if $Z(M) = \langle as^k \rangle$, then $k$ is the order of $\Phi$. This follows from the observation that an element $cs^l$ (with $c\in F$) belongs to $Z(M)$ if and only if the automorphism $\ad_c\varphi^l$ is trivial (see the proof of (\ref{center-1})$\Rightarrow$(\ref{center-2})). 
Finally, the element $a$ is unique as $Z(M) \cap F=\{1\}$. Note that $[a,s]=[as^k,s]=1$ as $as^k$ is central. Hence $\varphi(a)=a^s=a$, showing that $a$ is fixed by $\varphi$. 
\end{proof}

\begin{lemma}
\label{lemma:trivial-center}Let $M$ be a free-by-cyclic group. Then $M/Z(M)$ has trivial centre. 
\end{lemma}

\begin{proof}
This is clear if $M$ has trivial centre or if $M\simeq\mathbb{Z}^2$. We now assume otherwise. By Proposition~\ref{prop:periodic_implies_centre}, the centre $Z(M)$ is infinite cyclic, generated by an element $z\in M$ with $\chi(z)\neq 0$.

Let $u\in M$, and let $\overline{u}$ be its image in $M/Z(M)$. If $\overline{u}$ is central in $M/Z(M)$, then for every $g\in M$, we have $[g,u]\in Z(M)$. Applying $\chi$ shows that $[g,u]$ is trivial. Therefore $u\in Z(M)$, and $\overline{u}$ is trivial.
\end{proof}

Throughout, we refer to a free-by-cyclic group satisfying the equivalent conditions of Proposition~\ref{prop:periodic_implies_centre} as having finite-order monodromy, and we remind the reader that this refers to its outer class.

For future use, we will also record the existence of some well-understood splittings of these groups:

\begin{prop} 
\label{prop:periodic_trees}Let $M$ be a free-by-cyclic group with finite-order monodromy. There exists a splitting $T$ of $M$ such that: 
\begin{itemize}
    \item Edge and vertex stabilisers of $T$ are finitely generated abelian groups.
    \item $Z(M)$ acts on $T$ trivially.
     \item Any fibre of $M$ acts on $T$ cocompactly.
\end{itemize}
\end{prop}

\begin{proof}
If $M$ is abelian then it is enough to take $T$ to be a point. If $M\simeq K$ is the fundamental group of a Klein bottle then $K=\langle a,s \;|\; sas^{-1}=a^{-1} \rangle$ surjects onto \[D_\infty \simeq \langle a,s \;|\; sas^{-1}=a^{-1}, s^2=1 \rangle,\] and as $Z(K) = \langle s^2 \rangle$ we can take $T$ to be the induced action on the line coming from $D_\infty$. In the remaining cases the rank of any fibre is at least 2 and one can find $T$ by either viewing $M$ as a generalised Baumslag--Solitar group (Remark~\ref{rem:periodic_implies_GBS}) or by taking a fixed point for the action of the monodromy on the Outer space for a fibre. In this case $T$ is not a line, so $Z(M)$ acts trivially, and any fibre will act on $T$ without a global fixed point (as it is non-abelian). As such a fibre $F$ is normal it follows that its minimal subtree is the whole tree, and the action is then cocompact as the vertex and and edge stabilizers in $F$ are finitely generated. 
\end{proof}
 
\subsection{Groups of type $\VF$}\label{sec:vf}

A group is \emph{of type $\mathrm{F}$} if it admits a finite classifying space. Note that if $G$ is of type $\F$, then any finite-index subgroup of $G$ is of type $\F$. A group is \emph{of type $\mathrm{VF}$} if it has a finite-index subgroup of type $\F$.

\bigskip

\begin{prop}[{\cite[Theorem~7.3.4]{Geo}}]
\label{prop:actiononcontractiblespacetypeVF}Let $G$ be a group acting on a CW-complex $X$. Suppose that the action is cocompact and that every element of $G$ preserving a cell fixes it elementwise. If every cell stabiliser of $X$ is of type $\F$, then $G$ is of type $\F$. 
\end{prop}

If we have a short exact sequence \[1\to K\to G\to Q\to 1,\] such that $K$ and $Q$ are of type $\VF$, then it is not true in general that $G$ is of type $\VF$ (see \cite{LN}). However, it follows from the above proposition that being of type $\mathrm{F}$ is stable under extensions, and therefore we have the following result, which was recorded in e.g.\ \cite[Corollary~3.2]{GLMcCool}.

\bigskip

\begin{prop}
\label{prop:type_VF_ses}Suppose that \[ 1 \to K \to G \to Q \to 1\] is a short exact sequence of groups such that \begin{itemize}
\item there exists a finite-index subgroup $G^0$ of $G$ such that $K\cap G^0$ is of type $\F$, and
\item $Q$ is of type $\VF$. 
\end{itemize}
Then $G$ is of type $\VF$.
\end{prop}

\bigskip

\begin{prop}[Persistence of VF with virtually free kernels] \label{prop:virt_free_kernel_VF}Suppose that \[ 1 \to K \to G \to Q \to 1\] is a short exact sequence of groups such that $K$ is finitely generated and virtually free and $Q$ is of type $\VF$. If $G$ is virtually torsion free, then $G$ is also of type $\VF$.
\end{prop}

\begin{proof}
One can pick a finite-index subgroup $G^0\subseteq G$ which is torsion free and whose image in $Q$ is of type $\F$. Then $G^0\cap K$ is a finitely generated free group (so is of type $\F$), and therefore $G^0$ is of type $\F$.
\end{proof}

\subsection{Virtually polycyclic groups}

We will often use the fact that the class of virtually polycyclic groups is stable under extensions. We will also use the following lemma.

\bigskip

\begin{lemma}
\label{lemma:vpc}Let $G$ be a group that fits in the following short exact sequences:
\[ 1\to K_1\to G \stackrel{\pi_1}{\to} Q_1\to 1;\]
\[ 1\to K_2\to G \stackrel{\pi_2}{\to} Q_2\to 1.\]
If $Q_1$ is virtually polycyclic and $\pi_2(K_1)$ is virtually polycyclic, then $Q_2$ is virtually polycyclic.
\end{lemma}

\begin{proof}
The map $\pi_1:G\to Q_1$ factors through a surjective homomorphism \[\bar{\pi}_1:G/(K_1\cap K_2)\to Q_1,\] with kernel $K_1/(K_1\cap K_2)$. The image $Q_1$ of $\bar\pi_1$ is virtually polycyclic by assumption, and its kernel is isomorphic to $\pi_2(K_1)$, whence virtually polycyclic by assumption. Since being virtually polycyclic is stable under group extensions, it follows that $G/(K_1\cap K_2)$ is virtually polycyclic. Therefore so is its quotient $G/K_2\simeq Q_2$. 
\end{proof}

\section{Centralisers and automorphisms of free-by-cyclic groups}\label{sec:autos}

In this section we explain how Theorem~\ref{theo:centralisers-out} (regarding centralisers in $\Out(F)$ and $\Aut(F)$) follows from Theorem~\ref{theo:mccool} (regarding McCool subgroups of $\Out(M,F)$).

 Here, as in the previous section, $M$ is a free-by-cyclic group with fibre $F$ and monodromy $\varphi\in\Aut(F)$, with outer class $\Phi\in\Out(F)$, and $\chi:M\to\mathbb{Z}$ is the character with kernel $F$. And $s\in M$ is an element such that $\chi(s)=1$, so that $M=\langle F,s\rangle$, and the conjugation action of $s$ on $F$ induces the automorphism $\varphi$, i.e.\ $sgs^{-1}=\varphi(g)$ for every $g\in F$.

\subsection{Centralisers in \texorpdfstring{$\Out(F)$}{Out(F)}}\label{sec:out}

We reinterpret $C_{\Out(F)}(\Phi)$ in terms of the automorphism group of the free-by-cyclic group $M=F\rtimes_{\varphi}\ZZ$. 

The restriction to the fibre defines a homomorphism $\Aut^+(M,F)\to\Aut(F)$, which induces a homomorphism \[\kappa:\Aut^+(M,F)\to\Out(F).\] 
 We observe that the image of $\kappa$ is contained in $C_{\Out(F)}(\Phi)$. Indeed, let $\alpha\in\Aut^+(M,F)$, and write $\alpha(s)=us$ for some $u\in F$. Then for every $g\in F$, we have
\begin{equation}\label{eq:commutation}
\alpha|_{F}(\varphi(g))=\alpha(sgs^{-1})=us\alpha|_{F}(g)s^{-1}u^{-1}=\ad_u\circ\varphi(\alpha|_{F}(g)),
\end{equation}
showing that $\varphi$ and $\alpha|_{F}$ commute up to an inner automorphism of $F$.

After taking the quotient by the normal subgroup consisting of inner automorphisms by elements of the fibre, this descends to a homomorphism  \[\overline\kappa:\Aut^+(M,F)/\Inn(F)\to C_{\Out(F)}(\Phi).\]

\begin{prop}
\label{prop:centraliser-out}Suppose that $F$ is non-abelian. The homomorphism $\overline\kappa:\Aut^+(M,F)/\Inn(F)\to C_{\Out(F)}(\Phi)$ is an isomorphism. 
\end{prop}

\begin{proof}
As $F$ is non-abelian, applying Lemma~\ref{lemma:non-abelian-fibre} (with $H=M$) shows that $\overline\kappa$ is injective. To show surjectivity, it suffices to construct an inverse $\theta$ for $\overline\kappa$. Let $\Psi\in C_{\Out(F)}(\Phi)$. Let $\psi\in\Aut(F)$ be a representative of $\Psi$, and $u\in F$, such that $\psi\varphi=\ad_u\varphi\psi$. We let $\theta(\Psi)$ be the automorphism of $M$ given by $\theta(\Psi)(g)=\psi(g)$ for $g\in F$, and $\theta(\Psi)(s)=us$. One checks that this defines an element of $\Aut^+(M,F)$, and that changing the representative of $\Psi$ only changes this element by an inner automorphism from $F$. The map $\Psi\mapsto \theta(\Psi)$ is a group homomorphism. One can see from the definitions that $\overline\kappa \circ \theta$ is the identity, so $\overline\kappa$ is also surjective and therefore an isomorphism. 
\end{proof}

\begin{proof}[Proof of Theorem~\ref{theo:centralisers-out}(\ref{main-1}) using Theorem~\ref{theo:mccool}]
By Proposition~\ref{prop:centraliser-out}, our aim consists in proving that $\Aut^+(M,F)/\Inn(F)$ is of type VF. The natural homomorphism $\Aut(M)\to\Out(M)$ induces a surjective homomorphism $\Aut^+(M,F)/\Inn(F)\to\Out^+(M,F)$, so we have a short exact sequence \[1\to K\to \Aut^+(M,F)/\Inn(F)\to\Out^+(M,F)\to 1\] whose kernel $K$ is cyclic (possibly finite), generated by the inner automorphism of $M$ given by the conjugation by $s$. Theorem~\ref{theo:mccool}, applied with $\calh=\emptyset$, ensures that $\Out^+(M,F)$ is of type $\VF$. Since $\Aut^+(M,F)/\Inn(F)$ is isomorphic to a subgroup of $\Out(F)$ (Proposition~\ref{prop:centraliser-out}), it is virtually torsion free \cite{BT}. We can therefore apply Proposition~\ref{prop:virt_free_kernel_VF} and deduce that $\Aut^+(M,F)/\Inn(F)$ is of type VF, as desired. 
\end{proof}

\subsection{Centralisers in \texorpdfstring{$\Aut(F)$}{Aut(F)}} 

We make the following observation.

\begin{prop}
\label{prop:centraliser-aut}The restriction homomorphism $\Aut(M,F)\to\Aut(F)$ induces an isomorphism between $\Stab_{\Aut(M,F)}(s)$ and $C_{\Aut(F)}(\varphi)$.
\end{prop}

\begin{proof}
Given $\alpha\in\Stab_{\Aut(M,F)}(s)$, Equation~\eqref{eq:commutation} (with $u=1_F$) shows that $\alpha|_{F}$ commutes with $\varphi$. So the restriction homomorphism sends $\Stab_{\Aut(M,F)}(s)$ inside $C_{\Aut(F)}(\phi)$.

An inverse homomorphism $\theta:C_{\Aut(F)}(\varphi)\to\Stab_{\Aut(M,F)}(s)$ is defined as follows: for $\psi\in C_{\Aut(F)}(\varphi)$, we let $\theta(\psi)$ be the automorphism of $M$ given by $\alpha(g)=\psi(g)$ if $g\in F$, and $\alpha(s)=s$.
\end{proof}

\begin{proof}[Proof of Theorem~\ref{theo:centralisers-out}(\ref{main-2}) using Theorem~\ref{theo:mccool}]
The quotient map $\Aut(M)\to\Out(M)$ induces a surjective homomorphism $\Stab_{\Aut(M,F)}(s)\to\Out(M,F,[\langle s\rangle]^{(t)})$, whose kernel is equal to $\Inn(M) \cap \Stab_{\Aut(M,F)}(s)$. By Theorem~\ref{theo:mccool} applied with $\calh=\{[\langle s\rangle]\}$, the image of this homomorphism is of type VF. There are now two cases to consider.

\textbf{$\drsh$ Case 1.} If $Z(M)=\{1\}$, then
\[\Inn(M) \cap \Stab_{\Aut(M,F)}(s)\simeq C_M(s) = \Fix(\phi) \times \langle s \rangle\] is of type F since the free group $\Fix(\phi)$ is finitely generated~\cite{Ger, Cooper}. It follows that $\Stab_{\Aut(M,F)}(s)$ (which is isomorphic to $C_{\Aut(F)}(\varphi)$ by Proposition~\ref{prop:centraliser-aut}) is of type VF.

\textbf{$\drsh$ Case 2.} If $Z(M)$ is infinite, then $\Inn(M)\cap \Stab_{\Aut(M,F)}(s) \simeq C_M(s)/Z(M)$ is virtually free by Proposition~\ref{prop:periodic_implies_centre}, and as in the previous case it is finitely generated. Since $\Stab_{\Aut(M,F)}(s)$ is isomorphic to a subgroup of $\Aut(F)$, it is virtually torsion free \cite{BT}. As $\Stab_{\Aut(M,F)}(s)$ fits into a short exact sequence whose image is of type VF and whose kernel is a finitely generated virtually free group,  $\Stab_{\Aut(M,F)}(s)\simeq C_{\Aut(F)}(\varphi)$ is of type $\VF$ by Proposition~\ref{prop:virt_free_kernel_VF}.
\end{proof}

\section{Finite-order monodromy}\label{sec:finite-order}

In this section, we prove our main theorems in the case where the monodromy of the free-by-cyclic group $M$ has finite order. This will serve as a base case in our arguments in later sections. The main goal of this section is the following theorem:

\begin{thm}\label{th:section-4} 
\label{thm:Out(M)-type-VF-finite-order}Let $M=F\rtimes_\varphi\mathbb{Z}$ be a free-by-cyclic group with finite-order monodromy. There exists a finite-index subgroup $\Out^\per(M)$ of $\Out(M)$ such that, for every finite set $\mathcal{H}$ of conjugacy classes of subgroups of $M$, the following hold:
\begin{enumerate}
\item\label{finite-order-1} the group $\Out^\per(M) \cap \Out(M,\mathcal{H}^{(t)})$ is of type $\mathrm{F}$; 
\item\label{finite-order-2} the group $\Out^\per(M) \cap \Out(M,F) \cap \Out(M,\mathcal{H}^{(t)})$ is of type $\mathrm{F}$. 
\end{enumerate}
\end{thm}

This provides a uniform version of Theorem~\ref{theo:mccool} in the setting of finite-order monodromy, in the sense that $\Out^\per(M)$ is defined independently of $\calh$ (although it does depend on a choice of fibre $F$). 

The strategy is to reduce the result to a theorem of Bestvina--Feighn--Handel \cite{BFH2023} on finiteness properties of subgroups of $\Out(F)$ that follows from earlier work of Krsti\'c--Vogtmann \cite{KV93}. The result we require is given by Theorem~\ref{thm:BFH} below. This theorem is phrased in terms of McCool subgroups of centralisers of finite-order elements in $\Out(F)$, and the rest of Section~\ref{sec:finite-order} consists of the work needed to rephrase it in terms of McCool subgroups of $\Out(M)$.  In Section~\ref{sec:out_per_fibre-preserving_def} we define a subgroup $\Out^\per(M,F)$ of $\Out(M,F)$ which is suitable for the fibre-preserving situation and describe some of its salient properties ($\Out^\per(M,F)$ is chosen in order to control periodic behaviour of automorphisms). This allows us to prove that $\Out^\per(M,F)$, as well as its intersections with McCool subgroups, are of type $\mathrm{F}$ in Section~\ref{sec:finite_order_fibre_preserving}. In the final part of this section we define a subgroup $\Out^\per(M)$ of $\Out(M)$ which unless $M\simeq\mathbb{Z}^2$ satisfies $\Out^\per(M,F)=\Out^\per(M)\cap\Out(M,F)$. We then follow work of Levitt \cite{LevittGBS} and Andrew--Martino \cite{AMFbyZ} to complete the proof of Theorem~\ref{th:section-4} by a reduction to the fibre-preserving case. 

The definitions of $\Out^\per(M,F)$ and $\Out^\per(M)$ are a little complicated but explicit and can be found in Definition~\ref{dfn:out_0_finite_monodromy} and Definition~\ref{def:out_per_absolute} respectively.

\subsection{From Krsti\'c--Vogtmann and Bestvina--Feighn--Handel} \label{sec:BFH_type_F}

As mentioned above, we will reduce the proof of Theorem~\ref{thm:Mccoolgroupsfiniteordermonodromies} to the following theorem about McCool subgroups of centralisers of finite-order elements of $\Out(F)$, due to Bestvina--Feighn--Handel.

\begin{thm}[{Bestvina--Feighn--Handel \cite[Theorem 2.4(3)]{BFH2023}}]
\label{thm:BFH}Suppose $\Phi$ is a finite-order element of $\Out(F)$ and $\calp$ is a finite set of conjugacy classes of subgroups of $F$. Then every torsion-free, finite-index subgroup of $C_{\Out(F)}(\Phi) \cap \Out(F,\calp^{(t)})$ is of type $\F$.
\end{thm}

As the precise statement given above does not appear in \cite{BFH2023} (it instead asserts there exists \emph{some} finite-index subgroup of type $\F$), we give an outline of the proof below. As is noted in \cite{BFH2023}, the necessary ingredients can also be extracted from \cite{KV93}, although the precise statement also does not appear there.

\begin{proof}[Proof Outline] 
First, by~\cite[Theorem~1.5]{GLMcCool}, we can assume that every subgroup $H \subseteq F$ whose conjugacy class is in $\calp$ is finitely generated. 

By e.g.\ \cite[Lemma~2.5]{BFH2023}, we can find a finite list  \[ K=([\langle g_1 \rangle], [\langle g_2 \rangle],\ldots, [\langle g_n \rangle]) \] of conjugacy classes of cyclic subgroups of $F$ such that $\Out(F,\calp^{(t)})=\Out(F,K^{(t)})$, in particular $\Out(F,\calp^{(t)})$ is a finite-index subgroup of $\Out(F,K)$. Let $L_\Phi$ be the subspace of the spine of Outer space consisting of $\Phi$-invariant graphs with no separating edges. Consider the function: \begin{align*} f \colon L_\Phi &\to \mathbb{N}^n \\ \sigma &\mapsto (\|g_1\|_\sigma, \ldots, \|g_n\|_\sigma) \end{align*} where $\sigma \in L_\Phi$ and $\|g_i\|_\sigma$ is the combinatorial length of the reduced loop corresponding to $g_i$ in the graph $\sigma$.  The lexicographic order on $\mathbb{N}^n$ induces a partial order on $L_\Phi$ via $f$. A graph $\sigma \in L_\Phi$ is \emph{reduced} if $\sigma$ contains no non-trivial $\Phi$-invariant forest. Let $\mathcal{R}_{\Phi,K}$ be the set of reduced graphs in $L_\Phi$ which are minimal in the $f$-ordering, and let 
$L_{(\Phi,K)} \subset L_\Phi$ be the union of the stars of the elements of $\mathcal{R}_{\Phi,K}$.

The arguments in Section~7 of \cite{BFH2023} show that $L_{(\Phi,K)}$ is contractible: in the spirit of \cite{CV86,KV93} one extends $f$ to construct a combinatorial Morse function on $L_\Phi$ (precisely, this is the function $\|\sigma\|_G$ on $\mathbb{Z}^\mathcal{K}$ given on \cite[p.82]{BFH2023}). One then shows that descending links of this Morse function are either empty or contractible \cite[Proposition~7.5]{BFH2023}. This shows that $L_{(\Phi,K)}$ is a union of contractible components. Contractibility of $L_{(\Phi,K)}$ then follows the same argument as the proof of \cite[Corollary~6.4]{BFH2023}: one shows that if $L_{(\Phi,K)}$ had multiple components, the same would be true of $L_\Phi$, however $L_\Phi$ is contractible by \cite{KV93} (one can use the same argument to show that in fact descending links of the Morse function are all non-empty and hence contractible). 

A point in the quotient of $\mathcal{R}_{\Phi,K}$ by $C_{\Out(F)}(\Phi)\cap\Out(F,K)$ is determined by a core graph $\sigma$ with fundamental group $F$, a (reduced) action of $\mathbb{Z}/m\mathbb{Z}$ on $\sigma$ (where $m$ is the order of $\Phi$), and a finite collection of $n$ loops in $\sigma$ which are of a given length (determined by the value of $f$ on $\mathcal{R}_{\Phi,K}$). As this is a finite amount of data, the action of  $C_{\Out(F)}(\Phi) \cap \Out(F,K)$ on $L_{(\Phi,K)}$ is cocompact. It follows that any torsion-free, finite-index subgroup of $C_{\Out(F)}(\Phi) \cap \Out(F,\calp^{(t)})$ acts freely and cocompactly on $L_{(\Phi,K)}$, so is of type $\F$.
\end{proof}

\subsection{The finite-index subgroup $\Out^\per(M,F)$} \label{sec:out_per_fibre-preserving_def}

Let $\IA \subseteq \Out(F)$ be the kernel of the homomorphism 
\[\Out(F) \to \mathrm{GL}(\ZZ/3\ZZ)\] defined as the action of $\Out(F)$ on $H_1(F,\ZZ/3\ZZ)$. The group $\IA$ is torsion free \cite{BT}. Another important property of this subgroup, due to Handel--Mosher \cite[Theorem~II.4.1]{HM}, is that if $\Psi\in\IA$ has a power $\Psi^k$ (with $k\neq 0$) that preserves the conjugacy class of an element $g\in F$, then $\Psi$ itself preserves the conjugacy class of $g$. Similarly, conjugacy classes of free factors that are periodic under $\Psi$ are invariant under $\Psi$, although we will not need this fact in the work below. In this section, we define a finite-index subgroup $\Out^\per(M,F)$ of $\Out(M,F)$ which is similarly well behaved. 

\begin{defn}[The groups $\Aut^{\mathrm{IA}}(M,F)$ and $\Out^{\mathrm{IA}}(M,F)$]
\label{de:outia}We define:
\begin{itemize}
\item $\Aut^{\mathrm{IA}}(M,F)$ as the preimage of $\IA$ under the restriction map \[\kappa:\Aut^+(M,F)\to\Out(F).\]
\item $\Out^{\mathrm{IA}}(M,F)$ to be its image in $\Out(M,F)$. 
\end{itemize}
\end{defn} 

Note that $\Aut^{\mathrm{IA}}(M,F)$ is a finite-index subgroup of $\Aut^+(M,F)$ that contains $\Inn(F)$. Given a set $\calh$ of conjugacy classes of subgroups of $M$, we also let 
\[\Aut^{\mathrm{IA}}(M,F,\calh^{(t)}):=\Aut^{\mathrm{IA}}(M,F)\cap\Aut(M,\calh^{(t)});\]
\[\Out^{\mathrm{IA}}(M,F,\calh^{(t)}):=\Out^{\mathrm{IA}}(M,F)\cap\Out(M,\calh^{(t)}).\]  The group $\Out^{\mathrm{IA}}(M,F)$ is a finite-index subgroup of $\Out(M,F)$, and it has the following useful properties. 

 \begin{lemma}
 \label{lemma:ia}Let $M$ be a free-by-cyclic group with finite-order monodromy and non-abelian fibre. For every $\Psi\in\Out^{\mathrm{IA}}(M,F)$, the following hold. 
 \begin{enumerate}
     \item\label{ia-1} Any two representatives of $\Psi$ in $\Aut^{\mathrm{IA}}(M,F)$ differ by a conjugation by an element of $F$.
     \item\label{ia-2} For any representative  $\psi$ of $\Psi$ in $\Aut^{\mathrm{IA}}(M,F)$ and any $g\in F$, if the $M$-conjugacy class of $g$ is periodic under $\Psi$, then the $F$-conjugacy class of $g$ is fixed under $\psi$. 
     \item\label{ia-3} If $\Psi$ has a representative in $\Aut(M)$ acting as the identity on a subgroup $P \subset F$, then $\Psi$ has a representative in $\Aut^{\mathrm{IA}}(M,F)$ acting as the identity on $P$. 
     \item\label{ia-4} Let $\calh$ be a set of $M$-conjugacy classes of subgroups of $F$, and let $\calh|_{F}$ be the set of their $F$-conjugacy classes. The image of $\Aut^{\mathrm{IA}}(M,F,\calh^{(t)})$  under the restriction map $\kappa$ is a torsion-free finite-index subgroup of $C_{\Out(F)}(\Phi)\cap\Out(F,\calh|_{F}^{(t)})$.  
 \end{enumerate}
 \end{lemma}

 \begin{proof}
As ususal, we denote by $\varphi\in\Aut(F)$ the chosen monodromy of the free-by-cyclic group $M$, and by $\Phi\in\Out(F)$ its outer class. We also recall that $s\in M$ is an element such that $\langle F,s\rangle=M$, and the conjugation action of $s$ on $M$ induces the automorphism $\varphi$. We finally let $m>0$ be the order of $\Phi$ in $\Out(F)$. 

\smallskip
 
\textbf{Assertion 1.} Let $\psi_1,\psi_2\in\Aut^{\mathrm{IA}}(M,F)$ be two lifts of $\Psi$. Then there exist $g\in F$ and $k\in \{0,\dots,m-1\}$ such that $\psi_2=\ad_g\circ\ad_{s^k}\circ\psi_1$. Since $\psi_1,\psi_2$ and $\ad_g$ belong to $\Aut^{\mathrm{IA}}(M,F)$, it follows that $\ad_{s^k}\in\Aut^{\mathrm{IA}}(M,F)$. But this is only possible if $k=0$: otherwise $\kappa(\ad_{s^k})=\Phi^k$ has finite order and is non-trivial in $\Out(F)$, and therefore does not belong to the torsion-free subgroup $\IA$. Therefore $\psi_2=\ad_g\circ\psi_1$, and the first assertion follows.
 
\smallskip

\textbf{Assertion 2.}  We denote by $[\cdot]_M$ the conjugacy class in $M$, and by $[\cdot]_F$ the conjugacy class in $F$. Suppose that the $\psi$-orbit of $[g]_M$ is finite with cardinality $k \in \NN$. Since $\varphi^m$ is inner, it preserves $[g]_F$. Let $\ell\in\{1,\dots,m\}$ be the smallest positive integer such that $[\varphi^\ell(g)]_F=[g]_F$.
Then \[[g]_M=\coprod_{i=0}^{\ell-1}~[\varphi^i(g)]_F\] gives the partition of the $M$-conjugacy class of $g$ into $F$-conjugacy classes. Since $\psi^k$ preserves the fibre $F$ and preserves $[g]_M$, it preserves the above partition of $[g]_M$, possibly permuting the subsets of this partition (it has a power that preserves each subset of the partition). Thus $[g]_F$ is periodic under $\psi|_{F}$. Since the outer class of $\psi|_{F}$ belongs to $\IA$, it follows from \cite[Theorem~II.4.1]{HM} that $[g]_F$ is in fact fixed by $\psi|_{F}$. 
 
\smallskip

\textbf{Assertion 3.} 
Assume that $\Psi$ has a representative acting as the identity on $P$, and let $\psi \in \Aut^{\mathrm{IA}}(M,F)$ be some representative of $\Psi$. Then $\psi$ preserves the $M$-conjugacy class of every element of $P$. Assertion~2 then ensures that for every $g\in P$, the image $\psi(g)$ is $F$-conjugate to $g$. It thus follows from e.g.\ \cite[Lemma~6.2]{GLMcCool} that $\psi|_{F}$ acts on $P$ as a global conjugation by an element of $F$.
Precomposing $\psi$ by this inner automorphism gives the desired representative of $\Psi$ in $\Aut^{\mathrm{IA}}(M,F)$.

\smallskip

\textbf{Assertion 4.} Assertion~3 ensures that for every $\psi \in \Aut^{\mathrm{IA}}(M,F,\calh^{(t)})$, we have $\kappa(\psi)\in\Out(F,\calh|_{F}^{(t)})$. Since $\kappa$ takes its values in $C_{\Out(F)}(\Phi)$ (see Section~\ref{sec:out}), we have \[\kappa(\Aut^{\mathrm{IA}}(M,F,\calh^{(t)}))\subseteq C_{\Out(F)}(\Phi)\cap \Out(F,\calh|_{F}^{(t)})\subseteq\kappa(\Aut^+(M,F,\calh^{(t)}))\] (for the rightmost inclusion, any $\Psi\in C_{\Out(F)}(\Phi)$ lifts to $\Aut^+(M,F)$ by Proposition~\ref{prop:centraliser-out}, and if $\Psi$ acts trivially on $\calh|_{F}$, then the lift also acts trivially on $\calh$). It follows that these inclusions are finite index. The image is torsion free as $\IA$ is torsion free.\qedhere
\end{proof}

\begin{lemma}
\label{lem:finitelymanyroot}Let $M$ be a free-by-cyclic group with finite-order monodromy. For every $g\in Z(M)$, the set of roots of $g$ is contained in finitely many $F$-conjugacy classes.
\end{lemma}

\begin{proof}
We will assume that $M$ is not isomorphic to $\mathbb{Z}^2$, as otherwise the statement is clear. Let $g\in Z(M)$. If $g=1_M$, then $g$ has no proper roots as $M$ is torsion free. So we may assume that $g$ is non-trivial. Let $R(g)$ be the set of all roots of $g$, and for every $k\in\mathbb{N}$, let $R_k(g)$ be the set of all $k$-th roots of $g$. We have $R(g)=\bigcup_{k\in\mathbb{N}}R_k(g)$, so the lemma follows from the two facts below.

\smallskip

\textbf{Fact 1.} For all sufficiently large $k\in\mathbb{N}$, we have $R_k(g)=\emptyset$. 

Indeed, since $g$ is non-trivial and $M$ is not isomorphic to $\mathbb{Z}^2$, Proposition~\ref{prop:periodic_implies_centre} ensures that $\chi(g)\neq 0$ (where $\chi:M\to\mathbb{Z}$ is the character with kernel $F$). If $h\in R_k(g)$, then $\chi(h)=\frac{1}{k}\chi(g)$, so $k\le|\chi(g)|$. 

\smallskip

\textbf{Fact 2.} For a given $k\in\mathbb{N}$, the set $R_k(g)$ is contained in finitely many $F$-conjugacy classes.

Indeed, since $g\in Z(M)$, the image $\overline{R_k(g)}$ of $R_k(g)$ in the virtually free group $M/Z(M)$ consists of finite-order elements. Since virtually free groups have only finitely many conjugacy classes of finite-order elements, it follows that $\overline{R_k(g)}$ is contained in finitely many $(M/Z(M))$-conjugacy classes. Since $F$ embeds in $M/Z(M)$ under the quotient map, and its image has finite index, we deduce that $\overline{R_k(g)}$ is contained in finitely many $F$-conjugacy classes. It is therefore enough to show that if $g_1,g_2\in R_k(g)$ have $F$-conjugate images in $M/Z(M)$, then $g_1$ and $g_2$ are $F$-conjugate in $M$. In this case there exist $h \in F$ and $z \in Z(M)$ such that $g_1=hg_2h^{-1}z$. Raising to the $k$-th power gives $g=hgh^{-1}z^k$, and applying the character $\chi$ shows that $z=1_M$, so $g_1$ and $g_2$ are $F$-conjugate, as desired. 
\end{proof}

The group $\Aut(M,F)$ acts on the set of all $F$-conjugacy classes of elements of $M$. We consider the restriction of this action to $\Aut^{\mathrm{IA}}(M,F)$. 

Let $z$ be a generator of $Z(M)$. By Lemma~\ref{lem:finitelymanyroot}, the element $z$ has only finitely many $F$-conjugacy classes of roots in $M$, and as $z$ is fixed by any automorphism in $\Aut^+(M,F)$, the group $\Aut^{\mathrm{IA}}(M,F)$ permutes these $F$-conjugacy classes.

\begin{defn}[The finite-index subgroup $\Out^\per(M,F)$] \label{dfn:out_0_finite_monodromy}Let $M$ be a free-by-cyclic group with finite-order monodromy.
Let $z \in Z(M)$ be a generator of $Z(M)$. 

We let $\Out^\per(M,F)$ be the finite-index subgroup of $\Out^{\mathrm{IA}}(M,F)$ defined in the following way:
\begin{itemize}
\item If $F$ is non-abelian, then $\Out^\per(M,F)$
consists of elements $\Psi$ such that some (equivalently, any) representative of $\Psi$ in $\Aut^{\mathrm{IA}}(M,F)$ preserves the $F$-conjugacy class of every root of $z$. 

\item If $M$ is cyclic or the Klein bottle group, then $\Out^\per(M,F)=\{1\}$. 
\item If $M \simeq \mathbb{Z}^2$, then $\Out^\per(M,F)$ is the index-two subgroup of $\Out^+(M,F)$ acting trivially on $F \simeq \mathbb{Z}$ (in this case $\Out^\per(M,F)\simeq \mathbb{Z})$.
\end{itemize}
We let $\Aut^\per(M,F)\subseteq\Aut^{\mathrm{IA}}(M,F)$ be the full preimage of $\Out^\per(M,F)$ in $\Aut^{\mathrm{IA}}(M,F)$.
\end{defn}

\begin{rk}\label{rem:lifts_differ_by_els_of_F}
The `equivalently any' part of the first bullet point in the above definition is justified by the fact that any two lifts of $\Psi \in \Out^{\mathrm{IA}}(M,F)$ to $\Aut^{\mathrm{IA}}(M,F)$ differ by a conjugation by an element of $F$ (Lemma~\ref{lemma:ia}(1)).
\end{rk}

Given a set $\calh$ of conjugacy classes of subgroups of $M$, we let 
\[\Out^\per(M,F,\calh):=\Out^\per(M,F)\cap\Out(M,F,\calh),\]
\[\Out^\per(M,F,\calh^{(t)}):=\Out^\per(M,F)\cap\Out(M,F,\calh^{(t)}),\] \[\Aut^\per(M,F,\calh^{(t)}):=\Aut^\per(M,F)\cap\Aut(M,F,\calh^{(t)}).\]

\subsection{The proof in the fibre-preserving case} \label{sec:finite_order_fibre_preserving}

Having defined and proved key properties of $\Out^\per(M,F)$ in the previous section, in this section we show that this group has a finite classifying space. Our current aim is:

\begin{thm}
\label{thm:Mccoolgroupsfiniteordermonodromies}Let $M$ be a free-by-cyclic group with finite-order monodromy. Then for every finite set $\calh$ of conjugacy classes of subgroups of $M$, the group \[\Out^\per(M,F,\calh^{(t)})=\Out^\per(M,F)\cap\Out(M,\calh^{(t)})\] is of type $\mathrm{F}$.
\end{thm}

We will prove the above by reducing it to Theorem~\ref{thm:BFH}. In the sequel, when writing $[\cdot]$, we always mean the conjugacy class in $M$. A crucial observation is given by the following proposition.

\begin{prop} 
\label{prop:trivial_action_on_abelians_finite_case}Let $M$ be a free-by-cyclic group with finite-order monodromy and non-abelian fibre and let $H$ be a virtually abelian subgroup of $M$.
\begin{enumerate}
\item\label{item:1} If $\Psi\in\Out^\per(M,F)$ preserves $[H]$, then $\Psi$ has a representative in $\Aut^{\mathrm{IA}}(M,F)$ acting as the identity on $H$. In particular $\Out^\per(M,F,[H]) =\Out^\per(M,F,[H]^{(t)})$. 
\item\label{item:2} If $H \cap F \neq \{1\}$, then $\Out^\per(M,F,[H]^{(t)})=\Out^\per(M,F,[H\cap F]^{(t)})$. 
\end{enumerate}
\end{prop}

\begin{proof}
We will first treat the case where $Z(M)\subseteq H$, and explain at the end of this proof how to reduce to this case.

\medskip

\noindent\textbf{Case 1.} We have $Z(M)\subseteq H$.

\smallskip

\textbf{(i).} If $H \cap F = \{1\}$, then $H\simeq \mathbb{Z}$ and $\chi|_{H}$ is injective. Therefore $H$ is generated by a root of a generator of $Z(M)$. The $F$-conjugacy classes of such roots are preserved by representatives of elements in $\Out^\per(M,F)$ that belong to $\Aut^{\rm{IA}}(M,F)$, ensuring that the first part of the proposition holds in this case (the second part is void in this case).

\textbf{(ii).} We now assume that $H \cap F \neq \{1\}$. We will prove the following chain of inclusions:
\begin{equation*}\label{eq:inclusions}\tag{$\star$}
    \Out^\per(M,F,[H])  \subset \Out^\per(M,F,[H\cap F]^{(t)}) \subset  \Out^\per(M,F,[N_M(H\cap F)]^{(t)}),
\end{equation*}
 and in fact we will see that every element of $\Out^\per(M,F,[H])$ has a representative in $\Aut^\mathrm{IA}(M,F)$ acting as the identity on $N_M(H\cap F)$. 
As $H \subseteq N_M(H\cap F)$, the inclusions \[\Out^\per(M,F,[N_M(H\cap F)]^{(t)}) \subset \Out^\per(M,F,[H]^{(t)}) \subset \Out^\per(M,F,[H])\] are immediate, so it is enough to focus on the first chain. 

\medskip

\textbf{Leftmost inclusion in \eqref{eq:inclusions}.}  Let $\Psi \in \Out^\per(M,F,[H])$. We aim to find a representative $\psi$ of $\Psi$ in $\Aut^{\mathrm{IA}}(M,F)$ acting as the identity on $H \cap F$. So let $\psi_1$ be any representative of $\Psi$ in  $\Aut^{\mathrm{IA}}(M,F)$. Since the $M$-conjugacy class of any generator of the cyclic group $H\cap F$ has a finite orbit under $\Psi$, it follows from Lemma~\ref{lemma:ia}(\ref{ia-2}) that $\psi_1$ coincides on $H\cap F$ with the conjugation by an element $h\in F$. Then $\psi=\ad_{h^{-1}}\circ\psi_1$ is the desired representative of $\Psi$.

\medskip

\textbf{Rightmost inclusion in \eqref{eq:inclusions}.} Let $\Psi \in \Out^\per(M,F,[H\cap F]^{(t)})$. Our goal is to find a representative $\psi\in\Aut^{\mathrm{IA}}(M,F)$ acting as the identity on $N_M(H\cap F)$. By our proof of the leftmost inclusion in \eqref{eq:inclusions}, we can choose $\psi' \in \Aut^\per(M,F)$ representing $\Psi$ and acting as the identity on $H\cap F$. As $\psi'$ fixes $H \cap F$, this automorphism also preserves $N_M(H \cap F)$.

\begin{claim}\label{claim:inner}
The automorphism $\psi'$ acts on $N_M(H\cap F)$ as the conjugation by an element $g\in H\cap F$. 
\end{claim}

\begin{proof}[Proof of the claim]
As $H \cap F$ is infinite cyclic, the same is true of $N_M(H \cap F) \cap F$. As $N_M(H \cap F)$ contains $Z(M)$ this group has a finite index subgroup isomorphic to $\mathbb{Z}^2$. 
Since $\psi'$ belongs to $\Aut^+(M,F)$, it acts as the identity on the subgroup $Z(M)$ of $N_M(H\cap F)$. There are two cases to consider. 

\smallskip

\textbf{$\drsh$ 1.} Assume first that $N_M(H\cap F)$ is isomorphic to $\mathbb{Z}^2$. Then every automorphism of $N_M(H\cap F)$ acting trivially on the finite-index subgroup $\langle H\cap F, Z(M)\rangle$ is the identity, showing that $\psi'$ acts as the identity on $N_M(H \cap F)$.

\textbf{$\drsh$ 2.} Assume now that $N_M(H\cap F)$ is isomorphic to the fundamental group of a Klein bottle, generated by two elements $a,t$ with $tat^{-1}=a^{-1}$ and $a\in F$. The generator $z$ of $Z(M)$ is contained in $N_M(H\cap F)$, while the centre of $N_M(H\cap F)$ is generated by $t^2$, showing that $z=t^{2k}$ for some integer $k\neq 0$. It follows from the definition of $\Aut^\per(M,F)$ that $\psi'$ sends $t$ to an $F$-conjugate of $t$. Since $N_M(H\cap F)$ is the fundamental group of a Klein bottle, its outer automorphism group is isomorphic to $\mathbb{Z}/2\mathbb{Z}\times\mathbb{Z}/2\mathbb{Z}$, generated by the outer classes of the following automorphisms: \[\alpha_1: a\mapsto a, t\mapsto t^{-1} \quad ; \quad \alpha_2: a\mapsto a, t\mapsto ta\] (see e.g.\ \cite[Lemma~2.1]{GW}). 

We now check that none of the $(N_M(H\cap F))$-conjugates of the elements $t^{-1}$, $t^{-1}a$ and $ta$ (the respective images of $t$ by $\alpha_1,\alpha_1\alpha_2$, and $\alpha_2$) are $F$-conjugate to $t$ in $M$. For $t^{-1}$ and $t^{-1}a$, this follows from considering the image under the character $\chi:M\to\mathbb{Z}$ (with kernel $F$). The $N_M(H\cap F)$-conjugates of $ta$ are all of the form $(a^kt^p)ta(a^kt^p)^{-1}=ta^{-2k+(-1)^p}$ for some $k,p\in\ZZ$, and we need to check that this is never $F$-conjugate to $t$. This is because if $ta^{-2k+(-1)^p}=gtg^{-1}$ for some $g\in F$, then the automorphism $\kappa(\ad_{t^{-1}})$ would send $g$ to $a^{-2k+(-1)^p}g$, which is impossible since this automorphism has finite order in $\Out(F)$, while $-2k+(-1)^p\neq 0$. 

 The above ensures that $\psi'$ induces an inner automorphism of $N_M(H\cap F)$, i.e.\ it acts on $N_M(H\cap F)$ as the conjugation by an element $h\in N_M(H\cap F)$. Since $\psi'$ was chosen to be the identity on $H\cap F$, we have $\psi'(a)=a$, and therefore $h$ belongs to the centraliser of $a$ in $N_M(H\cap F)$, which is equal to $\langle a,t^2\rangle$. Let $p,q\in\mathbb{Z}$ be such that $h=a^p t^{2q}$. Since $t^2$ is central in $N_M(H\cap F)$, the conjugation by $h$ on $N_M(H\cap F)$ coincides with the conjugation by $a^p$. It follows that $\psi'$ coincides with the conjugation by $a^p$ on $N_M(H\cap F)$. 
\end{proof}

After composing $\psi'$ with the conjugation by $g^{-1}$, we get a representative of $\Psi$ contained in $\Aut^{\mathrm{IA}}(M,F)$ and which acts as the identity on $N_M(H\cap F)$. This concludes our proof in Case~1.

\medskip

\noindent\textbf{Case 2.} The general case.

\smallskip

Let $H'=\langle H,Z(M)\rangle$. Then $H'$ is virtually abelian and is preserved by any automorphism preserving $H$. So the first part of the proposition follows for $H$ because by Case~1 it holds for $H'$. If additionally $H\cap F\neq\{1\}$, then it is an infinite cyclic subgroup of $F$, and therefore it has finite index in $H'\cap F$. Using the uniqueness of roots in $F$, we deduce that $\Out(M,F,[H\cap F]^{(t)})=\Out(M,F,[H'\cap F]^{(t)})$, and the second part of the proposition follows for $H$ because by Case~1 it holds for $H'$. 
\end{proof}

\begin{lem}
\label{lem:replaceperipheralwithfibreperipheral}Let $M$ be a free-by-cyclic group with finite-order monodromy and non-abelian fibre, let $H\subseteq M$ be a subgroup, and let $\bar H=\langle H,Z(M)\rangle\cap F$. Then   \[\Out^\per(M,F,[H]^{(t)})=\Out^\per(M,F,[\bar H]^{(t)}).\]
\end{lem}

\begin{proof}
Notice that every automorphism in $\Aut^+(M,F)$ acts trivially on $Z(M)$, so \[\Out^\per(M,F,[H]^{(t)})=\Out^\per(M,F,[\langle H,Z(M)\rangle]^{(t)}).\] Thus, by replacing $H$ by $\langle H,Z(M)\rangle$, we may assume that $Z(M)\subseteq H$ and $\bar H=H\cap F$. There are three cases to consider. 

\textbf{$\drsh$ 1.} If $H \cap F$ is non-abelian, then Lemma~\ref{lemma:non-abelian-fibre} ensures that \[ \Out^\per(M,F,[H]^{(t)})=\Out^\per(M,F,[H \cap F]^{(t)}). \]

\textbf{$\drsh$ 2.} If $H \cap F$ is a non-trivial cyclic group, then $H$ is virtually abelian and the above equality is given in Proposition~\ref{prop:trivial_action_on_abelians_finite_case}(\ref{item:2}). 

\textbf{$\drsh$ 3.} If $H\cap F=\{1\}$, then $H$ is cyclic and contains $Z(M)$, so is generated by a root of a generator of $Z(M)$. The $F$-conjugacy classes of such roots are preserved by representatives of elements in $\Out^\per(M,F)$ which lie in $\Aut^{\rm{IA}}(M,F)$, therefore every element of $\Out^\per(M,F)$ acts trivially on $H$, and again the conclusion follows. 
\end{proof}

We now have all the tools we need to prove that $\Out^\per(M,F,\mathcal{H}^{(t)})$ is of type $\F$  for every finite set $\calh$ of $M$-conjugacy classes of subgroups of $M$.

\begin{proof}[Proof of Theorem~\ref{thm:Mccoolgroupsfiniteordermonodromies}]
As usual, we let $\chi:M\to\ZZ$ be the character with kernel $F$, and let $s\in M$ be an element such that $M=\langle F,s\rangle$ and conjugation by $s$ restricts to the automorphism $\varphi$ on $F$. 

When $M$ is cyclic or the Klein bottle group,  $\Out^\per(M,F)=\{1\}$. When $M \simeq \ZZ^2$ we have $\Out^\per(M,F)\simeq\mathbb{Z}$, so every subgroup of $\Out^\per(M,F)$ is of type $\F$.

Now suppose that $F$ is non-abelian. 
By Lemma~\ref{lem:replaceperipheralwithfibreperipheral}, when proving that $\Out^\per(M,F,\calh^{(t)})$ is of type $\F$, we may replace every subgroup $H$ from the collection $\mathcal{H}$ with the subgroup $\bar{H}=\langle H, Z(M) \rangle \cap F$. Hence for the remainder of the proof we will assume that $\mathcal{H}$ is a finite set of $M$-conjugacy classes of subgroups of $F$.

The natural projection $\Aut(M)\to\Out(M)$ induces a surjective homomorphism \[\theta:\Aut^\per(M,F,\calh^{(t)})/\Inn(F)\to\Out^\per(M,F,\calh^{(t)}),\] whose kernel is equal to $(\Aut^\per(M,F,\calh^{(t)})\cap \Inn(M))/\Inn(F)$. Since $\Inn(M)\simeq M/Z(M)$ contains $F$ as a finite-index subgroup, the kernel of $\theta$ is finite. It is therefore enough to prove that $\Aut^\per(M,F,\calh^{(t)})/\Inn(F)$ is of type $\F$ (indeed, this means in particular that it is torsion free, so $\theta$ will then be an isomorphism, thereby showing that $\Out^\per(M,F,\calh^{(t)})$ is of type $\F$).

So let us prove that $\Aut^\per(M,F,\calh^{(t)})/\Inn(F)$ is of type $\F$. Using Proposition~\ref{prop:centraliser-out}, we have $\Aut^\per(M,F,\calh^{(t)})/\Inn(F)\simeq\kappa(\Aut^\per(M,F,\calh^{(t)}))$.
Since the monodromy $\varphi$ of $M$ has finite order in $\Out(F)$, the set $\calh|_{F}$ (a set of $F$-conjugacy classes of subgroups of $F$) is finite. Since $\Aut^\per(M,F)\subseteq\Aut^{\mathrm{IA}}(M,F)$, Lemma~\ref{lemma:ia}(\ref{ia-4}) ensures that $\kappa(\Aut^\per(M,F,\calh^{(t)}))$ is a torsion-free finite-index subgroup of $C_{\Out(F)}(\Phi)\cap\Out(F,\calh|_{F}^{(t)})$. Therefore it is of type $\F$ by Theorem~\ref{thm:BFH},  as desired.  
\end{proof}

\subsection{$\Out^\per(M)$ and the complete proof of Theorem~\ref{th:section-4}} \label{sec:finite_order_absolute}

As with the rest of this section, let $M$ be a free-by-cyclic group with finite-order monodromy, with a preferred fibring $\chi:M \to \mathbb{Z}$. The centre $Z(M)$ is infinite, and the quotient group $\overline{M} \coloneq M/Z(M)$ is virtually free (possibly trivial if $M\simeq \mathbb{Z}^2$). The quotient $M \to \overline{M}$ induces a homomorphism  \[\tilde{p}_Z:\Aut(M)\to\Aut(\overline{M}),\] which descends to a homomorphism \[p_Z \colon \Out(M) \to \Out(\overline{M}).\] Assuming that $M$ is not virtually abelian, we will construct finite-index subgroups  $\Out^{\overline{\chi}}(M)$ and $\Out^{\overline{\chi}}(\overline{M})$ of $\Out(M)$ and $\Out(\overline{M})$ respectively fitting into a split short exact sequence \[
1 \to \ZZ^{\rank(M^{\ab}) - 1} \to \Out^{\overline{\chi}}(M) \to \Out^{\overline{\chi}}(\overline{M}) \to 1.\] 
Denoting by $\Aut^{\overline{\chi}}(M)$ and $\Aut^{\overline{\chi}}(\overline{M})$ their preimages in $\Aut(M)$ and $\Aut(\overline{M})$, we will see that every automorphism in $\Aut^{\overline{\chi}}(\overline{M})$ has a unique lift in $\Aut^+(M,F)$ (Lemma~\ref{lemma:unique-lift}), and this induces an isomorphism $\tau_Z:\Out^{\overline{\chi}}(\overline{M})\to\Out^+(M,F)$ that lifts the projection $p_Z$. The composition $R=\tau_Z\circ p_Z$ is a retraction from $\Out^{\overline{\chi}}(M)$ to $\Out^+(M,F)$ with kernel $\mathbb{Z}^{\rank(M^{\ab})-1}$ (Proposition~\ref{prop:retract}). Finally, given a set $\calh$ of conjugacy classes of subgroups, the retraction respects it (Lemma~\ref{lemma:andrewmartino}), allowing us to reduce the proof of  Theorem~\ref{th:section-4} to the fibre-preserving case, already established in Theorem~\ref{thm:Mccoolgroupsfiniteordermonodromies}.

\begin{rk}
In \cite[Section 4]{LevittGBS}, Levitt proves similar results from the viewpoint of automorphisms of generalised Baumslag--Solitar groups with trivial modulus (or more generally those that are unimodular) by constructing a certain fibration, and this approach was extended to incorporate peripheral subgroups in \cite[Section 5.3]{AMFbyZ}. Here we already have a fibration $\chi$, and we use this in our arguments instead, allowing us to work with the original fibre $F$. The cost is passing to a subgroup whose index is in general larger than that provided by Levitt's work. 
\end{rk}

\paragraph*{Studying the quotient from $\Out(M)$ to $\Out(\overline{M})$.}

This is based on the following general proposition.

\begin{prop}
\label{prop:quotient-general}Let $G$ be a group with centre $Z(G)$ and quotient $\overline{G}=G/Z(G)$. The kernel $K$ of the homomorphism \[ \Aut(G) \to \Aut(\overline{G}) \times \Aut(Z(G)) \]
consists of automorphisms of the form $g \mapsto gz(g)$, where $z:G\to Z(G)$ is a homomorphism which is trivial on $Z(G)$. It follows that $K \simeq \Hom(\overline{G},Z(G)).$ 

The subgroup $K \subset \Aut(G)$ maps onto the kernel $\overline{K}$ of the homomorphism \[ \Out(G) \to \Out(\overline{G}) \times \Out(Z(G)), \]  and $K\cap \Inn(G)$ is equal to the kernel of the natural map $\Inn(G) \to \Inn(\overline{G})$. It follows that if $Z(\overline{G})=\{1\}$ (or equivalently $\Inn(G) \simeq \Inn(\overline{G})$), then $K \cap \Inn(G)=\{1\}$ and $K \simeq \overline{K}$. 
\end{prop}

\begin{proof} As in the statement let $K \subseteq \Aut(G)$ be the subgroup of automorphisms acting trivially on both $\overline{G}$ and $Z(G)$. If $\psi \in K$, then as $\psi$ acts trivially on $\overline{G}$ then for every $g \in G$ there exists $z_{\psi}(g) \in Z(G)$ such that \[ \psi(g)=gz_{\psi}(g). \]
One can check that $\psi(gh)=ghz_\psi(g)z_\psi(h)$, so the map $z_\psi \colon G \to Z(G)$ is a homomorphism which is trivial on $Z(G)$ as $\psi|_{Z(G)}=\id_{Z(G)}$. We claim that $\psi \mapsto z_\psi$ is an injective homomorphism from $K$ to $\Hom(G,Z(G))$ with image equal to the subgroup of homomorphisms that are trivial on $Z(G)$. If $z_\psi$ is trivial then clearly $\psi$ is the identity, so $\psi \mapsto z_\psi$ is injective. If $g \in G$ and $\psi_1,\psi_2 \in K$ then \begin{align*} \psi_1\psi_2(g)&=\psi_1(gz_{\psi_2}(g)) &\text{by definition of $z_{\psi_2}$}\\
&=\psi_1(g)z_{\psi_2}(g) &\text{as $\psi_1$ is the identity on $Z(G)$} \\
&=gz_{\psi_1}(g)z_{\psi_2}(g) &\text{by definition of $z_{\psi_1}$},
\end{align*}
which shows that $\psi \mapsto z_\psi$ is a homomorphism. The image of $K$ under this map is the subgroup of $\Hom(G,Z(G))$ which maps $Z(G)$ to 0 as any such homomorphism $f$ determines an automorphism $\psi_f(g):=gf(g)$ in $K$, and we have already noted that if $f=z_\psi$ for $\psi \in K$ then $f$ is trivial on $Z(G)$. It follows that $K \simeq \Hom(\overline{G},Z(G))$ as $\Hom(\overline{G},Z(G))$ is isomorphic to the subgroup of $\Hom(G,Z(G))$ consisting of maps sending $Z(G)$ to 0. 

For the second part, let $\overline{K}$ be the kernel of the natural map \[\Out(G) \to \Out(\overline{G}) \times \Out(Z(G)).\] We claim that $K$ maps onto $\overline{K}$. Indeed let $\psi$ be any representative of $\Psi \in \overline{K}$. Then $\psi$ acts on $\overline{G}$ as $\ad_{\bar g}$ for some $\bar g\in\overline{G}$. Let $g\in G$ be a lift of $\bar g$. Then $\ad_g^{-1}\circ\psi$ is a representative of $\Psi$ in $K$ (for the $Z(G)$ factor, note that if some representative of an outer automorphism acts as the identity on $Z(G)$ then every representative acts as the identity on $Z(G)$). 

An inner automorphism $\ad_h$ belongs to $K$ if and only if it acts trivially on $\overline{G}$, which happens if and only if the image $\bar{h}$ of $h$ belongs to $Z(\overline{G})$. When $Z(\overline{G})=1$ this can only happen if $h \in Z(G)$ and the inner automorphism $\ad_h$ is trivial.
\end{proof}

Note that this proposition is very general: in addition to Levitt's work on generalised Baumslag--Solitar groups it has clear applications to the study of automorphism groups of nilpotent groups, and also appears in the study of automorphism groups of RAAGs with non-trivial centre (see Propositon~4.4 of \cite{CharneyVogtmann}).

\paragraph{The finite-index subgroup $\Out^{\overline{\chi}}(M)$.}

From now on we assume that $M$ has non-abelian fibre. The character $\chi:M \to \ZZ$ (with kernel $F$) induces a map $\overline{\chi}:\overline{M} \to \ZZ/m\ZZ$, where $m$ is the order of the monodromy $\Phi$. Since $F$ intersects $Z(M)$ trivially, it maps injectively to $\overline{M}$, and this image $\overline{F} \simeq F$ is the kernel of $\overline{\chi}$. The situation is illustrated in Figure~\ref{figure:fibrings}.

\begin{figure}

\centering
\begin{tikzcd}
& & 1 \ar[d] & 1\ar[d] \\
& & Z(M) \ar[r,"\simeq"] \ar[d] & m\ZZ \ar[d] \\
1 \ar[r]& F \ar[r] \ar[d,"\simeq"] & M \ar[r, "\chi"] \ar[d] & \ZZ \ar[r] \ar[d] & 1 \\
1 \ar[r]& \overline{F} \ar[r] & \overline{M} \ar[r, "\overline{\chi}"] \ar[d] & \ZZ/m\ZZ \ar[r] \ar[d]& 1 \\
& & 1 & 1 \\
\end{tikzcd}

\caption{The relationship between the fibrings $\chi$ and $\overline{\chi}$ of $M$ and $\overline{M}$ respectively, shown as a commutative diagram. All rows and columns are exact, all squares commute. The bottom exact sequence is not necessarily split: this turns out to depend exactly on if there is a choice of $\varphi$ representing $\Phi$ that is finite order as an automorphism.}\label{figure:fibrings}
\end{figure}  

Let $Z(M)=\langle z \rangle$ with $\chi(z) \neq 0$. Denote by $\Out^+(M)$ the finite-index subgroup of $\Out(M)$ fixing $Z(M)$ elementwise. Let $\Out^{\overline{\chi}}(\overline{M})$ be the subgroup of $\Out(\overline{M})$ that (acting by precomposition) preserves $\overline{\chi}$. This is finite index, since the target of $\overline{\chi}$ is a finite group; an equivalent description is that it preserves $\overline{F}$ and all of its cosets. Let $\Out^{\overline{\chi}}(M)$ be the preimage by $p_Z$ of $\Out^{\overline{\chi}}(\overline{M})$ in $\Out^+(M)$, a finite index subgroup of $\Out(M)$. Observe that the subgroup $\Out^+(M,F)$ is contained in $\Out^{\overline{\chi}}(M)$.

We also let $\Aut^{\overline{\chi}}(\overline{M})$ be the full preimage of $\Out^{\overline{\chi}}(\overline{M})$ in $\Aut(\overline{M})$, and $\Aut^{\overline{\chi}}(M)$ be the full preimage of $\Out^{\overline{\chi}}(M)$ in $\Aut(M)$.

We record the following consequence of the previous proposition.

\bigskip

\begin{lemma}[{cf. Levitt \cite[Lemma 4.2, Theorem 4.4]{LevittGBS}}]
\label{lem:kernel_projection}Let $M$ be a free-by-cyclic group with finite-order monodromy and non-abelian fibre.
The group $$N=\ker(p_Z \colon \Out^{\overline{\chi}}(M) \to \Out^{\overline{\chi}}(\overline{M}))$$ is isomorphic to $\ZZ^{\rank(M^{\ab}) - 1}$, where $\rank(M^{\ab})$ is the rank of the free part of the abelianisation of $M$.
\end{lemma}

\begin{proof}
As $M/Z(M)$ has trivial centre (Lemma~\ref{lemma:trivial-center}), Proposition~\ref{prop:quotient-general} ensures that $N$ is a finite-index subgroup of $\mathrm{Hom}(\overline{M},Z(M))$. Since $Z(M)$ is isomorphic to $\mathbb{Z}$, it follows that $\mathrm{Hom}(\overline{M},Z(M))$ (and hence $N$) is isomorphic to $\mathbb{Z}^{\rank(\overline{M}^\ab)}$.

It remains to prove that $\ZZ^{\rank(\overline{M}^\ab)}\simeq \ZZ^{\rank(M^\ab)-1}$. The natural surjection $M \to \overline{M}$ induces a surjection $M^\ab \to \overline{M}^\ab$ whose kernel is equal to the image of $Z(M)$ in $M^\ab$. As $Z(M)$ intersects the non-abelian free group $F$ trivially and is isomorphic to $\mathbb{Z}$ (Proposition~\ref{prop:periodic_implies_centre}), it maps injectively to an infinite cyclic subgroup of $M^\ab$. This shows that $\overline{M}^\ab \simeq M^{\ab}/Z(M)$  and $\ZZ^{\rank(\overline{M}^\ab)}\simeq \ZZ^{\rank(M^\ab)-1}$.
\end{proof}

\bigskip

Recall that $\tilde p_Z \colon \Aut(M) \to \Aut(\overline{M})$ is the homomorphism taking an automorphism of $M$ to its induced automorphism of the quotient $\overline{M}$.

\bigskip

\begin{lemma}[{cf. Levitt \cite[Lemma~4.1]{LevittGBS}}]
\label{lemma:unique-lift}Let $M$ be a free-by-cyclic group with finite-order monodromy and non-abelian fibre $F$.

Every $\alpha\in\Aut^{\overline{\chi}}(\overline{M})$ has a unique lift to an automorphism $\psi\in\Aut^+(M,F)$ such that $\tilde p_Z(\psi)=\alpha$.
\end{lemma}

\begin{proof}
Uniqueness is clear: if $\psi_1,\psi_2\in\Aut^+(M,F)$ project to the same automorphism $\alpha$, then for every $g\in M$, $\psi_1(g)$ and $\psi_2(g)$ differ by an element of the centre. As $\chi(g)=\chi(\psi_1(g))=\chi(\psi_2(g))$, it follows that $\psi_1(g)=\psi_2(g)$.

We now tackle existence. That is, given an automorphism $\alpha\in\Aut^{\overline{\chi}}(\overline{M})$, we want to define an automorphism $\psi$ that preserves $\chi$ and projects to $\alpha$. 

The fact that $\alpha$ leaves $\overline{\chi}$ invariant means it restricts to an automorphism of $\overline{F}$. Since the quotient $M \to \overline{M}$ restricts to an isomorphism $F \to \overline{F}$, we can use the inverse to define $\psi|_F=\alpha|_{\overline{F}}$. It remains to define $\psi(s)$. Towards this, note that the $\alpha$-invariance of $\overline{\chi}$ implies that $\alpha(\bar{s})=\bar{g}\bar{s}$ for some $\bar{g} \in \overline{F}$. There is a unique $g \in F$ that projects to $\overline{g}$, and so we set $\psi(s)=gs$.

We now need to check that $\psi$ induces a character-preserving automorphism of $M$.  We first check the relations are satisfied: let $h \in F$, and observe that since $\psi$ induces the automorphism $\alpha$ of $\overline{M}$, we must have that $\psi(s)\psi(h)\psi(s)^{-1}=\psi(\varphi(h))z^k$ in $M$. But by construction the exponent sum on $s$ in $\psi(s)\psi(h)\psi(s)^{-1}$ is $0$, so $0=k$, and it follows that $\psi(s)\psi(h)\psi(s)^{-1}=\psi(\varphi(h))$ holds, as required. It is clearly surjective, and it has a surjective inverse given by the lift of $\alpha^{-1}$. It is character preserving by construction: the exponent sum on $s$ does not change. 
\end{proof}

\bigskip

Recall that if $H \subseteq G$ then a homomorphism $R \colon G \to H$ is a \emph{retract} if $R$ is the identity on $H$.

\bigskip

\begin{prop} 
\label{prop:retract}Let $M$ be a free-by-cyclic group with finite-order monodromy and non-abelian fibre $F$. The lift $\alpha \mapsto \psi$ from Lemma~\ref{lemma:unique-lift} induces an isomorphism \[\tilde\tau_Z\colon \Aut^{\overline{\chi}}(\overline{M})\to\Aut^+(M,F)\] which is inverse to $\tilde p_{Z}|_{\Aut^+(M,F)}$.  It follows that the composition $\tilde R= \tilde \tau_Z \circ \tilde p_Z$ gives a retract \[ \tilde R \colon \Aut^{\overline{\chi}}(M) \to \Aut^+(M,F) \] with kernel equal to $\ker(\tilde p_Z)\cap \Aut^{\overline{\chi}}(M)$. The map $\tilde \tau_Z$ induces an isomorphism \[\tau_Z\colon \Out^{\overline{\chi}}(\overline{M})\to\Out^+(M,F),\] and the retract $\tilde R$ descends to a retract \[ R \colon \Out^{\overline{\chi}}(M) \to \Out^+(M,F) \] which satisfies $R=\tau_Z\circ p_Z$ and has kernel equal to $\ker(p_Z)\cap \Out^{\overline{\chi}}(M)$. This is the free abelian group $N$ given in Lemma~\ref{lem:kernel_projection}.
\end{prop}

\begin{proof}
For $\alpha \in \Aut^{\overline{\chi}}(\overline{M})$, let $\tilde \tau_Z(\alpha) \in \Aut^+(M,F)$ be the lift constructed in Lemma~\ref{lemma:unique-lift}. Then $\tilde \tau_Z$ is a homomorphism as $\tilde \tau_Z(\alpha_1)\tilde \tau_Z(\alpha_2)$ is a lift of $\alpha_1\alpha_2$, so by uniqueness of lifts $\tilde \tau_Z(\alpha_1)\tilde \tau_Z(\alpha_2)=\tilde \tau_Z(\alpha_1\alpha_2)$. If $\tilde\tau_Z(\alpha)=1$ then $\alpha$ lifts to the identity so applying $\tilde{p}_Z$ shows that $\alpha=1$.
It follows that $\tilde \tau_Z$ is injective. The composition $\tilde \tau_Z \circ \tilde p_Z$ is the identity on $\Aut^+(M,F)$ as $\psi \in \Aut^+(M,F)$ is a lift of $\tilde p_Z(\psi)$, so $\tilde R= \tilde \tau_Z \circ \tilde p_Z$ gives a retract \[ \tilde R \colon \Aut^{\overline{\chi}}(M) \to \Aut^+(M,F), \] with kernel equal to $\ker(\tilde p_Z)\cap \Aut^{\overline{\chi}}(M)$ as $\tilde \tau_Z$ is injective. We thus have an induced isomorphism $\tau_Z\colon \Out^{\overline{\chi}}(\overline{M})  \to \Out^+(M,F)$ as by the unique lifting property $\tilde \tau_Z$ maps $\Inn(\overline{M})$ isomorphically to $\Inn(M)$. Furthermore, we can also take the quotient by $\Inn(M)$ to obtain an induced retract \[ R \colon \Out^{\overline{\chi}}(M) \to \Out^+(M,F), \] such that $\ker(R)$ is the image of $\ker(\tilde R)$ in $\Out^{\overline{\chi}}(M)$, which is $\ker(p_Z)\cap \Out^{\overline{\chi}}(M)$.
\end{proof}

We finally need to understand the behaviour of the retraction with respect to peripheral subgroups. Given a subgroup $H\subseteq M$, we denote by $\overline{H}$ its image in $\overline{M}$. Given a set $\calh$ of conjugacy classes of subgroups of $M$, we let $\overline{\calh}$ be the induced set of conjugacy classes of subgroups of $\overline{M}$, i.e.\ the set of all $\overline{M}$-conjugacy classes of subgroups $\overline{H}$, with $H$ varying among all subgroups of $M$ whose conjugacy classes lie in $\calh$.

\begin{lemma}[{cf. Andrew--Martino \cite[Proposition 5.3.6]{AMFbyZ}}]
\label{lemma:andrewmartino}Let $M$ be a free-by-cyclic group with finite-order monodromy and non-abelian fibre $F$. For every set $\calh$ of conjugacy classes of subgroups of $M$, we have
\[p_Z(\Out^{\overline{\chi}}(M,\calh^{(t)}))=\Out^{\overline{\chi}}(\overline{M},\overline{\calh}^{(t)}) \quad \text{and} \quad 
R(\Out^{\overline{\chi}}(M,\calh^{(t)}))
=\Out^{\overline{\chi}}(M,F,\calh^{(t)}).\] 
\end{lemma}

\begin{proof}
 We will prove the following two inclusions:
\begin{equation}\label{inclusion1}
p_Z(\Out^{\overline{\chi}}(M,\calh^{(t)})) \subseteq \Out^{\overline{\chi}}(\overline{M},\overline{\calh}^{(t)});
\end{equation}
\begin{equation}\label{inclusion2}
\tau_Z(\Out^{\overline{\chi}}(\overline{M},\overline{\calh}^{(t)})) \subseteq \Out^{\overline{\chi}}(M,F,\calh^{(t)}).
\end{equation}
This will be sufficient to conclude, as we will then have 
\[
\begin{array}{rlrr}
\Out^{\overline{\chi}}(M,F,\calh^{(t)}) &= \tau_Z\circ p_Z(\Out^{\overline{\chi}}(M,F,\calh^{(t)})) &\quad&
\\ & \subseteq \tau_Z\circ p_Z(\Out^{\overline{\chi}}(M,\calh^{(t)})) & \quad&
\\ & \subseteq  \tau_Z(\Out^{\overline{\chi}}(\overline{M},\overline{\calh}^{(t)})) &\quad & \text{by \eqref{inclusion1}} \\ & \subseteq  \Out^{\overline{\chi}}(M,F,\calh^{(t)}) &\quad & \text{by \eqref{inclusion2}}
\end{array}
\]
and all these inclusions are then equalities. In particular, the inclusions in \eqref{inclusion1} and \eqref{inclusion2} are equalities, using for \eqref{inclusion1} that $\tau_Z$ is injective.

To prove \eqref{inclusion1}, note that if $\alpha$ represents an element of $\Out^{\overline{\chi}}(M,\calh^{(t)})$, then on each subgroup $H\subseteq M$ whose conjugacy class belongs to $\calh$ it agrees with some $\ad_g$. But then the same is true of the action of $\tilde p_Z(\alpha)$ on $\overline{H}$, which shows that \eqref{inclusion1} holds. 

To prove \eqref{inclusion2}, let $\overline{\alpha} \in \Aut(\overline{M})$ represent an element of $\Out^{\overline{\chi}}(\overline{M},\overline{\calh}^{(t)})$. Let $H\subseteq M$ be a subgroup whose conjugacy class belongs to $\calh$ (so that the conjugacy class of $\overline{H}$ belongs to $\overline{\calh}$). Let $\overline{g}_H$ be an element of $\overline{M}$ so that $\overline{\alpha}\circ\ad_{\overline{g}_H}$ restricts to the identity on $\overline{H}$, and let $g_H$ be a lift of $\overline{g}_H$ to $M$. Letting $\alpha:=\tilde\tau_Z(\overline{\alpha})$, the automorphism $\overline{\alpha}\circ\ad_{\overline{g}_H}$ lifts to $\alpha\circ\ad_{g_H}$, which sends each $h\in H$ to $hz^k$ for some $k\in\ZZ$. But in order to be fibre preserving, we must have $k=0$, and so $\alpha\circ\ad_{g_H}$ acts as the identity on $H$. This holds for each $H$ whose conjugacy class belongs to $\calh$, so $\alpha$ represents an outer automorphism in $\Out(M,F,\calh^{(t)})$, as needed.
\end{proof}

\begin{defn}[The finite-index subgroup $\Out^\per(M)$ of $\Out(M)$] 
\label{def:out_per_absolute}We let $\Out^\per(M)$ be the finite-index subgroup of $\Out(M)$ defined in the following way. 
\begin{itemize}
\item If $M$ be a free-by-cyclic group with finite-order monodromy and non-abelian fibre then we define
\[\Out^\per(M):=R^{-1}(\Out^\per(M,F)),\]
where $R \colon \Out^{\overline{\chi}}(M) \to \Out^+(M,F)$ is the retract constructed in Proposition~\ref{prop:retract} and $\Out^\per(M,F)$ is the finite-index subgroup of $\Out^+(M,F)$ from Defintion~\ref{dfn:out_0_finite_monodromy}.
\item If $M$ is the Klein bottle group or $M\simeq \mathbb{Z}$ we define $\Out^\per(M)=\{1\}$.
\item If $M \simeq \mathbb{Z}^2$ we define $\Out^\per(M)$ to be the torsion-free, finite-index subgroup consisting of matrices which are the identity in ${\mathrm{GL}}(2,\mathbb{Z}/3\mathbb{Z})$.
\end{itemize}
\end{defn}

As usual we adopt the notation 
\[\Out^\per(M,\calh):=\Out^\per(M)\cap\Out(M,\calh),\]
\[\Out^\per(M,\calh^{(t)}):=\Out^\per(M)\cap\Out(M,F,\calh^{(t)}).\]

Note that as $R$ is a retract we have $\Out^\per(M)\cap\Out^+(M,F)=\Out^\per(M,F)$ when $M$ has non-abelian fibre. 

\begin{proof}[Proof of Theorem~\ref{thm:Out(M)-type-VF-finite-order}]
Suppose first that $F$ is abelian. Then $M \simeq \ZZ, \ZZ^2$ or $K$, the fundamental group of a Klein bottle. For $\ZZ$ and $K$, the outer automorphism group is finite, and as we have set $\Out^\per(M)=\{1\}$ this suffices. If $M\simeq \ZZ^2$, let $\mathcal{H}$ be a finite set of conjugacy classes of finitely generated subgroups of $M$. Note that, if $\mathcal{H}\neq \varnothing$, then $\Out(M,\mathcal{H}^{(t)})$ is virtually cyclic. Thus, any torsion-free finite-index subgroup of $\Out(M,\mathcal{H}^{(t)})$ is of type $\F$. Similarly, as $\Out(M)$ is virtually free, any torsion-free finite-index subgroup of $\Out(M)$ is of type $\F$. Thus, any torsion-free finite-index subgroup $\Out^\per(M)$ of $\Out(M)$ gives the desired result (we chose a specific congruence kernel in Definition~\ref{def:out_per_absolute} for concreteness).

Now assume $F$ is non-abelian (equivalently, $M$ is not virtually abelian). In this case, as $\Out^\per(M)\cap\Out^+(M,F)=\Out^\per(M,F)$, the second point of the theorem is the content of Theorem~\ref{thm:Mccoolgroupsfiniteordermonodromies}, so we focus on the first. Let \[ R \colon \Out^{\overline{\chi}}(M) \to \Out^+(M,F) \] be the retract constructed in Proposition~\ref{prop:retract}, so that $\Out^\per(M)=R^{-1}(\Out^\per(M,F))$. Then restricting to $R$ to $\Out^\per(M)$ gives a short exact sequence

\[1\to \ZZ^{\rank(\overline{M}^\ab)} \to \Out^\per(M) \to \Out^\per(M,F) \to 1.\]

Notice that the rank of the kernel has not dropped as $\Out^\per(M)$ has finite index in $\Out^{\overline{\chi}}(M)$. Then for every finite set of conjugacy classes of subgroups $\calh$, we claim that \[ R\left(\Out^{\per}(M,\calh^{(t)})\right)=\Out^\per(M,F,\calh^{(t)}).\] As  $\Out^\per(M,F,\calh^{(t)}) \subseteq \Out^{\per}(M,\calh^{(t)})$ and $R$ is a retract we have $\Out^\per(M,F,\calh^{(t)}) \subseteq R(\Out^{\per}(M,\calh^{(t)}))$. Conversely, \[ R(\Out^\per(M,\calh^{(t)})) \subseteq \Out^\per(M,F) \cap \Out(M,F,\calh^{(t)})=\Out^\per(M,F,\calh^{(t)}) \]
by the definition of $\Out^\per(M)$ and Lemma~\ref{lemma:andrewmartino}, giving the opposite containment.

We can then restrict the retraction $R$ to $\Out^\per(M,\calh^{(t)})$ to give a (split) short exact sequence \[ 1 \to \ZZ^k \to \Out^\per(M,\calh^{(t)}) \to \Out^\per(M,F,\calh^{(t)}) \to 1, \] where $k \geq 0$.
Both image (by Theorem~\ref{thm:Mccoolgroupsfiniteordermonodromies}) and kernel are of type $\F$, and so this is also true of $\Out^\per(M,\calh^{(t)})$ by Proposition~\ref{prop:actiononcontractiblespacetypeVF}. 
\end{proof}

\section{More background on group actions on trees}\label{sec:arboriculture}

\subsection{Splittings, acylindricity, compatibility}\label{sec:def_splitting}
Let $G$ be a finitely generated group, let $\mathcal{H}$ be a set of conjugacy classes of subgroups of $G$, and let $\cala$ be a set of conjugacy classes of subgroups 
of $G$ which is closed under taking subgroups. A \emph{$(G,\calh;\cala)$-splitting} is a simplicial tree $T$ equipped with a minimal $G$-action (i.e.\ without any proper invariant subtree), such that: 
\begin{itemize}
    \item edge stabilisers of $T$ belong to $\cala$;
    \item every subgroup $H\subseteq G$ whose conjugacy class belongs to $\calh$ fixes a point in $T$.
\end{itemize}

We allow for the trivial action on a point in the definition. We also assume that the action is without \emph{edge inversion}: there does not exist an element of $G$ which preserves an edge while flipping its endpoints. In this paper $(G,\calh;\cala)$-splittings are always considered up to $G$-equivariant homeomorphism. A splitting is \emph{degenerate} if it is reduced to a point, and is \emph{non-degenerate} otherwise. 

We will mostly work with the following sets of subgroups for $\cala$:

\begin{itemize}
    \item $\cala$ consists of all subgroups of $G$. In that case, we will refer to a $(G,\calh;\cala)$-splitting simply as a $(G,\calh)$-splitting.
    \item $\cala$ consists of the trivial subgroup. In that case, a $(G,\calh;\cala)$-splitting is simply called a \emph{$(G,\calh)$-free splitting}.
    \item $\cala$ consists of all cyclic (possibly trivial) subgroups of $G$. Such a splitting is called a \emph{$(G,\calh)$-cyclic splitting}.
    \item $\cala$ consists of all slender subgroups of $G$ (recall that a group is \emph{slender} if all its subgroups are finitely generated). Such a splitting is called a \emph{$(G,\calh)$-slender splitting}.
    
    \item $\cala$ consists of all \emph{elementary subgroups} of $(G,\calp)$, in the case that $G$ is hyperbolic relative to a finite set $\calp$ of conjugacy classes of finitely generated subgroups. A subgroup is contained in $\cala$ if and only if it is virtually cyclic or contained in a subgroup whose conjugacy class is in $\calp$.
\end{itemize}

Let $K\geq 0$. A $(G,\calh)$-splitting $T$ is \emph{$K$-acylindrical} if there exists $C\ge 0$ such that the pointwise stabiliser of any geodesic segment in $T$ of length greater than $K$ has cardinality at most $C$. It is \emph{acylindrical} if it is $K$-acylindrical for some $K\geq 0$. For instance, all $(G,\calh)$-free splittings are $0$-acylindrical.

Let $T$ and $T'$ be two $(G,\calh)$-splittings. We say that \emph{$T$ refines $T'$} if $T'$ is obtained from $T$ by equivariantly collapsing some orbits of edges. We say that $T$ and $T'$ are \emph{compatible} if they have a common refinement.  We say that $T$ \emph{dominates} $T'$ if there exists a $G$-equivariant map $T\to T'$ (equivalently, every point stabiliser in $T$ is elliptic in $T'$).

\subsection{Stabilisers of trees for the action of the automorphism group}\label{sec:def_auto}

The group $\Aut(G,\calh)$ acts on the set of $(G,\calh)$-splittings by precomposition of the action. As $(G,\calh)$-splittings are considered up to $G$-equivariant homeomorphism, this action passes to the quotient to give an action of $\Out(G,\calh)$ on the set of $(G,\calh)$-splittings. 

Let $T$ be a $(G,\calh)$-splitting which is not a line. This implies in particular that $Z(G)$ fixes $T$ pointwise, so that $G \to \Isom(T)$ factors through $\Inn(G) \simeq G/Z(G) \to \Isom(T)$. We use $\Aut_T(G,\calh)$ and $\Out_T(G,\calh)$ to denote the stabilisers of $T$ in $\Aut(G,\calh)$ and $\Out(G,\calh)$ respectively. Let $\alpha \in \Aut_T(G,\calh)$. There exists an isomorphism $I_\alpha:T\to T$ which is $\alpha$-equivariant in the following sense: for every $g\in G$ and every $x\in T$, one has $I_{\alpha}(gx)=\alpha(g)I_\alpha(x)$ (see~\cite[Section~2.1]{Lev}). Since $T$ is not a line, $I_\alpha$ is unique. Therefore the map sending $\alpha$ to $I_\alpha$ defines a homomorphism $\Aut_T(G,\calh) \to \Isom(T)$ extending the action of $\Inn(G)$.

Let $\Out_T^0(G,\calh)$ be the maximal finite-index subgroup of $\Out_T(G,\calh)$ acting trivially on the quotient graph $T/G$. For every $v \in VT$, there exists a restriction homomorphism \[\Out_T^0(G,\calh) \to \Out(G_v)\] defined as follows. Let $\Phi \in \Out^0_T(G,\calh)$. Let $\phi_v$ be a representative of $\Phi$ such that $I_{\phi_v}$ fixes $v$. Then $\phi_v(G_v)=G_v$ and we send $\Phi$ to the outer class of $\phi_{v|G_v}$ in $\Out(G_v)$. As $\phi_v$ is well defined up to an element of $G_v$, this gives a well-defined homomorphism $\Out_T^0(G,\calh) \to \Out(G_v)$. For every $e \in ET$, we similarly construct a restriction homomorphism $\Out_T^0(G,\calh) \to \Out(G_e)$. 

If $v \in VT$, we let $\Inc_v$ be the finite set of $G_v$-conjugacy classes of incident edge stabilisers and let $\calh_v^0$ be the set of $G_v$-conjugacy classes of subgroups $H \subseteq G_v$ such that $v$ is the unique fixed point of $H$ and $[H] \in \calh$. The set $\calh_v^0$ is finite whenever $\mathcal{H}$ is (see \cite[Remark~4.6]{GLjsj}). The set $\Inc^\calh_v=\Inc_v \cup \calh_v^0$ is preserved by the image of $\Out_T^0(G,\calh)$ in $\Out(G_v)$. In other words, the image of $\Out_T^0(G,\calh) \to \Out(G_v)$ is contained in $\Out(G_v,\Inc^\calh_v)$.

\begin{lemma}
\label{lemma:restriction}Let $G$ be a group, $T$ be a $G$-tree, let $\calh$ be a set of conjugacy classes of subgroups of $G$, and let $K$ be the kernel of the homomorphism \[\Out^0_T(G,\calh^{(t)})\to\prod_{e\in ET/G}\Out(G_e).\]
Then for every $v\in VT$, the restriction homomorphism $K\to\Out(G_v)$ takes its values in $\Out(G_v,(\Inc^\calh_v)^{(t)})$.
\end{lemma}

\begin{proof}
Let $\Phi\in K$, and let $H\subseteq G_v$ be a subgroup whose conjugacy class belongs to $\Inc^\calh_v$. 

$\drsh$ Assume first that the $G_v$-conjugacy class of $H$ belongs to $\Inc_v$. Then we can find an edge $e\in ET$ incident on $v$ such that $H=G_e$. Let $\phi\in\Aut(G)$ be a representative of $\Phi$ that fixes $e$ (in particular $\phi$ fixes $v$). Since $\Phi\in K$, it follows that $\phi$ restricts to an inner automorphism on $H=G_e$. In particular, the restriction $\phi|_{G_v}$ acts by a conjugation (by an element of $G_e\subseteq G_v)$ on $H=G_e$. 

$\drsh$ Assume now that the $G$-conjugacy class of $H$ belongs to $\calh$, and $v$ is the unique vertex fixed by $H$. Let $\phi\in\Aut(G)$ be a representative of $\Phi$ that fixes $v$. Since $\Phi$ belongs to $\Out(G,\calh^{(t)})$, the automorphism $\phi$ restricts to conjugation by some $ g \in G$ on $H$, and since $\phi(G_v)=G_v$ we must have $gHg^{-1}\subseteq G_v$. But $gv$ is the unique vertex fixed by $gHg^{-1}$, so $g\in G_v$. Thus, the restriction $\phi|_{G_v}$ acts on $H$ by a conjugation by an element of $G_v$.
\end{proof}

\subsection{JSJ decompositions of groups}

In this section, we review general background regarding JSJ decompositions of groups. We mostly follow~\cite{GLjsj}. Let $T$ be a $(G,\calh;\cala)$-splitting. The tree $T$ is \emph{universally edge elliptic} (simply called ``universally elliptic'' in \cite{GLjsj}) if its edge stabilisers are elliptic in every $(G,\calh;\cala)$-splitting. A \emph{JSJ $(G,\calh;\cala)$-splitting} is a $(G,\calh;\cala)$-splitting which is universally edge elliptic and which dominates every universally edge elliptic $(G,\calh;\cala)$-splitting. 

\begin{rk}\label{rk:universally-vertex-elliptic}
Note that if a tree $T$ is \emph{universally vertex elliptic} in the sense that every vertex stabiliser of $T$ is elliptic in every $(G,\calh;\cala)$-splitting, then automatically $T$ is a JSJ $(G,\calh;\cala)$-splitting. The converse is not true in general.
\end{rk}

JSJ $(G,\calh;\cala)$-splittings exist whenever $G$ is finitely presented and $\calh$ consists of a finite set of conjugacy classes of finitely generated subgroups, see~\cite[Theorem~2.20]{GLjsj}. JSJ $(G,\calh;\cala)$-splittings are in general not unique but they all belong to the same \emph{deformation space} (see~\cite{GL-defo}), that is, they all dominate each other. The set of all JSJ $(G,\calh;\cala)$-splittings is called the \emph{JSJ deformation space of $(G,\calh;\cala)$}. The JSJ deformation space is $\Aut(G,\calh)$-invariant whenever $\cala$ is. 

Let $T$ be a JSJ $(G,\calh;\cala)$-splitting and let $v \in VT$. We say that $v$ is \emph{rigid} if $G_v$ is elliptic in every $(G,\calh;\cala)$-splitting; otherwise it is \emph{flexible}. In particularly nice situations, each flexible vertex $G_v$ is \emph{quadratically hanging (QH) with finite fibre}: there is a compact hyperbolic surface $\Sigma_v$ and a finite subgroup $\mathrm{Fin} \subseteq G_v$ such that $G_v$ fits into an exact sequence $1 \to \mathrm{Fin} \to G_v \to \pi_1(\Sigma_v) \to 1$ and the image in $\pi_1(\Sigma_v)$ of every incident edge group is conjugate into a boundary subgroup of $\Sigma_v$. For instance, flexible vertices are QH with finite fibre when $G$ is one-ended and $\cala$ is the set of all virtually cyclic subgroups (\cite[Theorem~9.21]{GLjsj}), or when $G$ is one-ended relatively hyperbolic and $\cala$ is the set of conjugacy classes of elementary subgroups of $G$ (\cite[Theorem~9.18]{GLjsj}). Notice that if $T$ is non-degenerate and a vertex $v \in VT$ is QH with finite fibre, then $G_v$ is virtually free.

\subsection{Transverse coverings and their skeleta}\label{sec:skeleta}

Let $T$ be a $(G,\calh)$-splitting. Recall from \cite[Definition~4.6]{Guirardellimit} that a \emph{transverse covering} of $T$ is a $G$-invariant set $\caly$ of non-degenerate subtrees of $T$ such that any two distinct subtrees in $\caly$ intersect in at most one point, and every non-degenerate arc of $T$ is covered by finitely many subtrees in $\caly$. 

Recall from \cite[Definition~4.8]{Guirardellimit} that the \emph{skeleton} of $\caly$ is then the bipartite $G$-tree $T_\caly$ with vertex set $VT_\caly=V_0\amalg V_1$, having one vertex $v\in V_0$ per vertex $x\in VT$ that does not belong to exactly one subtree in $\caly$, one vertex $w\in V_1$ per subtree $Y\in\caly$, with an edge joining $v\in V_0$ (representing $x$) and $w\in V_1$ (representing $Y$) whenever $x\in Y$.

As $\caly$ is $G$-invariant, the skeleton of $\caly$ also comes equipped with an action of $G$. By \cite[Lemma~4.9]{Guirardellimit}, this action is minimal. More generally, if $H$ is a subgroup of $\Aut_T(G,\calh)$ preserving $\caly$, then the action of $H$ on $T$ induces an action of $H$ on $T_\caly$. 

As shown by the following lemma, acylindricity is preserved under passing to a skeleton.

\begin{lemma}
\label{lemma:skeleton-acylindrical}If the $G$-action on $T$ is $K$-acylindrical and $\mathcal{Y}$ is a transverse covering of $T$, then the $G$-action on the skeleton of $\mathcal{Y}$ is $(2K+2)$-acylindrical.
\end{lemma}

\begin{proof}
Let $\gamma$ be an edge path of length $2K+3$ in the skeleton $S$. Since $S$ is bipartite, $\gamma$ contains a subpath of length $2K+2$ whose endpoints are in $V_0$; denote the vertices along this path by $v_0, Y_1, v_1, \dots, Y_{K+1}, v_{K+1}$. In $T$, the arc $[v_{i-1}, v_i]$ is contained in $Y_i$, so there is no overlap between $[v_{i-1}, v_i]$ and $[v_i, v_{i+1}]$. It follows that the geodesic $[v_0,v_{K+1}]$ is the concatenation of the geodesics $[v_i,v_{i+1}]$; in particular $[v_0,v_{K+1}]$ contains each $v_i$. The vertices $v_i$ are pairwise distinct, so $[v_0,v_{K+1}]$ has length at least $K+1$. But any element $g$ stabilising $\gamma$ also stabilises $[v_0,v_{K+1}]$, and so this stabiliser has bounded cardinality (with the same bound as for $T$).
\end{proof}

\subsection{Trees of cylinders}\label{sec:tree-of-cyl}

The material in this section is taken from \cite{GL-cylinders}. Let $\cale \subseteq \cala$ be a subset of conjugacy classes of subgroups. We do not insist that $\cale$ is stable under taking subgroups, contrary to the set $\cala$. A \emph{$(G,\calh;\cale)$-splitting} is then a $(G,\calh;\cala)$-splitting with edge stabilisers in $\cale$. An equivalence relation $\sim$ on $\cale$ is \emph{admissible} if the following hold:
\begin{itemize}
    \item if $A \sim B$ and $g \in G$, then $gAg^{-1} \sim gBg^{-1}$;
    \item if $A \subseteq B$, then $A \sim B$;
    \item if $T$ is a $(G,\calh;\cale)$-splitting, if $A \sim B$ and $A$, $B$ fix respectively $a,b \in VT$, then, for every edge $e \subseteq [a,b]$, we have $A \sim G_e \sim B$. 
\end{itemize}

For instance, if $\cale$ consists of all infinite virtually cyclic subgroups of $G$, then commensurability is an admissible equivalence relation.

Now let $T$ be a $(G,\calh;\cale)$-splitting and let $\sim$ be an admissible equivalence relation on $\cale$. We define an equivalence relation $\sim_T$ on the set of edges of $T$ by setting $e \sim_T e'$ if and only if $G_e \sim G_{e'}$. A \emph{cylinder of $T$} is a subforest of $T$ consisting of all edges in a given $\sim_T$-equivalence class. A cylinder is in fact a subtree of $T$ (see~\cite[Lemma~4.2]{GL-cylinders}). The \emph{tree of cylinders of $T$ relative to $\cale$}, denoted by $T_c$, is the skeleton of the transverse covering consisting of all the cylinders. In general edge stabilisers of $T_c$ do not belong to $\cale$. The \emph{collapsed tree of cylinders of $T$}, denoted by $T_c^*$, is the tree obtained from $T_c$ by collapsing every edge of $T_c$ whose stabiliser is not in $\cale$.

\begin{lemma}
\label{lemma:tree-of-cylinders}Let $T$ be a $(G,\calh;\cale)$-splitting, let $\sim$ be an admissible equivalence relation on $\cale$, let $T_c$ be the associated tree of cylinders and let $T_c^*$ be the collapsed tree of cylinders. Assume that all cylinders in $T$ are bounded. Then:
\begin{enumerate}
\item\label{collapse1} The trees $T$, $T_c$ and $T_c^*$ have the same elliptic subgroups.
\item\label{collapse2} Every vertex stabiliser of $T_c^*$ is either a vertex stabiliser in $T$ or belongs to $\mathcal{E}$.
\item\label{collapse3} Every edge stabiliser of $T_c^*$ is $\sim$-equivalent to an edge stabiliser in $T$.
\end{enumerate}
\end{lemma}

\bigskip

\begin{proof}
The first point is \cite[Proposition~5.2 and Corollary~5.10]{GL-cylinders}. The proof of Assertions~(\ref{collapse2}) and~(\ref{collapse3}) follows \cite[Proposition~5.9]{GL-cylinders} closely; we include a proof for the convenience of the reader.

Let $C\subseteq T$ be a cylinder. Since $C$ is bounded, $G_C$ fixes a point $y\in C$, namely the circumcentre of $C$. Let now $\varepsilon=(x,C)$ be an edge in $T_c$ having $C$ as one of its endpoints (the other being given by a vertex $x\in VT$ belonging to $C$). Note that $G_\varepsilon=G_x \cap G_C$. We first show Assertion~(\ref{collapse3}), which is equivalent to saying that if $\varepsilon$ is not collapsed in $T_c^*$ then $G_\varepsilon$ is $\sim$-equivalent to an edge stabiliser in $T$.

$\drsh$ If $y\neq x$, let $e\subseteq C$ be the initial edge of the segment $[x,y]$. Then $G_\varepsilon\subseteq G_e$. On the other hand $G_e$ fixes $x$ and preserves $C$, so $G_e\subseteq G_\varepsilon$. Hence $G_e=G_\varepsilon$ (in particular $\varepsilon$ is not collapsed in $T_c^*$). 

$\drsh$ If $y=x$, then $G_C\subseteq G_x$ and therefore $G_\varepsilon=G_C$. In particular, if $e\in ET$ is an edge contained in the cylinder $C$, then as $G_e\subseteq G_C$ we have $G_e \subseteq G_\varepsilon$. So if $\varepsilon$ is not collapsed in $T_c^*$, then $G_e\sim G_\varepsilon$ by the second bullet point in the definition of admissible.

The above proves Assertion~(\ref{collapse3}); it also shows, as in \cite[Proposition~5.9]{GL-cylinders}, that at most one edge with endpoint $C$ can be collapsed in $T_c^*$ -- namely the edge $(y,C)$, if it exists (i.e.\ if $y$ belongs to another cylinder), and in this case $G_C\subseteq G_y$.

To prove the second assertion, let $v\in VT_c^*$, and let $X_v$ be the preimage of $v$ under the collapse map $T_c\to T_c^*$. 

$\drsh$ We first assume that $X_v$ contains a vertex $x$ corresponding to a point of $T$. In this case, the above shows that $X_v$ is contained in the star of $x$, and that every cylinder vertex $C$ in $X_v$ satisfies $G_C\subseteq G_x$. In particular $G_v=G_x$ is a point stabiliser in $T$.

$\drsh$ Otherwise, $X_v$ is reduced to one vertex corresponding to a cylinder $C$. In this case, either the circumcenter $y$ of $C$ is contained in a unique cylinder (in which case $G_C=G_y$), or else there is an edge joining $C$ to $y$ in $T_c$ that does not get collapsed in $T_c^*$, and by the above this implies that $G_C\in\mathcal{E}$. So $G_v=G_C$ is either a point stabiliser of $T$ or belongs to $\mathcal{E}$.
\end{proof}

\bigskip

\subsection{The canonical JSJ tree of a relatively hyperbolic group}

As recalled in Section~\ref{sec:growth_prelim}, free-by-cyclic groups with exponentially-growing monodromy are relatively hyperbolic. Our arguments in this case will rely on the existence of a canonical JSJ tree, whose properties are recorded in the following statement. Recall that given a group $G$ and a set $\calh$ of conjugacy classes of subgroups of $G$, one says that $G$ is \emph{one-ended} relative to $\calh$ if there is no non-degenerate splitting of $G$ relative to $\calh$ with finite edge stabilisers.

\bigskip

\begin{thm}[{\cite[Theorem~9.18]{GL-cylinders}}]
\label{thm:tree_rel_aut}Let $G$ be a group which is hyperbolic relative to a finite set $\calp$ of conjugacy classes of non-virtually-cyclic finitely generated subgroups. Let $\calh$ be a finite set of conjugacy classes of finitely generated subgroups of $G$, and assume that $G$ is one-ended relative to $\calp\cup\calh$. Let $\cala$ be the 
set of all elementary subgroups of $(G,\calp)$.

Then there exists a unique JSJ $(G,\calp\cup\calh;\cala)$-splitting $U$ that is equal to its own tree of cylinders. The tree $U$ is $\Aut(G,\calp\cup\calh)$-invariant, and it enjoys the following properties:
\begin{enumerate}
\item\label{jsj-1} \emph{\cite[Proposition~6.1]{GL-cylinders}} The tree $U$ is acylindrical as a $G$-tree.
    \item\label{jsj-2} \emph{\cite[Lemma~3.7]{GLsplittings}} Edge stabilisers of $U$,  viewed as a $G$-tree, are finitely generated.
    \item\label{jsj-3} The vertex set of $U$ has a bipartition $VU=V_0\sqcup V_1$ such that 
    \begin{itemize}
    \item the stabiliser of every vertex $v\in V_0$ is either rigid (i.e.\ elliptic in every $(G,\calp\cup\calh;\cala)$-splitting), or QH with finite fibre; 
        \item the stabiliser of every vertex in $V_1$ is either virtually cyclic, or its conjugacy class belongs to $\calp$. 
    \end{itemize}
       \item\label{jsj-4} \emph{\cite[Proposition~4.7]{GLsplittings}} For every rigid vertex $v\in V_0$, the image of the restriction homomorphism $\Out^0_U(G,(\calp\cup\calh)^{(t)})\to\Out(G_v)$ is finite.
    \item\label{jsj-5} \emph{\cite[Lemma~3.8]{GLsplittings}} For every edge $e\in EU$, the image of the restriction homomorphism $\Out^0_U(G,(\calp\cup\calh)^{(t)})\to\Out(G_e)$ is finite.
\end{enumerate}
\end{thm}

\begin{rk}
The contents of \cite[Lemma~3.8]{GLsplittings} is that the stabiliser $G_v$ of any vertex $v\in V_0$ is hyperbolic relative to the set $\calp_v=\Inc_v\cup\calp_v^0$ introduced in Section~\ref{sec:def_auto}. In particular, incident edge stabilisers are malnormal in $G_v$. This is why the fifth property stated above follows from the fourth when the incident vertex in $V_0$ is not QH (if it is QH, then the edge group is virtually cyclic and the fifth assertion also holds). 
\end{rk}

\begin{rk}\label{rk:no-qh}
In the case where $G$ is a free-by-cyclic group, QH vertices never arise \cite[Proposition~2.11]{dahmani2016suspensions}.
\end{rk}

\section{Limiting trees and their collapsed trees of cylinders}\label{sec:limiting-tree}

Let $M$ be a free-by-cyclic group with polynomial growth of degree $d\ge 1$. Our goal in this section is to build an $\Aut(M)$-invariant tree in which vertex stabilisers are  free-by-cyclic with polynomial growth of smaller degree (possibly virtually abelian). 
This tree will serve as an input for the constructions developed in Sections \ref{sec:tree-of-cylinders} and \ref{sec:dévissage}. 

The section is organised as follows. In Section~\ref{sec:limit}, we use a limiting construction that comes from the work of Cohen--Lustig and Bestvina--Feighn--Handel to build an $M$-tree $T$ where the point stabilisers have smaller growth. The tree $T$ is called \emph{the limiting tree of $\Phi$}; its construction depends on a choice of monodromy $\Phi$ for $M$. We prove the acylindricity of $T$ as an $M$-tree in Section~\ref{sec:acylindricity}, and give a fine description of its point and arc stabilisers in Section~\ref{sec:vertex-edge-stabs}. 

In Section~\ref{sec:slow} we show that $T$ is a slender JSJ tree of $M$ relative to an $\Aut(M)$-invariant set of subgroups with finitely many conjugacy classes in $M$, namely maximal free-by-cyclic subgroups with strictly smaller growth (this family is called $\Slow^{\max}(M)$ in the sequel). This shows that even though our chosen limiting tree depended on a choice of monodromy, every limiting tree of $M$ lies in the same JSJ deformation space. In Section~\ref{sec:collapsed} we obtain the required $\Aut(M)$-invariant tree as the \emph{collapsed tree of cylinders} $T^*_c$ of the limiting tree. The key properties of $T^*_c$ are given in Theorem~\ref{theo:collapsed-tree-of-cylinders}, which is the main result of this section. Section~\ref{sec:bonus} contains additional properties of the tree $T_c^*$, like its invariance under taking finite-index subgroups, that are not needed in the rest of the paper.

\subsection{A \texorpdfstring{$\Phi$}{Phi}-invariant limiting tree \texorpdfstring{$T$}{T}}\label{sec:limit}

Recall that an $F$-tree $T$ is \emph{very small} if arc stabilisers are either trivial or maximal cyclic, and $T$ has trivial tripod stabilisers.

\begin{thm}[Cohen--Lustig \cite{CL95}, Bestvina--Feighn--Handel \cite{BFHTits2}] 
\label{thm:limiting_tree_existence}Let $\Phi \in \Out(F)$ be an outer automorphism with polynomial growth of degree $d \geq 1$. Then there exists a positive power $\Psi=\Phi^k$ of $\Phi$ such that for every metric simplicial tree $S$ with a free action of $F$ and every $g\in F$, the limit \[ \|g\|_T=\lim_{n \to \infty}\frac{\|\Psi^n(g)\|_S}{n^d}, \] exists, and $g\mapsto \|g\|_T$ is the length function of a metric simplicial tree $T$. Furthermore:
\begin{itemize}
\item The underlying simplicial tree determined by $T$ is independent of $S$ (and choice of such $k$) and $\Phi$--invariant.
\item As an $F$-tree, $T$ is a free splitting if $d \geq 2$, and is a very small $F$-tree with $\mathbb{Z}$-edge stabilisers if $d=1$.
\item If $g \in F$ then the conjugacy class of $g$ has polynomial growth of degree $d$ under $\Phi$ if and only if $g$ is hyperbolic in $T$.
\end{itemize}
\end{thm}

The underlying simplicial tree of $T$ will be called the \emph{limiting tree} of $\Phi$.

\begin{proof} The final bullet point of the theorem follows directly from the definition of $\|g\|_T,$ as the growth rate of a conjugacy class is the same under both $\Phi$ and $\Psi=\Phi^k$. Given independence of the simplicial structure, one can deduce that the simplicial tree determined by $T$ is $\Phi$--invariant by taking limits as $\Psi^n(\Phi(S))=\Phi(\Psi^n(S))$. If $k$ and $k'$ are two powers of $\Phi$ for which the result holds, by considering the limit of  $\Phi^{kk'}$ one can see the obtained tree is independent of the choice of power $k$. The rest of the theorem can be extracted from \cite{CL95} and \cite{BFHTits2}; we give precise references below, referring the reader to the papers for the majority of the definitions.

When $d=1$, every linearly growing automorphism has a power $\Psi$ which is represented as an \emph{efficent Dehn twist} in the language of \cite{CL99}, and called \emph{proper with all twistors non-trivial} in \cite{CL95}. The important information we need for this proof is that an efficient Dehn twist is determined by a very small graph of groups $\mathcal{G}$ with edge groups isomorphic to $\mathbb{Z}$ and a choice of non-trivial \emph{twistor} $z_e$ in each edge group. In \cite[Theorem 13.2]{CL95}, Cohen and Lustig show that the rescaled length functions of $\Psi^n(S)$ converge to the Bass--Serre tree of $\mathcal{G}$, endowed with the metric where each edge $e$ has length $\|z_e\|_S$ (the translation length of the twistor in the original tree $S$). Although the edge lengths depend on $S$, the underlying simplicial tree does not.

When $d \geq 2$, this is a result of Bestvina--Feighn--Handel, where $\Psi$ is taken to be a \emph{UPG} power of $\Phi$. Convergence is given in \cite[Theorem 4.22]{BFHTits2}. Independence of the simplicial structure is a consequence of \cite[Lemma 4.33]{BFHTits2}: it tells us that an explicit description of the tree $T$ is obtained by assigning a length $L_G(E)$ to every edge in a \emph{UR train track} $G$, passing to the universal cover $\tilde G$, and collapsing the edges assigned length zero (in particular the limit tree $T$ is a free splitting). It is therefore enough to show that the property that $L_G(E)>0$ is independent of the tree $S$.

Writing the UR train track as $f:G \to G$, every edge $E$ satisfies \[ f^n(E)= E \cdot u_E \cdot [f(u_E)] \cdot [f^2(u_E)] \cdots [f^{n-1}(u_E)], \] where each $[f^j(u_E)]$ is a loop in $G$. Let $g_n(E)$ be an element in the conjugacy class of the loop $u_E \cdot [f(u_E)] \cdot [f^2(u_E)] \cdots [f^{n-1}(u_E)]$.

Unpacking Definition~4.27 of \cite{BFHTits2} tells us that the length of each edge in T is \[L_G(E)=\lim_{n \to \infty}  \frac{\|g_n(E)\|_S}{n^d}.\]
It follows that $L_G(E)>0$ if and only if the sequence of conjugacy classes $(\|g_n(E)\|_S)_{n \in \mathbb{N}}$ grows like $n^d$. This is a quasi-isometry invariant property, so independent of the chosen starting point $S$.
\end{proof}

\subsection{Acylindricity of the limiting tree}\label{sec:acylindricity}

From now on we will consider the limiting tree $T$ constructed in the previous section as a simplicial tree. As $T$ is $\Phi$--invariant, the $F$-action on $T$ extends to an action of $M$ as explained in Section~\ref{sec:def_auto}. It is relatively easy to show the $F$-action on $T$ is acylindrical (see Lemma~\ref{lem:very-small-acylindricity}, below). The main goal of this section is to prove Corollary~\ref{cor:acylindrical}, which states that the $M$-action is also acylindrical.

\begin{rk}\label{rem:subdivision}
Note that the $M$-action on $T$ may invert edges. For convenience we will allow this for the time being, but note that one would need to subdivide inverted edges to find the graph of groups structure for this action. In general, when we say `edge stabiliser' or write $M_e$ without any further qualification we are referring to elements that fix every point in an edge; the subgroup of elements that preserve but possibly flip an edge will be called the 
\emph{setwise edge stabiliser}. 
\end{rk}

We first note that the $F$-action on $T$ is acylindrical. 

\begin{lem}
\label{lem:very-small-acylindricity}Let $T$ be a very small simplicial $F$-tree. 
The fixed point set of every elliptic element of $F$ is either a vertex or a finite arc.  The $F$-action on $T$ is acylindrical.
\end{lem}

\begin{proof}
 We claim that if two edges of $T$ have the same infinite cyclic stabiliser, and lie in the same $F$-orbit, then they must be equal. Indeed, let $e \in ET$ and $g \in F$ be such that $F_e \cap F_{ge}$ is non-trivial. Since $F_e$ is a maximal cyclic (hence malnormal) subgroup of $F$, and $F_{ge}= gF_eg^{-1}$, we must have that $g$ is an element of $F_e$ and hence $e=ge$.
 
 The above claim implies that if $e,e'$ are distinct edges in the same $F$-orbit, then $F_e\cap F_{e'}=\{1\}$, which gives an acylindricity constant of (at worst) the number of edges in the quotient graph $T/F$. As $T$ has trivial tripod stabilisers, the fixed set of every elliptic element is then either a unique point or a finite arc.
\end{proof}

We now aim to extend acylindricity to the $M$-action. 

\begin{lemma}
\label{lemma:edge-stab}Let $G$ be a group, let $M=G\rtimes_\Phi\mathbb{Z}$, and let $\chi:M\to\mathbb{Z}$ be a character with kernel $G$. Let $T$ be an $M$-tree, and assume that the $G$-action on $T$ is cocompact.  

Then for every edge $e\in ET$, we have $\chi(M_e)\neq\{0\}$.
\end{lemma}

\begin{proof}
As $G$ is normal in $M$, the group $M$ permutes the $G$-orbits of edges. Let $M^0\subseteq M$ be the subgroup consisting of automorphisms that preserve the $G$-orbit of every edge. As the $G$-action on $T$ is cocompact, $M^0$ has finite index in $M$ (in particular it is infinite). 

Let $e\in ET$, and let $t\in M^0$ with $\chi(t)\neq 0$. Then $te$ belongs to the $G$-orbit of $e$, so there exists $g\in G$ such that $te=ge$. So $g^{-1}t$ belongs to $M_e$, and it satisfies $\chi(g^{-1}t)\neq 0$.  
\end{proof}

Our next result, based on \cite{AGH}, gives a general criterion to extend acylindricity to the $M$-action.

\begin{prop}
\label{prop:suspensions-acylindricity}Let $G$ be a torsion-free finitely generated group, acting minimally and acylindrically on a tree $T$ with cyclic (maybe trivial) edge stabilisers. Let $\Phi\in \Out(G)$, let $M=G \rtimes_\Phi \ZZ$, and suppose $T$ is $\Phi$-invariant. Then the following are equivalent: 
\begin{enumerate}
\item\label{acyl-1} $T$ is acylindrical as an $M$-tree;
\item\label{acyl-2} for every $\alpha\in\Aut(G)$ representing a power $\Phi^k$ with $k\ge 1$, the fixed subgroup $\Fix(\alpha)$ is elliptic in $T$.
\end{enumerate}
\end{prop}

\begin{proof}
\underline{(\ref{acyl-1})$\Rightarrow$(\ref{acyl-2}).} We prove the contrapositive statement.
Let $\alpha\in\Aut(G)$ be a representative of $\Phi^k$ for some $k\ge 1$, and assume that $\Fix(\alpha)$ is not elliptic in $T$. Then, as the action of $G$ on $T$ is acylindrical, by e.g.\ \cite[Lemma~4.3]{AGH} the group $\Fix(\alpha)$ contains a hyperbolic element $h$. Since $h$ belongs to $\Fix(\alpha)$, it is contained in a subgroup of $M$ isomorphic to $\mathbb{Z}^2$. This subgroup preserves the axis of $h$, and therefore contains a $\mathbb{Z}$ subgroup that fixes it pointwise. Hence the action of $M$ on $T$ is not acylindrical.

\medskip

\noindent \underline{(\ref{acyl-2})$\Rightarrow$(\ref{acyl-1}).} We start with two special cases before treating the general case. 
In what follows, we let $s\in M$ be an element with non-trivial projection to $\mathbb{Z}$, such that conjugation by $\mathbb{Z}$ induces a monodromy $\varphi$ (with outer class $\Phi$).

$\drsh$ \textbf{Case~1.} First suppose all edge stabilisers for the $F$-action on $T$ are isomorphic to $\ZZ$. In this case, suppose that $gs^k$ fixes an edge of $T$. As a consequence of Lemma~\ref{lemma:edge-stab}, the $M$-edge stabilisers are virtually $\ZZ^2$. So we find an element $h\in G_e$ which commutes with some power $(gs^k)^\ell$ with $\ell\neq 0$ (in other words $h$ is fixed by $\ad_g\circ\phi^\ell$). In particular $h$ is contained in $G_v$ for any vertex $v$ incident on $e$. So we can apply \cite[Lemma 4.5]{AGH} with $J=1$ to  conclude that the $M$-action is acylindrical (with constant $K+2$, where $K$ was the constant for the $G$-action).

$\drsh$ \textbf{Case~2.} Now suppose all edge stabilisers for the $F$-action on $T$ are trivial. In this case let $K$ be the number of edge orbits for the $G$-action. We will prove the contrapositive statement:  if $T$ is not $2K$-acylindrical as an $M$-tree, then (\ref{acyl-2}) fails. So suppose that $gs^k$ fixes an arc of length $2K+1$. Then the arc contains three edges in the same $G$-orbit: $e$, $xe$, $ye$ (with $x,y\in G$). Note that $x$ (and by symmetry $y$) commute with $gs^k$: indeed $\langle(gs^k)^x\rangle$ has non-trivial (hence finite-index) intersection with $\langle gs^k \rangle$; since conjugation preserves the exponent sum of $s$, this implies that $(gs^k)^x=gs^k$, as claimed. Setting $\alpha := \ad_g \phi^k$, we see that $\langle x,y \rangle \subseteq \Fix(\alpha)$. We claim that the subgroup $\langle x,y \rangle$ contains a hyperbolic element: indeed, as the three edges $e,xe,ye$ are in a line, at least two of them have coherent orientations, so either $x$, $y$ or $xy^{-1}$ is hyperbolic, showing that (\ref{acyl-2}) fails. 

$\drsh$ \textbf{The general case.} Observe that we have two collapses of $T$: one to $T'$ keeping only the trivially-stabilised edges, and one to $T''$ keeping only those with infinite cyclic stabiliser. If all fixed subgroups are elliptic in $T$ this remains true in both collapses, and so by the previous two cases both are acylindrical as $M$-trees. But this implies the same of $T$: since all edges survive in one of the collapsed trees, a sufficiently long arc in $T$ projects to an arc whose length exceeds the acylindricity constant in at least one of $T'$ and $T''$, and so must have trivial stabiliser.
\end{proof}

\begin{cor}
\label{cor:acylindrical}Let $M=F\rtimes_\Phi\mathbb{Z}$ be a free-by-cyclic group with polynomial growth of degree $d\ge 1$, and let $T$ be the limiting tree of $\Phi$.

As an $M$-tree, $T$ is acylindrical.
\end{cor}

\begin{proof}
The $F$-action on $T$ is acylindrical by Lemma~\ref{lem:very-small-acylindricity}. Let $\alpha\in\Aut(F)$ be a representative of a power $\Phi^k$ with $k\ge 1$. If $g\in F$ is contained in $\Fix(\alpha)$, then $g$ is elliptic in $T$ since its conjugacy class does not grow with degree $d$ under $\Phi$. Since fixed subgroups of automorphisms of free groups are finitely generated \cite{Ger} this implies that $\Fix(\alpha)$ is elliptic. So Proposition~\ref{prop:suspensions-acylindricity} allows us to conclude that $T$ is also acylindrical as an $M$-tree.
\end{proof}

\subsection{Point and arc stabilisers of the limiting tree}\label{sec:vertex-edge-stabs} 

Recall from Section~\ref{sec:growth_prelim} that the degree of a polynomially-growing free-by-cyclic group $M$ is the degree of any of its monodromies. Although we now know the $M$-action on the limiting tree is acylindrical, we will require a more detailed description of stabilisers for the $M$-action on $T$, which is the aim of this section. Proposition~\ref{prop:limiting_tree_0} is a general statement, whereas Lemma~\ref{lem:linear-image-edge-finite} is specialised to the case where the monodromies of $M$ have linear growth.

\begin{prop} 
\label{prop:limiting_tree_0}Let $M=F\rtimes_\Phi\mathbb{Z}$ be a free-by-cyclic group with polynomial growth of degree $d \geq 1$, and let $T$ be the limiting tree of $\Phi$, where $M$ acts possibly with inversions. 
\begin{enumerate}
    \item\label{limit1,5} If $H$ is a non-cyclic free-by-cyclic subgroup of $M$ with strictly smaller growth then $H$ is elliptic in $T$.
    \item\label{limit1}  Every point stabiliser for the $M$-action on $T$ is either slender or a maximal free-by-cyclic subgroup of $M$ of strictly smaller growth. 
    \item\label{limit2} Pointwise edge stabilisers for the $M$-action on $T$ are isomorphic to $\mathbb{Z}$ when $d>1$, and isomorphic to $\mathbb{Z}^2$ when $d=1$.
\end{enumerate}
\end{prop}

\begin{proof}
\textbf{Assertion~1.}
 Suppose that $H$ is a non-cyclic free-by-cyclic subgroup of $M$ with strictly smaller growth, in particular $H \cap F \neq \{1\}$. By \cite[Lemma~3.5]{MutanguhaPG}, $H\cap F$ is finitely generated, so $\chi|_H$ is a fibring of $H$. If $g \in H \cap F$ is hyperbolic, then it has growth of degree $d$ under the monodromy induced by $\chi$ by Theorem~\ref{thm:limiting_tree_existence}, so  has growth of degree $d$ under the monodromy induced by $\chi|_H$ by Lemma~\ref{lem:growth_subgroups}. As $H$ has strictly smaller growth than $M$, it follows that every element of $H \cap F$ is elliptic. As $H \cap F$ is finitely generated it has a global fixed point, and as the $F$-action on $T$ is acylindrical, the fixed subtree of $H \cap F$ is bounded. As $H \cap F$ is a normal subgroup of $H$, the group $H$ acts on the fixed point set of $H \cap F$, so preserves the circumcenter of this bounded subtree. Hence $H$ has a global fixed point for its action on $T$.

\textbf{Assertion~2.} Note that any point stabiliser $M_v$ for the $M$-action on $T$ is free-by-cyclic by Theorem~\ref{thm:limiting_tree_existence},  i.e.\ $M_v\cap F$ is finitely generated. Moreover, $M_v \cap F$ consists of elements of growth at most $d-1$ under iteration by $\Phi$. Therefore, by Lemma~\ref{lem:growth_subgroups}, the group $M_v$ is a polynomially-growing free-by-cyclic subgroup of degree at most $d-1$. Since every non-cyclic free-by-cyclic subgroup $H$ of $M$ of degree at most $d-1$ is elliptic in $T$ (Assertion~1), and since edge stabilisers of $T$ are slender  (Theorem~\ref{thm:limiting_tree_existence}), any non-slender point stabiliser of $T$ is maximal with respect to free-by-cyclic subgroups of $M$ of degree at most $d-1$. Assertion~2 follows.

\smallskip

\textbf{Assertion~3.} By Lemma~\ref{lemma:edge-stab}, the restriction of the character $\chi$ to every vertex and pointwise edge stabiliser is non-trivial. Consequently, each pointwise edge stabiliser for the $M$-action is isomorphic to $\mathbb{Z}$ when $d >1$ and isomorphic to $\mathbb{Z}^2$ or to the Klein bottle group when $d=1$. We will rule out Klein bottle edge stabilisers in Lemma~\ref{lem:linear-image-edge-finite}.
\end{proof}

As preparation for understanding the tree of cylinders, we record a lemma on the stabilisers of the $M$-action in the linear case. For the purposes of this lemma, let $\SetStab(\cdot)$ denote the setwise stabiliser in $M$ and $\PtStab(\cdot)$ the pointwise stabiliser in $M$.

\begin{lem}
\label{lem:linear-image-edge-finite}Let $M=F\rtimes_\Phi\mathbb{Z}$ be a linearly-growing free-by-cyclic group and let $T$ be the limiting tree of $\Phi$, where $M$ acts possibly with inversions.
\begin{enumerate}
    \item\label{stab-1} For every vertex $v \in VT$, its stabiliser $M_v$ is a periodic free-by-cyclic group with non-abelian fibre and infinite cyclic centre $\langle z_v \rangle.$
    \item\label{stab-2} If $v$ and $w$ are distinct then $\langle z_v \rangle \cap \langle z_w \rangle = \{1\}$.
\item\label{stab-3} If $a\in F$ has a fixed point set equal to a non-degenerate finite arc $\beta$ in $T$, then the following hold:
\begin{enumerate}
\item    $\SetStab(\beta)=N_M(\langle a\rangle)$ and $\PtStab(\beta)=C_M(a)$;
\item $C_M(a)\simeq \mathbb{Z}^2$; 
\item If $N_M(\langle a\rangle)\neq C_M(a)$ then $N_M(\langle a\rangle)$ is isomorphic to the fundamental group of a Klein bottle and acts by inversion on the arc $\beta$ via the quotient map $N_M(\langle a\rangle)/C_M(a)\simeq \mathbb{Z}/2\mathbb{Z}$. 
\item The pointwise stabiliser of each edge in $\beta$ is equal to $C_M(a)\simeq \mathbb{Z}^2$.
\end{enumerate}
\end{enumerate}
\end{lem}

\begin{proof}By Proposition~\ref{prop:limiting_tree_0}, each vertex group $M_v$ is a periodic free-by-cyclic group. As $T$ is very small and has infinite cyclic $F$-edge stabilisers, the $F$-stabiliser $F_v$ of any vertex is non-abelian (and equal to the fibre of $M_v$ under $\chi$). So by Proposition~\ref{prop:periodic_implies_centre}, the group $M_v$ has non-trivial infinite cyclic centre $\langle z_v \rangle$, and the final part of Proposition~\ref{prop:periodic_implies_centre} ensures that $z_v$ has non-trivial image under our given character $\chi:M\to\ZZ$. We can and do choose all the $z_v$ to have positive image under $\chi$ (such a choice is compatible with the $M$-action, since conjugations are character preserving, but does depend on $\chi$). If $v$ and $w$ are distinct, then $\langle z_v \rangle \cap \langle z_w \rangle = \{1\}$: if $z_v^k=z_w^\ell$ for some (non-zero) powers $k$ and $\ell$, the group $\langle F_v,F_w\rangle$ would consist entirely of periodic conjugacy classes and contain a hyperbolic element, contradicting the fact that hyperbolic elements in $F$ have linear growth. This completes the proof of the first two points in the lemma.

Let $a \in F$ and suppose that the fixed set of $a$ is a non-degenerate finite arc $\beta$. The unique root property in $F$ implies that, for every $b \in F$, we have $N_M(\langle b^k\rangle)=N_M(\langle b\rangle)$ and $C_M(b)=C_M(b^k)$. As the $F$-stabiliser of each edge is root-closed there is no loss of generality in assuming $a$ is not a proper power.

As $\beta$ consists of the edges in $T$ with $F$-stabiliser $\langle a \rangle$, an element in $M$ preserves $\beta$ if and only if it preserves $\langle a \rangle$, which shows that $\SetStab(\beta)=N_M(\langle a\rangle)$. In particular, for every edge $e \subseteq \beta$, we have $\SetStab(e) \subseteq \SetStab(\beta)$.

Fix such an edge $e \subseteq \beta$. As $\beta$ is a finite arc, any element of $\PtStab(e)$ is contained in $\PtStab(\beta)$. Hence $\SetStab(e)$ has an index-at-most-two subgroup that is contained in $\PtStab(\beta)$. Therefore, by Lemma~\ref{lemma:edge-stab}, the group $\chi(\PtStab(\beta))$ is non-trivial. Since $\PtStab(\beta)\cap F=\langle a \rangle$ and $\chi(\PtStab(\beta))$  is non-trivial, we have $\PtStab(\beta)=\langle a,t\rangle$ for some $t \in M$.

Suppose $g \in N_M(\langle a\rangle)$ acts non-trivially on $\beta$, so there exist distinct vertices $v$ and $w$ swapped by $g$. As we have chosen the $z_v$ so that $\chi(z_v)>0$, the element $g$ conjugates $z_v$ to $z_w$ (and vice versa). As $z_v$ and $z_w$ are conjugate, we have $\chi(z_v)=\chi(z_w)$. As $z_v, z_w \in C_M(a)$, they preserve $\beta$ so $z_v^2$ and $z_w^2$ fix $\beta$ (as $\beta$ is an arc, the only way an element of $\SetStab(\beta)$ can act non-trivially is as an inversion). We can then find $i,j,k$ such that $z_v^2=a^it^k$ and $z_w^2=a^jt^k$, where $i-j \neq 0$ (as $\langle z_v \rangle \cap \langle z_w \rangle = \{1\}$). Then $$ga^{i-j}g^{-1}=gz_v^2z_w^{-2}g^{-1}=z_w^2z_v^{-2}=a^{j-i},$$
which shows that $g \in N_M(\langle a\rangle)\setminus C_M(a)$ and $C_M(a) \subset \PtStab(\beta)$. In particular for each $v \in \beta$ the element $z_v$ belongs to $\PtStab(\beta)$.

It remains to show that $\PtStab(\beta)\subset C_M(a)$. It is sufficient to show that $\PtStab(\beta)\simeq \mathbb{Z}^2$, as if $N_M(\langle a\rangle) \neq C_M(a)$ then $N_M(\langle a\rangle)$ is isomorphic to the fundamental group of a Klein bottle (and in particular is not abelian). Similarly to above, as $\beta$ contains at least two vertices and $\langle z_v \rangle \cap \langle z_w \rangle = \{1\}$ for distinct vertices in $T$, we can find $v$ such that $z_v=a^it^k$ with $i \neq 0$. As $t$ fixes $\beta$ we have $tz_vt^{-1}=z_v$, so $ta^it^{-1}=a^i$ and therefore $tat^{-1}=a$. Therefore $\PtStab(\beta)\simeq \mathbb{Z}^2$ and is equal to $C_M(a)$.
\end{proof}

\subsection{Slow subgroups, and viewing the limiting tree as a JSJ tree}\label{sec:slow}

In this section, we introduce the set $\Slow^{\max}(M)$ of maximal slow-growing subgroups of $M$, show that it is finite and $\Out(M)$-invariant, and that $T$ can be reinterpreted as a slender JSJ splitting of $M$ relative to $\Slow^{\max}(M)$. In particular, even if $T$ is not $\Out(M)$-invariant, its deformation space is, which will allow us to obtain an $\Out(M)$-invariant tree in the next section by passing to a tree of cylinders.  

\begin{defn}[Slow subgroups]
Given a polynomially-growing free-by-cyclic group $M$ of degree $d\geq 1$, we let $\Slow(M)$ be the set of non-cyclic free-by-cyclic subgroups $H\subseteq M$ of degree at most $d-1$, and $\Slow^{\max}(M)$ be the set of maximal (for inclusion) subgroups in $\Slow(M)$. 
\end{defn}

Note that non-cyclic slender subgroups of $M$ are automatically free-by-cyclic groups with polynomial growth of degree $0$, so they are automatically in $\Slow(M)$. Our next lemma establishes some basic properties of the sets $\Slow(M)$ and $\Slow^{\max}(M)$, reminiscent of Levitt's work on polynomially-growing subgroups \cite[Proposition~1.4]{Levittgrowth}, or \cite[Theorem~1.4]{CHHL}.

\begin{lemma}
\label{lemma:slow}Let $M$ be a polynomially-growing free-by-cyclic group with degree $d\ge 1$. Then:
\begin{enumerate}
\item\label{slow-1} $\Slow^{\max}(M)$ consists of finitely many conjugacy classes of subgroups of $M$;
\item\label{slow-2} $\Slow(M)$ and $\Slow^{\max}(M)$ are $\Aut(M)$-invariant;
\item\label{slow-3} every subgroup in $\Slow(M)$ is contained in a subgroup in $\Slow^{\max}(M)$.
\end{enumerate}
\end{lemma}

\begin{proof}
Let $\Phi\in\Out(F)$ be a monodromy of $M$, and let $T$ be the limiting tree of $\Phi$.

\textbf{Assertion 1.} By Proposition~\ref{prop:limiting_tree_0}, every subgroup in $\Slow^{\max}(M)$ is a point stabiliser in $T$. In particular, $\Slow^{\max}(M)$ consists of finitely many conjugacy classes. 

\textbf{Assertion~2.} Let $H \in \Slow(M)$ and let $\alpha \in \Aut(M)$. Then $\alpha(H)$ is a free-by-cyclic subgroup of $M$. Since $H$ is non-cyclic, so is $\alpha(H)$. As $\alpha|_H \colon H \to \alpha(H)$ is a quasi-isometry, the group $\alpha(H)$ is polynomially growing of degree at most $d-1$ by the work of Macura~\cite{Macura}, so $\alpha(H)\in\Slow(M)$. 

If $H$ belongs to $\Slow^{\max}(M)$, we now prove that $\alpha(H)$ belongs to $\Slow^{\max}(M)$. Indeed, if $\alpha(H)\subseteq K$ and $K\in\Slow(M)$, then $H\subseteq\alpha^{-1}(K)$, and as before $\alpha^{-1}(K)\in\Slow(M)$. By maximality of $H$, we have $H=\alpha^{-1}(K)$ and hence $\alpha(H)=K$. Altogether, the group $\alpha(H)$ belongs to $\Slow^{\max}(M)$, showing that $\Slow^{\max}(M)$ is $\Aut(M)$-invariant.

\textbf{Assertion~3.} Let $H\in\Slow(M)$. By Proposition~\ref{prop:limiting_tree_0}(\ref{limit1,5}), $H$ is elliptic in $T$. 

$\drsh$ If $H$ is not slender, then $H$ fixes a unique vertex $v$ of $T$, and by Proposition~\ref{prop:limiting_tree_0}(\ref{limit1}), $M_v$ belongs to $\Slow^{\max}(M)$. As $H\subseteq M_v$ we are done.

$\drsh$ We now assume that $H$ is slender. 
Let $\hat H$ be a maximal slender subgroup that contains $H$.
Since $\hat H\cap F$ contains $H\cap F$ as a finite-index subgroup, it is elliptic in $T$. By acylindricity, and using that $H\cap F\neq\{1\}$, the fixed point set of $\hat H\cap F$ is bounded, so $\hat H$ fixes its circumcenter. Since $H$ (and hence $\hat H$) is non-cyclic, Proposition~\ref{prop:limiting_tree_0}(\ref{limit1}) ensures that $\hat H\in\Slow(M)$. If $\hat H$ belongs to $\Slow^{\max}(M)$, we are done. Otherwise, there exists $\hat H'\in\Slow(M)$ such that $\hat H\subsetneq\hat H'$, and by maximality of $\hat H$ among slender subgroups, the group $\hat H'$ is non-slender. Proposition~\ref{prop:limiting_tree_0}(\ref{limit1,5}) ensures that $\hat H'$ fixes a point $v$ in $T$, and Proposition~\ref{prop:limiting_tree_0}(\ref{limit1}) ensures that the conjugacy class of $M_v$ belongs to $\Slow^{\max}(M)$. As $H\subseteq M_v$, the proof is complete. 
\end{proof}

\begin{prop}
\label{prop:elliptic-characterisation}Let $M=F\rtimes_\Phi\mathbb{Z}$ be a polynomially-growing free-by-cyclic group with degree $d\ge 1$, and let $T$ be the limiting tree of $\Phi$.
    
A subgroup $H\subseteq M$ is elliptic in $T$ if and only if $H$ has a finite-index subgroup contained in a subgroup in $\Slow(M)$.
\end{prop}

\begin{proof}
Subgroups in $\Slow(M)$ are elliptic in $T$ by Proposition~\ref{prop:limiting_tree_0}(\ref{limit1,5}), and therefore so is any subgroup that is virtually contained in a subgroup in $\Slow(M)$.

Conversely, assume that $H\subseteq M$ fixes a point $v$ in $T$. By Proposition~\ref{prop:limiting_tree_0}(\ref{limit1}), either $M_v\in\Slow(M)$ or $M_v\simeq\mathbb{Z}$. In the former case we are done, so let us assume that $M_v\simeq\mathbb{Z}$.  

Consider a path from $v$ to a vertex $w$ such that $M_w$ is non-cyclic (in particular $M_w\in\Slow(M)$), while for every vertex $u\neq w$ along the path, one has $M_u\simeq\mathbb{Z}$. All stabilisers (of edges or vertices) along this path (except $M_w$) are commensurable, so $M_v$ has a finite-index subgroup contained in $M_w$. Altogether $H$ has a finite-index subgroup contained in $M_w$ and $M_w\in\Slow(M)$.
\end{proof}

The above proposition ensures that the deformation space of the limiting tree $T$ is $\Aut(M)$-invariant, in fact the following is true.

\begin{cor}
\label{cor:limit-is-jsj}Let $M=F\rtimes_\Phi\mathbb{Z}$ be a polynomially-growing free-by-cyclic group with degree $d\ge 1$, and let $T$ be the limiting tree of $\Phi$.
    
    Then $T$ is a slender JSJ tree relative to $\Slow^{\max}(M)$ (or equivalently, relative to $\Slow(M)$).
\end{cor}

\begin{proof}
In view of the above proposition, every point stabiliser $M_v$ is elliptic in any slender splitting of $M$ relative to $\Slow(M)$ (or equivalently $\Slow^{\max}(M)$ by Lemma~\ref{lemma:slow}(\ref{slow-3})). This is enough to conclude in view of Remark~\ref{rk:universally-vertex-elliptic}. 
\end{proof}

\subsection{The collapsed tree of cylinders of the limiting tree}\label{sec:collapsed}

Recall that given an action of a group $G$ on a simplicial tree $T$, an \emph{inversion point} is a point $x\in T$ which is the midpoint of an edge, whose stabiliser contains an element that permutes the two incident half-edges. From now on, in order to stick to the setting of Section~\ref{sec:arboriculture}, we will always add all inversion points to the simplicial structure of the trees we consider (notably, the limiting tree $T$) to ensure that the action has no inversion (i.e.\ if some element stabilises an edge, then it fixes it pointwise).

For our next definition, the reader is refered to Section~\ref{sec:tree-of-cyl} for background on (collapsed) trees of cylinders.

\begin{defn}[The (collapsed) tree of cylinders of the limiting tree]
\label{de:collapsed-tree-of-cyl}Let $M=F\rtimes_\Phi\mathbb{Z}$ be a free-by-cyclic group with polynomial growth $d\ge 1$, and let $T$ be the limiting tree of $\Phi$. We define the \emph{tree of cylinders} $T_c$ and the \emph{collapsed tree of cylinders} $T_c^*$ in the following way:
\begin{itemize}
    \item If $d=1$, then $T_c$ is the tree of cylinders of $T$ for the relation of commensurability on the class $\mathcal{E}$ of subgroups of $M$ that have a finite-index subgroup isomorphic to $\mathbb{Z}^2$, and $T_c^*$ is its collapsed tree of cylinders.
    \item If $d\ge 2$, then $T_c$ is the tree of cylinders of $T$ for the relation of commensurability on the class $\mathcal{E}$ of infinite cyclic subgroups of $M$, and $T_c^*$ is its collapsed tree of cylinders.
\end{itemize}
\end{defn}

\begin{rk}
Note that in both cases, edge stabilisers of $T$ are indeed in the class $\mathcal{E}$ considered, as established in Proposition~\ref{prop:limiting_tree_0}. See also Remark~\ref{rk:KB} for a more detailed discussion of our choice of class $\mathcal{E}$ in this definition. 
\end{rk}

We now collect properties of the tree $T_c^*$ which will be crucial in the sequel.

\begin{thm}
\label{theo:collapsed-tree-of-cylinders}Let $M=F\rtimes_\Phi\mathbb{Z}$ be a free-by-cyclic group with polynomial growth of degree $d\ge 1$, and let $T$ be the limiting tree of $\Phi$. Then $T_c^*$ satisfies the following properties:
\begin{enumerate}
\item\label{properties1} $T_c^*$ has the same elliptic subgroups as $T$; in particular, a subgroup $H\subseteq M$ is elliptic in $T_c^*$ if and only if $H$ has a finite-index subgroup contained in a subgroup in $\Slow(M)$.
\item\label{properties2} $T_c^*$ is $\Out(M)$-invariant.
    \item\label{properties3} If $g \in F$ then the conjugacy class of $g$ has polynomial growth of degree $d$ under the monodromy if and only if $g$ is hyperbolic in $T_c^*$.
    \item \label{item:smaller_growth_vertex_groups} Every vertex group in $T_c^*$ is either slender or belongs to $\Slow^{\max}(M)$.
    \item\label{properties5} As an $M$-tree, $T_c^*$ is acylindrical.
    \item \label{item:limiting_tree_as_an_F_tree} As an $F$-tree, $T_c^*$ is a free splitting if $d \geq 2$, and is a cyclic splitting if $d=1$.
\end{enumerate}
\end{thm}

\begin{proof}
(1) Since $T$ is acylindrical as an $M$-tree (Corollary~\ref{cor:acylindrical}), the cylinders (in the sense of having commensurable stabilisers) are bounded. It thus follows from Lemma~\ref{lemma:tree-of-cylinders}(\ref{collapse1}) that $T$ and $T_c^*$ have the same elliptic subgroups. Their characterisation now follows from Proposition~\ref{prop:elliptic-characterisation}.

(2) We now prove that $T_c^*$ is $\Out(M)$-invariant. By Corollary~\ref{cor:limit-is-jsj}, the tree $T$ is a slender JSJ splitting of $M$ relative to $\Slow(M)$. Recall in addition that $\Slow(M)$ is $\Aut(M)$-invariant (Lemma~\ref{lemma:slow}). Therefore, for every $\Phi\in\Out(M)$, the tree $\Phi\cdot T$ is a slender JSJ splitting of $M$ relative to $\Slow(M)$, so it is in the same deformation space as $T$. It follows from \cite[Corollary~4.10]{GL-cylinders} that $\Phi\cdot T_c^*=(\Phi\cdot T_c)^*=T_c^*$, showing that $T_c^*$ is $\Out(M)$-invariant.

(3) Since by the first assertion $T_c^*$ has the same elliptic subgroups as $T$, Assertion~(\ref{properties3}) follows from Theorem~\ref{thm:limiting_tree_existence}.

(4) By Lemma~\ref{lemma:tree-of-cylinders}(\ref{collapse2}), every vertex stabiliser in $T_c^*$ is either a vertex stabiliser in $T$ or slender. Assertion~(\ref{item:smaller_growth_vertex_groups}) thus follows from Proposition~\ref{prop:limiting_tree_0}.

(5) Recall from Corollary~\ref{cor:acylindrical} that $T$ is acylindrical as an $M$-tree. Since by the first assertion $T_c^*$ has the same elliptic subgroups as $T$, acylindricity is inherited (see \cite[Section~4]{GL-defo}).

(6) By Lemma~\ref{lemma:tree-of-cylinders}(\ref{collapse3}), every edge of $T_c^*$ has stabiliser commensurable to an edge stabiliser in $T$. Assertion~\ref{item:limiting_tree_as_an_F_tree} thus follows from Theorem~\ref{thm:limiting_tree_existence}.
\end{proof}

We conclude with the analogue of Lemma~\ref{lem:linear-image-edge-finite} for the collapsed tree of cylinders in the linear case, as well as studying the possible automorphisms of edge groups induced by $\Out(M)$. Denote by $K$ the fundamental group of a Klein bottle.

\begin{lem}
\label{lem:linear_T_c_edges}Let $M=F\rtimes_\Phi\mathbb{Z}$ be a linearly-growing free-by-cyclic group, $T$ the limiting tree of $\Phi$, $T_c$ its tree of cylinders and $T_c^*$ its collapsed tree of cylinders. 
\begin{enumerate}
\item Each cylinder $C$ in $T$ is a finite arc, with stabiliser isomorphic to either $\mathbb{Z}^2$ or $K$. In particular, $T_c=T_c^*$.
\item Every point $v \in T$ contained in at least two cylinders has stabiliser $M_v$ equal to a periodic free-by-cyclic group with non-abelian fibre.
\item At most one edge $\varepsilon=(v,C)$ adjacent to $C \in T_c$ has $M$-stabiliser isomorphic to $K$, otherwise $M_\varepsilon\simeq \mathbb{Z}^2$. Each edge group $M_\varepsilon$ is index at most two in the adjacent vertex group $M_C$. 
\item Under the $\Aut(M)$-action on $T_c$, the image of $\Stab_{\Aut(M)}(C)$ in $\Aut(M_C)$ is virtually cyclic (possibly finite), and therefore the same is true for the image of $\Stab_{\Aut(M)}(\varepsilon)$ in $\Aut(M_\varepsilon)$ for every edge $\varepsilon$ in $T_c^*$. 
\end{enumerate}
\end{lem}

\begin{proof}
\textbf{Assertion~1.} Let $C$ be a cylinder for commensurability in $T$. If $e$ is an edge in $C$ and $F_e=\langle a\rangle$, then $F_{e'}$ is commensurable to $F_e$ for every other edge $e'$ in $C$. As edge stabilisers for the $F$-action on $T$ are root-closed, it follows that $F_{e'}=\langle a\rangle $ for every edge $e'$ in $C$, so $C$ is contained in the finite arc $\beta$ fixed by $a$. In Lemma~\ref{lem:linear-image-edge-finite}(\ref{stab-3}) we showed that the pointwise stabiliser of each edge in $\beta$ is the centraliser $C_M(a) \simeq \ZZ^2$. This implies that every edge in $\beta$ lies in $C$ and $C=\beta$. 

By Lemma~\ref{lem:linear-image-edge-finite}(\ref{stab-3}) again, the setwise stabiliser of $\beta$ is the normaliser $N_M(\langle a\rangle)$. Then $M_C=N_M(\langle a\rangle)$, which we showed in Lemma~\ref{lem:linear-image-edge-finite}(\ref{stab-3}) is isomorphic to $\mathbb{Z}^2$ (when $C_M(a)=N_M(\langle a\rangle)$) or $K$ (when $C_M(a)$ is a proper subgroup of $N_M(\langle a\rangle)$). In particular, as every edge stabiliser of $T_c$ is contained in a cylinder stabiliser, $T_c$ is already an $M$-slender splitting and $T_c=T_c^*$. This proves the first assertion.

\smallskip

\textbf{Assertion~2.} Notice that we have shown in the first assertion that cylinders for commensurability for the $M$-action on $T$ are exactly the cylinders for commensurability for the $F$-action on $T$. It follows that if two cylinders intersect they do so at a vertex $v$ of valence at least three. We have seen above (e.g.\ Lemma~\ref{lem:linear-image-edge-finite}(\ref{stab-1})) that $M_v$ is a periodic free-by-cyclic group with non-abelian fibre. The second assertion follows.

\smallskip

\textbf{Assertion~3.} Note that for any edge $\varepsilon=(v,C)$ in $T_c$ we have $M_\varepsilon=M_C \cap M_v$. As $M_C=N_M(\langle a\rangle)$ acts on $C$ by inversions under the quotient $N_M(\langle a\rangle)/C_M(a)\simeq \mathbb{Z}/2\mathbb{Z}$ (Lemma~\ref{lem:linear-image-edge-finite}(3)), the midpoint of $C$ is fixed by $N_M(\langle a\rangle)$ and the stabiliser of every other point $v$ under the $M_C$-action is isomorphic to $C_M(a) \simeq \mathbb{Z}^2$. This shows that $M_{\varepsilon}=C_M(a)\simeq \mathbb{Z}^2$ unless $N_M(\langle a\rangle)=K$ and the vertex $v$ at the other end of $\varepsilon$ is the midpoint of $C$. As $M_\varepsilon$ is equal to either $N_M(\langle a\rangle)$ or $C_M(a)$, it is of index at most 2 in $M_C$.

\textbf{Assertion~4.} Note that if $M_C\simeq K$ the group $\Aut(K)$ is virtually cyclic (see e.g.\ \cite[Lemma~2.1]{GW}) and the result is immediate. We can therefore assume that $M_C\simeq \mathbb{Z}^2$, and suppose $\alpha \in \Aut(M)$ fixes $C\in T_c$. As $C$ is a finite arc (viewed as a subset of $T$), the vertex determined by $C$ in $T_c^*$ has finite valence, so some finite power $\alpha^k$ fixes all edges adjacent to $C$ in $T_c^*$. As the $M$-action on $T_c^*$ is minimal, there are at least two distinct adjacent edges which we call $\varepsilon=(v,C)$ and $\varepsilon'=(w,C)$. As $\alpha^k$ fixes the vertices $v$ and $w$, it preserves $M_v$ and $M_w$ and therefore also preserves their centres $\langle z_v \rangle$ and $\langle z_w\rangle$ given by Lemma~\ref{lem:linear-image-edge-finite}(\ref{stab-1}). These belong to $M_C$ and have trivial intersection by  Lemma~\ref{lem:linear-image-edge-finite}(\ref{stab-2}). Any automorphism preserving two distinct cyclic subgroups of $\mathbb{Z}^2$ has finite order (e.g.\ by considering Smith normal forms), so a further power of $\alpha$ is trivial. We have shown that the image of  $\Stab_{\Aut(M)}(C)$ in $\Aut(M_C)={\mathrm{GL}}(2,\mathbb{Z})$ is torsion, so must be finite.
\end{proof}

\begin{cor}
\label{coro:image-edge-finite}Let $M$ be a polynomially-growing free-by-cyclic group with infinite order monodromy and let $T_c^*$ be its collapsed tree of cylinders. For every edge $\varepsilon \in ET_c^*$, the image of $\Stab_{\Aut(M)}(\varepsilon)$ in $\Aut(M_\varepsilon)$ given by the restriction homomorphism is virtually cyclic.
\end{cor}
\begin{proof}
Combine Lemma~\ref{lem:linear_T_c_edges} in the linear growth case with the fact that in higher growths the edge stabilisers $M_\varepsilon$ are cyclic and so any subgroup of $\Aut(M_\varepsilon)$ is finite. 
\end{proof}

\subsection{Further discussion}\label{sec:bonus}

We conclude this section by establishing a few extra properties of the tree $T_c^*$ and comparing our work with the existing literature. The contents of this section are not needed for the rest of the article. 

 Since the (collapsed) tree of cylinders of $T$ depends only on the deformation space containing $T$ and not the precise tree \cite[Corollary 4.10]{GL-cylinders}, it is reasonable to describe the collapsed tree of cylinders of the limiting tree as a ``canonical'' (not merely invariant) $M$-tree. In particular, the collapsed tree of cylinders is independent of the choice of fibre (and orientation) which is used to build a limiting tree. In fact it is even a commensurability invariant, in the following sense. Recall that the degree of a polynomial free-by-cyclic group is a commensurability (even quasi-isometry) invariant \cite{Macura}.

\begin{prop}
\label{prop:commensurability}Suppose $G$ is a finite-index subgroup of a polynomially-growing free-by-cyclic group $M$. Let $T_1$ be the collapsed tree of cylinders of a limiting tree of $G$ and let $T_2$ be the collapsed tree of cylinders of a limiting tree of $M$ with the action restricted to $G$. Then as $G$-trees, $T_1$ and $T_2$ are equivariantly homeomorphic.
\end{prop}

\begin{proof} 
We will first prove that the $G$-trees $T_1$ and $T_2$ are in the same deformation space, i.e.\ they have the same elliptic subgroups. By Theorem~\ref{theo:collapsed-tree-of-cylinders}(\ref{properties1}), applied to the $M$-action on $T_2$, a subgroup $H\subseteq M$ is elliptic for its action on $T_2$ if and only if $H$ is virtually contained in a subgroup in $\Slow(M)$. The work of Macura \cite{Macura} implies that the polynomial growth rate of a free-by-cyclic group is a commensurability invariant (see Theorem~\ref{thm:Macura}), therefore a subgroup of $M$ is contained in $\Slow(M)$ if and only if its intersection with $G$ is contained in $\Slow(G)$. It follows that a subgroup $K$ of $G$ is elliptic for its action on $T_2$ if and only if $K$ is virtually contained in a subgroup in $\Slow(G)$. Using Theorem~\ref{theo:collapsed-tree-of-cylinders}(\ref{properties1}) again, this time for the $G$-action on $T_1$, it follows that $T_1$ and $T_2$ belong to the same deformation space of $G$-trees.

Now \cite[Corollary~4.10]{GL-cylinders} implies that $T_1$ and $T_2$, as $G$-trees, have the same collapsed tree of cylinders. By \cite[Corollary~5.8]{GL-cylinders}, $T_1$ is its own collapsed tree of cylinders for the relation of commensurability, and $T_2$, as an $M$-tree, is its own collapsed tree of cylinders for the relation of commensurability. Since for this relation cylinders do not change when passing to a finite-index subgroup, it follows that $T_2$, as a $G$-tree, is also its own collapsed tree of cylinders. Therefore $T_1=T_2$, as desired.
\end{proof}

\begin{rk}\label{rk:KB}
In the case $d=1$, edge stabilisers of the limiting tree $T$ of the monodromy are isomorphic to $\mathbb{Z}^2$ (Proposition~\ref{prop:limiting_tree_0}). This allows us to also view $T$ as a JSJ splitting relative to $\Slow(M)$ with abelian edge stabilisers (instead of slender). When passing to a tree of cylinders (in Definition~\ref{de:collapsed-tree-of-cyl}), we could thus alternatively have chosen to define $T_c$ and $T_c^*$ by using the relation of commensurability on the class of subgroups of $M$ isomorphic to $\mathbb{Z}^2$, instead of virtually abelian. Doing so, the resulting collapsed tree of cylinders has abelian edge stabilisers instead of slender ones. By Lemma~\ref{lem:linear_T_c_edges}(3), it is obtained from the slender collapsed tree of cylinders (i.e.\ our $T_c^*$) by collapsing at most one edge $\varepsilon=(v,C)$ for every cylinder vertex $C$. If such an edge $\varepsilon$ exists, then $M_\varepsilon=M_C$, which implies that the abelian collapsed tree of cylinders and the slender collapsed tree of cylinders are in the same deformation space. 

For the purposes of the rest of the paper, working with the abelian collapsed tree of cylinders instead of the slender one would have been possible. However, one of the advantages of working with the slender collapsed tree of cylinders is that it is a commensurability invariant (Proposition~\ref{prop:commensurability}), which is not the case of the abelian collapsed tree of cylinders.
\end{rk}

\begin{rk}
We finish this section by commenting upon the relationship of our work with other constructions from the literature.

Using limiting arguments in the study of automorphisms of free groups finds its roots in the celebrated Bestvina--Paulin argument. As mentioned above, our limiting tree $T$ appeared for $d=1$ in work of Cohen--Lustig \cite{CL95} and for $d\ge 2$ in work of Bestvina--Feighn--Handel \cite{BFHTits2}. More recently, these trees and some variations have been significantly exploited in the context of free-by-cyclic groups, see below. To our knowledge, using the tools of deformation spaces and (collapsed) trees of cylinders to make the trees $\Out(M)$-invariant is new in this context.

More specifically, in the case $d=1$ the tree $T_c^*$ agrees with the one studied in \cite{AMFbyZ,DT24}, and many of the conclusions can be extracted from the combination of those papers. Full $\Out(M)$-invariance in the non-UPG case is new, as is our more direct argument for acylindricity. 

In the superlinear case, we believe that in fact $T_c^*=T$. The tree $T$ has been used in the literature on free-by-cyclic groups many times before (e.g. \cite{AHKtorsionhomology,KV24}). Verifying that $T_c^*=T$ needs a more careful analysis of the train-track map giving rise to the limiting tree, after which one could conclude (following \cite{KV24}) that in the superlinear case the limiting tree must be 2-acylindrical, every edge is adjacent to (not just at finite distance from) a vertex stabilised by an element of $\Slow(M)$, and as a result $T$ is its own collapsed tree of cylinders. Even so, the fact that $T$ is an $M$-tree (not just an $M_0$-tree for some finite index subgroup) and the $\Out(M)$-invariance (even in the UPG case) is new in all but some low-complexity cases. 
\end{rk}

\section{A transverse covering construction}\label{sec:tree-of-cylinders}

In this section we work in a more general setting with a finitely generated group $G$ (to be assumed torsion free later), and a finite set $\calh$ of conjugacy classes of finitely generated subgroups of $G$. We let $T$ be a simplicial tree equipped with a minimal action of $G$. 

The goal of the section is the following: starting from the action of $G$ on $T$ (soon to be assumed acylindrical), we will describe a procedure to construct a new action of $G$ on a tree $U$ in which the subgroups from the set $\calh$ become elliptic.  The new tree $U$ is $\Out_T(G,\calh^{(t)})$-invariant, and we have enough control of the new vertex and edge stabilisers so as to get some understanding of $\Out_T(G,\calh^{(t)})$, by getting a fine description of its action on vertex and edge groups of the tree $U$. 

The construction of $U$ proceeds as follows: we look at the set of subtrees $T_H$ where $H$ varies among subgroups of $G$ whose conjugacy class belongs to $\calh$ and which are not elliptic in $T$. We then glue this collection of trees together along non-degenerate overlaps to form a  transverse covering of $T$. The tree $U$ is the skeleton of this transverse covering. Its key properties are summarised in Theorem~\ref{thm:properties-tree-cyl}. 

This construction does not require the acting group $G$ to be free-by-cyclic; we specify our results to this particular setting in Section~\ref{sec:tree-of-cyl-f-by-z}.  
 
\subsection{Transverse coverings from minimal invariant subtrees}

Let $\calh^\infty_T\subseteq\calh$ be the subset consisting of all conjugacy classes $[H] \in \calh$ such that $H$ does not fix any point in $T$. Since all subgroups whose conjugacy class belongs to $\calh$ were assumed finitely generated, for every subgroup $H\subseteq G$ whose conjugacy class belongs to $\calh^\infty_T$, there exists a minimal $H$-invariant subtree $T_H\subseteq T$, which is also the union of all axes of elements of $H$ acting hyperbolically on $T$, see \cite[Proposition~3.1]{CM} and the paragraph below it.

Let $\sim$ be the equivalence relation on $ET$ generated by: $e\sim_0 e'$ if there exists a subgroup $H\subseteq G$ with $[H]\in\calh^\infty_T$, such that $e$ and $e'$ both belong to the minimal subtree $T_H$. The equivalence relation $\sim$ is $G$-invariant, and we claim that its classes are subtrees. Indeed, if two distinct edges $e,e'$ are $\sim$-equivalent, there exist edges $e=e_1,\dots,e_k=e'$ such that for every $i\in\{1,\dots,k-1\}$, the edges $e_i,e_{i+1}$ both belong to some subtree $T_{H_i}$ with $[H_i]\in\calh^\infty_T$. Then every edge on a geodesic joining $e_i$ to $e_{i+1}$ belongs to $T_{H_i}$ and is therefore equivalent to $e_i$ and $e_{i+1}$. It follows that there is a path joining $e$ to $e'$ consisting of equivalent edges, thereby proving our claim.

Thus, the set $\caly_\calh$ consisting of all $\sim$-equivalence classes is a transverse covering of $T$. We call it the \emph{$\calh$-transverse covering} of $T$.

\begin{rk}\label{rk:covering}
If an edge $e\in ET$ is not contained in any subtree $T_H$ with $[H]\in\calh^\infty_T$, then its equivalence class is reduced to itself. On the other hand, if $e\subseteq T_H$, then the $\sim$-equivalence class of $e$ contains all edges in $T_H$, and is therefore infinite.
\end{rk}

We denote by $U_{T,\calh}$ the skeleton of the transverse covering $\caly_\calh$ of $T$.
As recalled in Section~\ref{sec:skeleta}, it is a minimal $G$-splitting which comes equipped with a $G$-invariant bipartition $VU_{T,\calh}=V_0 \amalg V_1$ of its vertex set, where $V_0$ consists of vertices of $T$ that do not belong to a unique subtree in $\caly_\calh$, and $V_1$ consists of subtrees in $\caly_\calh$. Let $V_1^\infty\subseteq V_1$ be the subset consisting of all infinite subtrees in $\caly_\calh$. 
\begin{ex}
Let $F$ be a finitely generated free group, let $g \in F$ and let $T$ be a Cayley graph of $F$ with respect to some basis of $F$. Let $\calh=\{[\langle g \rangle]\}$. One way to define the \emph{Whitehead graph} $\mathrm{Wh}(g,T)$ of $g$ with respect to $T$ is as follows: there is one vertex for every edge of $T$ adjacent to $1$, and there is an edge between two such vertices if there exists $h \in F$ such that $h\mathrm{Ax}(g)$ contains both edges. We chose the identity vertex in $T$, but one can similarly define a Whitehead graph at every vertex (the graphs at each vertex are isomorphic as $F$ acts transitively on the vertex set of $T$). 

Suppose that $h_1 \mathrm{Ax}(g)$ and $h_2 \mathrm{Ax}(g)$ are two axes in $T$ which contain $1$. Then $h_1 \mathrm{Ax}(g)$ and $h_2 \mathrm{Ax}(g)$ overlap on some edge in $T$ if and only if the edges they define in $\mathrm{Wh}(g,T)$  have a common endpoint. Therefore, the connected components of $\mathrm{Wh}(g,T)$ are in one-to-one correspondence with the subtrees in $\caly_\calh$ which contain $1$. One can use this to show that the tree $U_{T,\calh}$ is non-degenerate (i.e.\ not reduced to one point) if and only if the Whitehead graph $\mathrm{Wh}(g,T)$ is disconnected. Moreover, if $\mathrm{Wh}(g,T)$ is disconnected, then the vertex $1$ of $T$ defines a vertex in $U_{T,\calh}$ and the edges of $U_{T,\calh}$ adjacent to $1$ are in one-to-one correspondence with the connected components of $\mathrm{Wh}(g,T)$. 

Since vertex stabilisers in $T$ are trivial, edge stabilisers in $U_{T,\calh}$ are trivial, that is, $U_{T,\calh}$ is a free splitting where $g$ is elliptic. Therefore, the fact that $U_{T,\calh}$ is non-degenerate implies that $g$ is contained in a proper free factor of $F$. However, the converse is not true: if an element is contained in a proper free factor then it has a disconnected Whitehead graph for \emph{some} basis, but not for every basis.
\end{ex}

Recall that $\Aut_T(G,\calh)$ is the subgroup of $\Aut(G,\calh)$ consisting of automorphisms that preserve $T$. The group $\Aut_T(G,\calh)$ preserves the transverse covering $\caly_\calh$, so the tree $U_{T,\calh}$ is $\Aut_T(G,\calh)$-invariant. More explicitly, every automorphism $\alpha\in\Aut_T(G,\calh)$, when viewed as an isometry of $T$, sends every subtree $Y$ in $\caly_\calh$ to a subtree $Y'=\alpha(Y)$ in $\caly_\calh$; in addition, if $x\in Y$, then $\alpha(x)\in Y'$. The assignment $x\mapsto\alpha(x)$ and $Y\mapsto \alpha(Y)$ is then an $\alpha$-equivariant isometry of $U_{T,\calh}$. This yields an inclusion $\Aut_T(G,\calh)\subseteq\Aut_{U_{T,\calh}}(G,\calh)$, and in the quotient $\Out_T(G,\calh)\subseteq\Out_{U_{T,\calh}}(G,\calh)$.

If $Y\in\caly_\calh$ is infinite, then it contains some $T_H$, so there exists an element of $G_Y$ acting hyperbolically on $Y$. This implies that the minimal subtree of $G_Y$ in $T$ exists (and is contained in $Y$); the next two lemmas show that the minimal subtree of $G_Y$ is equal to $Y$. 

\begin{lemma}
\label{lem:H_i-in-stab}Let $Y \in \caly_\calh$, and let $H\subseteq G$ be a subgroup such that $[H] \in \calh^\infty_T$. 

Then $H\subseteq G_Y$ if and only if $T_H \subseteq Y$. 
\end{lemma}

\begin{proof}
If $H\subseteq G_Y$, then $H$ preserves $Y$, so by minimality of $T_H$, we have $T_H \subseteq Y$. 

Conversely, if $T_H\subseteq Y$, then for every $h\in H$, the intersection $hY\cap Y$ contains an edge. Since $\caly_\calh$ is $G$-invariant, we have $hY\in\caly_\calh$, and it follows from the definition of a transverse covering that $hY=Y$. Hence $h\in G_Y$. This shows that $H\subseteq G_Y$.
\end{proof}

\begin{lemma}
\label{lem:minsubtree}Let $Y \in \caly_\calh$ be infinite. Then $Y$ is the minimal subtree of $G_Y$ in $T$.
\end{lemma}

\begin{proof}
Let $T_{G_Y}$ be the minimal subtree of $G_Y$ in $T$. Since $Y$ is $G_Y$-invariant, we have $T_{G_Y}\subseteq Y$. We now prove the reverse inclusion.

Let $e\subseteq Y$ be an edge. By Remark~\ref{rk:covering}, the edge $e$ is contained in the minimal subtree $T_H$ of some subgroup $H\subseteq G$ with $[H] \in \calh^\infty_T$. By construction of the transverse covering, we have $T_H\subseteq Y$, and thus $H\subseteq G_Y$ by Lemma~\ref{lem:H_i-in-stab}. Therefore $T_{G_Y}$ contains the minimal subtree $T_H$ of $H$, so $e\subseteq T_{G_Y}$. It follows that $Y\subseteq T_{G_Y}$, as desired.
\end{proof}

\begin{lemma}
\label{lem:stabedge}Let $Y\in \caly_\calh$ and let $e\in ET$. Suppose that $e \subseteq Y$. Then $G_e \subseteq G_Y$.
\end{lemma}

\begin{proof}
Let $g \in G_e$. Since $ge=e$, the intersection $gY \cap Y$ contains the edge $e$. Since the transverse covering $\caly_\calh$ is $G$-invariant, we have $gY\in\caly_\calh$. By definition of a transverse covering, we deduce that $gY=Y$, so $g\in G_Y$.
\end{proof}

\begin{lemma}
\label{lem:stability_conjugation}Let $Y\in \caly_\calh$. Let $H\subseteq G$ be a subgroup such that $[H] \in \calh^\infty_T$, let $g \in G$, and assume that $H,gHg^{-1} \subseteq G_Y$. Then $g \in G_Y$.
\end{lemma}

\begin{proof}
By Lemma~\ref{lem:H_i-in-stab}, the trees $T_H$ and $gT_H$ are contained in $Y$. Hence $gY \cap Y$ contains an edge, and therefore $gY=Y$.
\end{proof}

\subsection{Generating the vertex and edge groups of \texorpdfstring{$U_{T,\mathcal{H}}$}{UT,H}}

As above, $T$ is a tree equipped with a minimal action of a finitely generated group $G$, and $\calh$ is a finite set of conjugacy classes of finitely generated subgroups of $G$. Consider the $\calh$-transverse covering $\caly_\calh$ described above and the associated skeleton $U_{T,\calh}$. For $Y\in V_1^\infty$, we let $\calh_Y\subseteq\calh$ and  $\calh^\infty_Y\subseteq\calh^\infty_T$ be the subsets consisting of conjugacy classes having a representative $H\subseteq G_Y$ (equivalently by Lemma~\ref{lem:H_i-in-stab}, such that $T_H\subseteq Y$). Note that Lemma~\ref{lem:stability_conjugation} allows us to blur the distinction between $G$-conjugacy classes and $G_Y$-conjugacy classes of subgroups in $\calh_Y$ (or $\calh_Y^\infty$). We will thus also allow ourselves to view $\calh_Y$ and $\calh_Y^\infty$ as sets of $G_Y$-conjugacy classes of subgroups of $G_Y$.

Let $\tau\subseteq Y$ be a subtree, and $H\subseteq G_Y$ be a subgroup with $[H]\in\calh^\infty_Y$. We say that $\tau$ is \emph{$H$-spanning} if $H\cdot (\tau\cap T_H)=T_H$. More generally, if $X$ is a set of subgroups of $G_Y$, we say that $\tau$ is \emph{$X$-spanning} if it is $H$-spanning for every $H\in X$. Note that if $X$ is a finite set of finitely generated subgroups of $G_Y$, then as each $H \in X$ acts cocompactly on $T_H$ there always exists an $X$-spanning finite subtree $\tau\subseteq Y$. 

\begin{lemma}
\label{lemma:spanning-genset}Let $Y\in V_1^\infty$. Let $X$ be a finite set of subgroups of $G_Y$ which is a set of representatives of the conjugacy classes in $\mathcal{H}^\infty_Y$, and let $\tau\subseteq Y$ be an $X$-spanning subtree. 

Then $G_Y$ is generated by the subgroups in $X$, together with the set of all elements $h\in G$ such that $h\tau\cap\tau$ contains an edge.
\end{lemma}

\begin{proof} 
Let $S_Y$ be the set consisting of all elements in the subgroups in $X$, together with all elements $h\in G$ such that $h\tau\cap\tau$ contains an edge. We first show that for every edge $f\subseteq Y$ there exist $g \in \langle S_Y \rangle$ and an edge $e \subseteq \tau$ such that $ge=f$. We proceed by induction on the minimal length of a chain \[ e_0,e_1, \ldots, e_m=f,\] such that $e_0 \in \tau$ and $e_j\sim_0 e_{j+1}$ for every $j\in\{0,\dots,m-1\}$ (recall that $e\sim_0e'$ means that the edges belong to some common $T_H$). Such a chain exists by construction of the transverse covering. If $m=0$ then $f \subseteq \tau$ so certainly belongs to the $\langle S_Y \rangle$-orbit of $\tau$. Now suppose $f \in Y$ and there is a minimal length chain as above of length $m$. By induction there exists $g_{m-1} \in \langle S_Y \rangle$ and $e' \subseteq \tau$ such that $g_{m-1}e'=e_{m-1}$. As $e_{m-1}\sim_0f$ there exists a subgroup $H \subseteq G_Y$ whose conjugacy class belongs to $\mathcal{H}^\infty_Y$ such that $e_{m-1},f \subseteq T_H$. Let $h\in G$ be such that $H_0:=hHh^{-1}$ belongs to $X$. 
Since $he_{m-1},hf \subseteq T_{H_0}$ and $\tau$ contains a fundamental domain for the action of $H_0$ on $T_{H_0}$, there exist $h_{m-1},h_f \in H_0$ such that $h_{m-1}he_{m-1} \subseteq\tau$ and $h_fhf\subseteq \tau$. As $h_{m-1}he_{m-1}=h_{m-1}hg_{m-1}e'$ is contained in $\tau$, we have $h_{m-1}hg_{m-1} \in S_Y$. Therefore $h \in \langle S_Y \rangle$ as $h_{m-1} \in H_0$ and $g_{m-1} \in \langle S_Y \rangle$. It follows that $h_fh \in \langle S_Y \rangle$ as $h_f \in H_0$, and $f$ is in the $\langle S_Y \rangle$-orbit of an edge in $\tau$. 

The work above shows that the  $\langle S_Y\rangle$-translates of $\tau$ cover $Y$. Given that, it is routine to deduce that $S_Y$ generates $G_Y$: if $g \in G_Y$ and $e \subseteq \tau$ then $ge \subseteq Y$ so there exist $g' \in \langle S_Y\rangle$ and $e' \subseteq \tau$ such that $g'e'=ge$. Then $g^{-1}g' \tau \cap \tau$ contains an edge, so belongs to $S_Y$. Hence $g \in \langle S_Y \rangle$ as $g^{-1}g' \in S_Y$ and $g' \in \langle S_Y \rangle$.
\end{proof}

\begin{lemma}
\label{lemma:edge-stab-fg}Assume that edge stabilisers of $T$ are finitely generated. Then vertex and edge stabilisers of $U_{T,\calh}$ are finitely generated.
\end{lemma}

\begin{proof}
We first prove that the stabiliser of every vertex $v\in VU_{T,\calh}$ is finitely generated.

Suppose first that $v \in V_0$. Then $G_v$ is a vertex stabiliser in $T$. As $G$ is finitely generated and edge stabilisers of $T$ are finitely generated, so are vertex stabilisers of $T$ (see, e.g.\ \cite{HaglundWise}); in particular $G_v$ is finitely generated. 

Suppose now that $v \in V_1$ and let $Y \in \caly_\calh$ be the corresponding tree. By Remark~\ref{rk:covering}, either $Y$ is an edge or else $Y$ is infinite. 

$\drsh$ If $Y$ is an edge, then $G_Y$ is an edge stabiliser in $T$, so it is finitely generated by assumption. 

$\drsh$ If $Y$ is infinite, then $v \in V_1^\infty$. Let $X$ be a finite set of subgroups of $G_Y$ which is a set of representatives of the conjugacy classes in $\mathcal{H}^\infty_Y$, and let $\tau\subseteq Y$ be a finite $X$-spanning subtree. By Lemma~\ref{lemma:spanning-genset}, the group $G_Y$ is generated by the subgroups in $X$ (which are finitely generated), together with the set $\mathcal{X}$ of all elements $h\in G$ such that $h\tau\cap\tau$ contains an edge. The set $\mathcal{X}$ is a finite union of cosets of edge stabilisers: indeed, letting $\mathcal{E}$ be the set of pairs $(e,f)$ of edges of $\tau$ in the same $G$-orbit, and choosing an element $g_{e,f}\in G$ sending $e$ to $f$ for every $(e,f)\in \mathcal{E}$, one has \[\mathcal{X}=\bigcup_{(e,f)\in \mathcal{E}}g_{e,f}G_e.\] Since the subgroups $G_e$ are finitely generated by assumption, so is $G_Y$. 

\medskip

We now prove that the stabiliser of every edge $\varepsilon \in EU_{T,\calh}$ is finitely generated. Write $\varepsilon=(Y,x)$, where $Y\in\caly_\calh$ is a subtree and $x \in VT$. Thus, $G_\varepsilon=\Stab_{G_Y}(x)$. We now consider the action of $G_Y$ on its minimal subtree $T_{G_Y}$, which is equal to $Y$ by Lemma~\ref{lem:minsubtree}. For the action of $G_Y$ on $Y$, the stabiliser (in $G_Y$) of every edge $e\subseteq Y$ is the same as the stabiliser of $e$ for the $G$-action on $T$ (Lemma~\ref{lem:stabedge}). In particular, the action of $G_Y$ on $Y$ has finitely generated edge stabilisers, and as $G_Y$ is finitely generated, it follows that it also 
has finitely generated vertex stabilisers. So $\Stab_{G_Y}(x)$ is finitely generated, as desired. 
\end{proof}

\subsection{The action of \texorpdfstring{$\Aut_T(G,\calh)$}{AutT(G,H)} on the edge groups of \texorpdfstring{$U_{T,\calh}$}{UT,H}} 

Our next lemma shows that, in order to understand the action of $\Aut_T(G,\calh)$ on the edge groups of $U_{T,\calh}$, it suffices to understand the action of $\Aut_T(G,\calh)$ on the vertex groups of $U_{T,\calh}$. To state it, we introduce the following notation.

\begin{notation} 
\label{notation:Out_1}Recall that $\Out^0_{U_{T,\calh}}(G,\calh)$ denotes the finite-index subgroup acting trivially on the quotient graph $U_{T,\calh}/G$. We let $\Out^1_T(G,\calh):=\Out_T(G,\calh)\cap\Out^0_{U_{T,\calh}}(G,\calh)$.
\end{notation}

\begin{lemma}
\label{lem:factorisation-restriction}Let $G$ be a finitely generated group, let $T$ be a $G$-tree and let $\mathcal{H}$ be a finite set of conjugacy classes of finitely generated subgroups of $G$. 
Let $\varepsilon\in EU_{T,\calh}$ and let $Y \in VU_{T,\calh}$ be its endpoint in $V_1$. 

The restriction homomorphism $\Out_T^1(G,\calh) \to \Out(G_\varepsilon)$ factors through the restriction homomorphism $\Out_T^1(G,\calh) \to \Out(G_Y)$.
\end{lemma}

\begin{proof}
Our goal is to prove that if $\Phi\in\Out_T^1(G,\calh)$ acts trivially on $G_Y$ (i.e.\ it belongs to the kernel $K_Y$ of the restriction homomorphism to $\Out(G_Y)$), then it also acts trivially on $G_\varepsilon$.

The edge $\varepsilon$ corresponds to a pair $(Y,w)$ where $Y\in \caly_\calh$ and $w\in VT$. Since $\Phi\in\Out^0_{U_{T,\calh}}(G,\calh)$, we can choose a representative $\varphi\in\Aut(G)$ that fixes $\varepsilon$. In particular $\varphi(G_Y)=G_Y$. Since $\Phi\in K_Y$, the restriction $\varphi|_{G_Y}$ is an inner automorphism of $G_Y$, of the form $\ad_{h}$ for some $h\in G_Y$. We aim to prove that $h\in G_\varepsilon$, which will imply that $\Phi$ acts trivially on $G_\varepsilon$. 

Recall that $Y$ is the minimal subtree of $G_Y$ in $T$ (Lemma~\ref{lem:minsubtree}). The fact that $\varphi$ fixes $\varepsilon$ implies that $\varphi$ acts on $Y$ fixing $w$. 
But $\varphi$ acts on $Y$ by translation by $h$, so $h$ must belong to the $G_Y$-stabiliser of $w$, which is exactly $G_\varepsilon$. 
\end{proof}

\subsection{The action of \texorpdfstring{$\Aut_T(G,\calh^{(t)})$}{AutT(G,Ht} on the vertex groups of \texorpdfstring{$U_{T,\calh}$}{UT,H}}\label{sec:action-on-vertices}

In this section, we assume that $G$ is torsion free and that the $G$-action on $T$ is acylindrical and non-elementary. In particular $Z(G)=\{1\}$. 

\paragraph*{Canonical lifts.} Given a finitely generated subgroup $H\subseteq G$ that is not elliptic in $T$, we define a finite-index subgroup 
$\Out^{H}_T(G,[H]^{(t)})\subseteq\Out_T(G,[H]^{(t)})$ as follows: 
\begin{itemize}
\item If $H$ is not cyclic we keep $\Out^{H}_T(G,[H]^{(t)})=\Out_T(G,[H]^{(t)})$. 
\item If $H$ is cyclic, then $T_H$ is a line. Letting $\Stab^+_{\Aut_T(G)}(H)$ be the subgroup of $\Stab_{\Aut_T(G)}(H)$ of index at most $2$ acting as the identity on $H$, we get a translation homomorphism \[\tau_H:\Stab^+_{\Aut_T(G)}(H)\to \mathbb{Z}.\] We let $\Out^{H}_T(G,[H]^{(t)})$ be the image in $\Out_T(G,[H]^{(t)})$ of $\tau_H^{-1}(\tau_H(H))$ -- here we are using that $G$ is torsion free to view $H$ as a subgroup of $\Aut_T(G)$ consisting of inner automorphisms.
\end{itemize}
Note that the subgroup $\Out^H_T(G,[H]^{(t)})$ only depends on the conjugacy class of $H$ in $G$. Given a finite set $\calh$ of conjugacy classes of finitely generated subgroups of $G$, we let
\[\Out^\calh_T(G,\calh^{(t)}):=\bigcap_{[H]\in\calh_T^\infty}\Out^H_T(G,[H]^{(t)}).\]

\begin{lemma}
\label{lemma:lift}Assume that $G$ is torsion free and that the $G$-action on $T$ is acylindrical and non-elementary, and let $H\subseteq G$ be a finitely generated subgroup that is not elliptic in $T$. Then every $\Phi\in\Out^{H}_T(G,[H]^{(t)})$ has a unique lift in $\Aut(G)$ acting as the identity on $H$ and on $T_H$.
\end{lemma}

\begin{proof}
The uniqueness of such a lift follows from the acylindricity of $T$. We now prove its existence. So let $\Phi\in\Out^H_T(G,[H]^{(t)})$, and let $\varphi\in\Aut(G)$ be a lift acting as the identity on $H$. In particular $\phi$ preserves $T_H$. 

$\drsh$ Assume first that $H$ is cyclic. By definition of $\Out^H_T(G,[H]^{(t)})$, there exists $h\in H$ such that $\tau_H(\varphi)=\tau_H(h)$. Then $\ad_{h^{-1}}\circ\varphi$ is the desired lift of $\Phi$. 

$\drsh$ Assume now that $H$ is not cyclic. Then $\varphi$ preserves the axis of every element $h\in H$ acting hyperbolically on $T$. Since $T$ is acylindrical and $H$ is not cyclic (and thus not virtually cyclic as $G$ is torsion free), the tree $T_H$ is not a line. This forces $\varphi$ to act as the identity on $T_H$, so $\varphi$ is the desired lift of $\Phi$.
\end{proof}

We let $\theta_H:\Out^H_T(G,[H]^{(t)})\to\Aut(G)$ be the homomorphism sending $\Phi$ to its unique lift acting as the identity on $H$ and on $T_H$. 

\paragraph*{Restriction to edge groups.} Let $H\subseteq G$ be a finitely generated subgroup, and let $e\subseteq T_H$ be an edge. We have a homomorphism  \[\res_{H,e}:\Out^{H}_T(G,[H]^{(t)})\to\Aut(G_e)\] defined as follows: for every $\Phi\in\Out^{H}_T(G,[H]^{(t)})$, the lift $\theta_{H}(\Phi)\in\Aut_T(G,[H]^{(t)})$ fixes the edge $e$ so preserves $G_e$, and we let $\res_{H,e}(\Phi)=\theta_{H}(\Phi)|_{G_e}$. We let \[K_{H,e}:=\ker(\res_{H,e}).\]

\paragraph*{Comparison homomorphisms.}
Now let $H_1,H_2\subseteq G$ be two finitely generated subgroups. Let $e\subseteq T$ be an edge, and assume that $e\subseteq T_{H_1}\cap T_{H_2}$. Let $K_{H_1,H_2,e}:=K_{H_1,e}\cap K_{H_2,e}$. Since $Z(G)=\{1\}$,
the group $Z(G_e)$ embeds as a subgroup of $G/Z(G)\simeq\Inn(G)$.

\bigskip

\begin{lemma}
For every $\Phi\in K_{H_1,H_2,e}$, the automorphism $\theta_{H_1}(\Phi)\theta_{H_2}(\Phi)^{-1}$ belongs to $Z(G_e)$, viewed as a subgroup of $\Inn(G)$.
The map
\begin{displaymath}
\begin{array}{cccc}
K_{H_1,H_2,e} & \to & Z(G_e)\\
\Phi & \mapsto & \theta_{H_1}(\Phi)\theta_{H_2}(\Phi)^{-1}
\end{array}
\end{displaymath}
is a group homomorphism.
\end{lemma}

In the sequel, we will denote by $\Comp_{H_1,H_2,e}: K_{H_1,H_2,e}  \to Z(G_e)$ the homomorphism given by the lemma.

\begin{proof}
The automorphisms $\theta_{H_1}(\Phi)$ and $\theta_{H_2}(\Phi)$ differ by the inner automorphism by an element $g_\Phi\in G$ fixing $e$, and by definition of $K_{H_1,H_2,e}$ this inner automorphism acts as the identity on $G_e$, so $g_\Phi\in Z(G_e)$. 

The fact that $\Comp_{H_1,H_2,e}$ is indeed a group homomorphism is a consequence of the following computation: for $\Phi,\Psi \in K_{H_1,H_2,e}$, we have  
\[\begin{array}{rcl}
    \theta_{H_1}(\Phi\Psi) \circ \theta_{H_2}(\Phi\Psi)^{-1}&=&\theta_{H_1}(\Phi)\circ\theta_{H_1}(\Psi)\circ\theta_{H_2}(\Psi)^{-1}\circ\theta_{H_2}(\Phi)^{-1} \\ &=&\theta_{H_1}(\Phi)\circ\ad_{g_\Psi}\circ\theta_{H_2}(\Phi)^{-1}\\ &=&\theta_{H_1}(\Phi)\circ\theta_{H_2}(\Phi)^{-1}\circ\ad_{g_\Psi} \\ &=&\ad_{g_\Phi}\circ\ad_{g_\Psi},
\end{array}\]
where the third equality comes from the fact that $g_\Psi \in G_e$ and $\theta_{H_2}(\Phi)$ is the identity on $G_e$.
\end{proof}

Given a set $X$ of finitely generated subgroups of $G$, and a subtree $\tau\subseteq T$ we finally let 
\[K_{X,\tau}:=\bigcap_{H_1,H_2\in X}\left(\bigcap_{e\subseteq\tau\cap T_{H_1}\cap T_{H_2}}\ker(\Comp_{H_1,H_2,e})\right).\] 

\bigskip

\begin{lemma}
\label{lemma:spanning-K-trivial}Assume that $G$ is torsion free and that the $G$-action on $T$ is acylindrical and non-elementary. Let $Y\in V_1^\infty$. 

Then there exist a finite set $X$ of subgroups of $G_Y$ whose conjugacy classes belong to $\calh_Y^\infty$, and an $X$-spanning finite subtree $\tau\subseteq Y$, such that every $\Phi\in K_{X,\tau}$ has a representative $\varphi\in\Aut(G)$ such that $\varphi|_{G_Y}=\mathrm{id}_{G_Y}$.
\end{lemma}

\begin{proof}
We make the following claim.
\begin{claim}\label{claim:spanning-tree}
There exists a finite subtree $\tau\subseteq Y$ and a finite set $X$ of subgroups of $G_Y$ whose conjugacy classes belong to $\calh_Y^\infty$, satisfying the following properties:
\begin{enumerate}
    \item\label{claim:1} $\tau$ is $X$-spanning, and every edge of $\tau$ belongs to some $T_H$ with $H\in X$;
    \item\label{claim:2} for any $H,H'\in X$, there exists a sequence $H=H^0,\dots,H^n=H'$ of subgroups in $X$ such that for every $i\in\{0,\dots,n-1\}$, the intersection $T_{H^i}\cap T_{H^{i+1}}\cap \tau$  contains an edge;
    \item\label{claim:3} there exist a finite subtree $\tau'\subseteq\tau$ and a subset $X'\subseteq X$ such that 
    \begin{enumerate}
        \item\label{claim:3a} $\tau'$ is $X'$-spanning;
        \item\label{claim:3b} $X'$ contains a representative of every conjugacy class in $\mathcal{H}_Y^\infty$;
        \item\label{claim:3c} if $e,f\subseteq\tau'$ are in the same $G$-orbit, then there exists $k_{e,f}\in G_Y$ such that the following two hold:
        \begin{itemize}
        \item $k_{e,f}e=f$, and \item $k_{e,f}Hk_{e,f}^{-1}\in X$ for every $H\in X'$.
        \end{itemize}
    \end{enumerate}
\end{enumerate}
\end{claim}

\begin{proof}[Proof of the claim]
We first choose an enumeration $H_1,\dots,H_\ell$ of a set of representatives of $\calh_Y^\infty$ (with $|\calh_Y^\infty|=\ell$), with the property that for each $i$ the tree $T_{H_i}$ contains an edge $e_i$ of $T_{H_1}\cup T_{H_2} \cup \cdots \cup T_{H_{i-1}}$. This is done in the following way:
\begin{itemize}
    \item choose $H_1$ arbitrarily;
    \item for every $i\in\{2,\dots,\ell\}$, note that by definition of $Y$, the union $T_{H_1}\cup\dots\cup T_{H_{i-1}}$ must intersect a tree $T_{H'}$ in at least one edge $e_i$, for some subgroup $H'$ whose conjugacy class belongs to $\calh_Y^\infty\setminus\{[H_1],\dots,[H_{i-1}]\}$. We then let $H_i:=H'$.
\end{itemize}

Let $X'=\{H_1, \ldots, H_\ell\}$.  For each $i$, let $\tau_i'$ be a $H_i$-spanning subtree of $T_{H_i}$ that contains $e_i$. Let $\tau'$ be the smallest subtree of $Y$ containing each $\tau_i'$. As $T_{H_1}\cup T_{H_2} \cup \cdots \cup T_{H_\ell}$ is a subtree of $Y$, the tree $\tau'$ is a subtree of $T_{H_1}\cup T_{H_2} \cup \cdots \cup T_{H_\ell}$. Note that $\tau'$ and $X'$ satisfy Properties~(\ref{claim:3a}) and~(\ref{claim:3b}).

\smallskip

For a set $X$ of subgroups whose conjugacy classes belong to $\calh^\infty_Y$ and a subtree $\tau \subset Y$ we define $\Gamma(X,\tau)$ to be the graph with \begin{itemize} \item vertex set $X$, and \item there is an edge between $H,H'\in X$ in $\Gamma(X,\tau)$ if and only if the intersection $T_H \cap T_{H'}\cap \tau$ contains an edge. 
\end{itemize} The set $(X,\tau)$ satisfies Property~(\ref{claim:2}) if and only if $\Gamma(X,\tau)$ is connected. Note that if $\tau_1 \subset \tau_2$ and $X_1 \subset X_2$ then $\Gamma(X_1,\tau_1)$ is a subgraph of $\Gamma(X_2,\tau_2)$. We first prove that the set $(X',\tau')$ satisfies Property~(\ref{claim:2}). For this, we prove by induction on $i \in \{1,\ldots,\ell\}$ that $\Gamma(\{H_1,\ldots,H_i\},\tau')$ is connected.

$\drsh$ The case $i=1$ is immediate.

$\drsh$ Now let $i \geq 2$. Notice that, by construction of $X'$, the edge $e_i \subset T_{H_i}$ is contained in some tree $T_{H_j}$ for $j<i$.  So $T_{H_i} \cap T_{H_j} \cap \tau'$ contains the edge $e_i$ and $H_i$ is connected by an edge to $H_j$ in $\Gamma(\{H_1, \ldots, H_i\},\tau')$. As $\Gamma(\{H_1, \ldots, H_{i-1}\},\tau')$ is connected by our induction hypothesis, so is $\Gamma(\{H_1, \ldots, H_i\},\tau')$.

\smallskip

Let $\mathcal{E}$ be the finite set of all pairs $(e,f)$ of edges in $\tau'$ in the same $G$-orbit (and hence in the same $G_Y$-orbit because the translates of $Y$ form a transverse covering of $T$). For any pair $(e,f)\in\mathcal{E}$, we choose an element $k_{e,f}\in G_Y$ such that $k_{e,f}e=f$ (with $k_{e,e}=1$). We let $X$ be the set of all subgroups of the form $k_{e,f}Hk_{e,f}^{-1}$ with $H\in X'$ and $(e,f)\in\mathcal{E}$. Note that $X'\subseteq X$, and that Property~(\ref{claim:3c}) is automatically satisfied.

We let $\tau$ be the union of all subtrees $k_{e,f}\tau'$ with $(e,f)\in\mathcal{E}$. Notice that our definition of $\mathcal{E}$ ensures that $\tau$ is connected, and as $\mathcal{E}$ is finite $\tau$ is a finite subtree. Since $\tau'$ is $X'$-spanning, $\tau$ is $X$-spanning. Since every edge of $\tau'$ is contained in some $T_H$ with $H\in X'$, every edge of $\tau$ is contained in some $T_H$ with $H\in X$. So Property~(\ref{claim:1}) is satisfied.

\smallskip

We finally prove Property~(\ref{claim:2}). Let $(e,f) \in \mathcal{E}$. Then there exist $i,j$ such that $e \subset T_{H_i}$ and $f \subset T_{H_j}$. As $f$ is an edge of $T_{k_{e,f}H_ik_{e,f}^{-1}}$ the subgroups $H_j$ and $k_{e,f}H_ik_{e,f}^{-1}$ are connected by an edge in $\Gamma(X,\tau)$. We showed above that $\Gamma(X',\tau')$ is connected, so by translation the same is true of $\Gamma(k_{e,f}X'k_{e,f}^{-1},k_{e,f}\tau')$. It follows that $\Gamma(k_{e,f}X'k_{e,f}^{-1},\tau)$ is connected as $k_{e,f}\tau' \subset \tau$. Therefore, any element of $k_{e,f}X'k_{e,f}^{-1}$ is connected by a path in $\Gamma(X,\tau)$ to $k_{e,f}H_ik_{e,f}^{-1}$, which in turn is connected by an edge to $H_j$. As $\Gamma(X',\tau)$ is connected and the choice of $(e,f) \in \cale$ was arbitrary this shows that $\Gamma(X,\tau)$ is connected, which is equivalent to Property~(\ref{claim:2}).
\end{proof}

We now fix $\tau,X,\tau',X'$ as provided by the claim, and get back to the proof of Lemma~\ref{lemma:spanning-K-trivial}. Property~(\ref{claim:2}) from Claim~\ref{claim:spanning-tree}, together with the definition of the homomorphisms $\Comp_{H_1,H_2,e}$, ensures that all homomorphisms $\theta_{H}$ (with $H$ varying in $X$) agree on $K_{X,\tau}$. We denote by $\theta:K_{X,\tau}\to\Aut_T(G,\calh^{(t)})$ the common homomorphism, which is injective. 

Let $\Phi\in K_{X,\tau}$, and let $\varphi:=\theta(\Phi)$. Let $S\subseteq G_Y$ be the subset consisting of all elements in the subgroups in $X'$, together with all elements $h\in G_Y$  such that $h\tau'\cap\tau'$ contains an edge. Since $\tau'$ is $X'$-spanning, Lemma~\ref{lemma:spanning-genset} ensures that the set $S$ generates $G_Y$. So it is enough to prove that $\varphi(s)=s$ for every $s\in S$. 

Let $H\in X$ and $h\in H$. Then $\varphi$ fixes $T_H$ pointwise, so for every $x\in VT_H$, we have $hx=\varphi(hx)=\varphi(h)\varphi(x)=\varphi(h)x$. By acylindricity, and using the fact that $G$ is torsion free, we have $\varphi(h)=h$. So $\varphi$ is the identity on every subgroup in $X$.

Now let $e,f\subseteq\tau'$, and $k\in G$ be such that $ke=f$. Let $k_{e,f}\in G_Y$ be as given by point~(\ref{claim:3c}) from Claim~\ref{claim:spanning-tree}. Since \[K_{X,\tau}\subseteq\bigcap_{H\in X}\left(\bigcap_{ e'\subseteq \tau\cap T_H} \ker(\res_{H,e'})\right),\] we have $\varphi|_{G_{e'}}=\mathrm{id}_{G_{e'}}$ for every edge $e'\subseteq\tau$ (Property~(\ref{claim:1}) from Claim~\ref{claim:spanning-tree} tells us that every edge of $\tau$ is contained in some $T_H$ with $H\in X$). In particular  $\varphi(k_{e,f}^{-1}k)=k_{e,f}^{-1}k$, in other words $\varphi(k)=\varphi(k_{e,f})k_{e,f}^{-1}k$. There remains to prove that $\varphi(k_{e,f})=k_{e,f}$.

Let $H\in X'$ and $h\in H$. We have $k_{e,f}Hk_{e,f}^{-1}\in X$ by property~(\ref{claim:3c}) from Claim~\ref{claim:spanning-tree}. The automorphism $\varphi$ restricts to the identity on $H$ and on $k_{e,f}Hk_{e,f}^{-1}$ (by the definition of $\theta(\Phi)$). We thus have \[\varphi(k_{e,f})h\varphi(k_{e,f})^{-1}=\varphi(k_{e,f}hk_{e,f}^{-1})=k_{e,f}hk_{e,f}^{-1},\] so $\varphi(k_{e,f})=k_{e,f}k'$ for some $k'\in C_G(H)$. As $\varphi$ fixes $e$ and $f$, the element $\varphi(k_{e,f})$ sends $e$ to $f$, showing that $k'\in G_e$. Since $T_H$ is infinite, $G$ is torsion free, and the $G$-action on $T$ is acylindrical, no non-trivial element of $C_G(H)$ can fix a point in $T$. So $k'=1$ and $\varphi(k_{e,f})=k_{e,f}$, as desired. 
\end{proof}

\subsection{Summary}

We summarise the results of the present section in the following theorem.

\begin{thm}
\label{thm:properties-tree-cyl}Let $G$ be a finitely generated torsion-free group, let $T$ be an acylindrical non-elementary $G$-tree and let $\mathcal{H}$ be a finite set of conjugacy classes of finitely generated subgroups of $G$. Let $U_{T,\calh}$ be the skeleton of the $\calh$-transverse covering of $T$, coming with a bipartition of its vertex set as $VU_{T,\calh}=V_0\sqcup V_1$. Then $U_{T,\calh}$ is $\Aut_T(M,\calh)$-invariant and satisfies the following properties: 
\begin{enumerate}
    \item\label{prop1} for every vertex $v\in V_0$, the group $G_v$ is a point stabiliser in $T$;
    \item\label{prop2} for every vertex $v \in V_1\setminus V_1^\infty$, the group $G_v$ is an edge stabiliser in $T$;
    \item\label{prop3} for every vertex $Y\in V_1^\infty$ with $G_Y$ non-cyclic, there exists a finite set $\mathcal{E}_Y$ of pairs $(e,H)$, where $H$ is conjugate to an element of $\calh^\infty_Y$ and $e \in T_H$, such that the image $I_Y$ of $\Out_T^1(G,\calh^{(t)})$ in $\Out(G_Y)$ has a finite-index subgroup $I_Y^0$ which fits into an exact sequence \[1 \to K_Y \to I_Y^0 \xrightarrow{\res} \prod_{(e,H) \in \mathcal{E}_Y} \Aut(G_e),\] where 
    \begin{enumerate}
    \item $K_Y$ is isomorphic to a subgroup of $\prod_{(e,H) \in\mathcal{E}_Y} Z(G_e)$, and 
    \item for every $(e,H)\in\mathcal{E}_Y$, there exists a lift $s_{H}:I_Y^0\to\Aut_Y(G_Y)$ such that for every $\Psi\in I_Y^0$, the automorphism $s_{H}(\Psi)$ fixes the edge $e$, and the projection of $\res(\Psi)$ to $\Aut(G_e)$ coincides with the restriction of $s_{H}(\Psi)$ to $G_e$; 
    \end{enumerate}
    \item\label{prop4} if all edge groups of $T$ are finitely generated, then all vertex and edge groups of $U_{T,\calh}$ are finitely generated;
    \item\label{prop5} the $G$-tree $U_{T,\calh}$ is acylindrical.
\end{enumerate}
\end{thm}

\begin{proof}
The first two items follow from the definition of $U_{T,\calh}$, the fourth is given by Lemma~\ref{lemma:edge-stab-fg} and the fifth follows from Lemma~\ref{lemma:skeleton-acylindrical}. We now focus on the third one.

As a consequence of Lemma~\ref{lem:stability_conjugation}, note that the image of $\Out^1_T(G,\calh^{(t)})$ in $\Out(G_Y)$ is contained in $\Out(G_Y,\calh_Y^{(t)})$, and in fact in $\Out_Y(G_Y,\calh_Y^{(t)})$.

The group $G_Y$ is torsion free, and the action of $G_Y$ on $Y$ is acylindrical, and non-elementary as soon as $G_Y$ is non-cyclic. We can apply the content of Section~\ref{sec:action-on-vertices} to the $G_Y$-action on $Y$. We get a finite-index subgroup $\Out^{\calh_Y}_Y(G_Y,\calh_Y^{(t)})$ of $\Out_Y(G_Y,\calh_Y^{(t)})$, and for every subgroup $H\subseteq G_Y$ whose conjugacy class belongs to $\calh_Y^\infty$, a lift $\theta_H:\Out^{\calh_Y}_Y(G_Y,\calh_Y^{(t)})\to\Aut(G_Y)$, sending every $\Psi\in\Out^{\calh_Y}_Y(G_Y,\calh_Y^{(t)})$ to the unique lift that acts trivially on $H$ and on $T_H$. 

Lemma~\ref{lemma:spanning-K-trivial} ensures that we can find a finite set $X$ of subgroups of $G_Y$ whose conjugacy classes belong to $\calh_Y$, and an $X$-spanning finite subtree $\tau\subseteq Y$, such that $K_{X,\tau}$, defined as a subgroup of $\Out_Y(G_Y,\calh_Y^{(t)})$, is trivial. 

Given $H_1,H_2\in X$, let $\mathcal{E}_{H_1,H_2}$ be the set of all edges in $\tau\cap T_{H_1}\cap T_{H_2}$, and let $\cale_Y$ be the set of all pairs $(e,H_1)$ (and $(e,H_2)$) with $H_1,H_2\in X$ and $e\in\cale_{H_1,H_2}$. 
Let $\res$ be the product of all maps $\res_{H,e}$ with $(e,H)\in\cale_Y$. Since $K_{X,\tau}$ is trivial, the kernel $K_Y$ of $\res$ maps injectively via the product of all comparison homomorphisms, to the product of the centers $Z(G_e)$ with $(e,H)$ varying in $\mathcal{E}_Y$, as desired. 
\end{proof}

\subsection{Application to the case of free-by-cyclic groups}\label{sec:tree-of-cyl-f-by-z}

Let $M=F \rtimes_\Phi \ZZ$ be a polynomially-growing free-by-cyclic group with infinite-order monodromy and let $\calh$ be a finite set of conjugacy classes of finitely generated subgroups of $M$. Let $T_c^*$ be the collapsed tree of cylinders of the limiting tree of $\Phi$, as in Definiton~\ref{de:collapsed-tree-of-cyl}. By Theorem~\ref{theo:collapsed-tree-of-cylinders}, the tree $T_c^*$ is $\Aut(M)$-invariant. In this section, we apply the above construction to get a better understanding of the action of the group $\Out(M,\calh^{(t)})$ on the set of vertices of $U_{T_c^*,\calh}$. This study benefits from the fact that edge groups in $T_c^*$ are virtually abelian. Recall from Notation~\ref{notation:Out_1} that \[\Out_{T_c^*}^1(M,\calh)=\Out_{T_c^*}(M,\calh) \cap \Out_{U_{T_c^*,\calh}}^0(M,\calh).\]

\begin{prop}
\label{prop:tree-cyl-F-by-Z-case}Let $M=F \rtimes_\Phi \ZZ$ be a polynomially-growing free-by-cyclic group with infinite-order monodromy and let $\calh$ be a finite set of conjugacy classes of  finitely generated subgroups of $M$. Let $T_c^*$ be the collapsed tree of cylinders of the limiting tree of $\Phi$. Let $U_{T_c^*,\calh}$ be the associated skeleton, coming with the bipartition $VU_{T_c^*,\calh}=V_0\sqcup V_1$ of its vertex set. 

\begin{enumerate}
    \item For every $v \in V_1$, the restriction homomorphism $\Out_{T_c^*}^1(M,\calh^{(t)}) \to \Out(M_v)$ has virtually polycyclic image.
    \item For every $e \in EU_{T_c^*,\calh}$, the restriction homomorphism $\Out_{T_c^*}^1(M,\calh^{(t)}) \to \Out(M_e)$ has virtually polycyclic image.
\end{enumerate}
\end{prop}

\begin{proof}
We prove the first assertion. Let $v\in V_1$. 

$\drsh$ If $v\in V_1\setminus V_1^\infty$, then $M_v$ is an edge stabiliser of $T_c^*$. The assertion then follows from Corollary~\ref{coro:image-edge-finite}. 

$\drsh$ If $v\in V_1^\infty$, then it corresponds to an infinite subtree $Y$ of the transverse covering $\caly_\calh$ of $T_c^*$, and we have $M_v=M_Y$. We consider the image $I_Y$ of $\Out^1_{T_c^*}(M,\calh^{(t)})$ in $\Out(M_Y)$. Let $I_Y^0$ be the finite-index subgroup of this image given by Theorem~\ref{thm:properties-tree-cyl}(\ref{prop3}). Then there exists a finite set $\cale_Y$ of pairs $(e,H)$, where $e\subseteq T_c^*$ is an edge, and $H$ is a subgroup whose conjugacy class belongs to $\calh_Y$, such that $I_Y^0$ fits into an exact sequence \[1 \to K_Y \to I_v^0 \xrightarrow{\res} \prod_{(e,H) \in \mathcal{E}_Y} \Aut(M_e),\] where
\begin{enumerate}
\item $K_Y$ is isomorphic to a subgroup of $\prod_{(e,H) \in\mathcal{E}_Y} Z(M_e)$, and 
\item for every $(e,H)\in\mathcal{E}_Y$, there exists a lift $s_{H}:I_Y^0\to\Aut_Y(G_Y)$ such that for every $\Psi\in I_Y^0$, the automorphism $s_{H}(\Psi)$ fixes the edge $e$, and the projection of $\res(\Psi)$ to $\Aut(G_e)$ coincides with the restriction of $s_{H}(\Psi)$ to $G_e$;  
\end{enumerate}

By Theorem~\ref{theo:collapsed-tree-of-cylinders}(\ref{item:limiting_tree_as_an_F_tree}), every edge stabiliser in $T_c^*$ is either abelian or isomorphic to the fundamental group of a Klein bottle. Thus, $\prod_{e \in\mathcal{E}_Y} Z(M_e)$ is a finitely generated abelian group. Therefore, in order to prove that $I_Y$ is virtually polycyclic, it suffices to prove that, for every $(e,H)\in \cale_Y$, the image of $I_Y^0$ in $\Aut(M_e)$ is virtually polycyclic.

From the definition of $I_Y$ as a quotient of $\Out_T(G,\calh^{(t)})$, we see that for any $\Psi\in I_Y^0$, the image $s_H(\Psi)$ extends to an automorphism in $\Aut_T(G)$. Therefore, the homomorphism $I_Y^0 \to \Aut(M_e)$ factors through the homomorphism $\Stab_{\Aut_{T}(M,\calh)}(e) \to \Aut(M_e)$ coming from the action on edge groups. By Corollary~\ref{coro:image-edge-finite}, we see that the image of $\Stab_{\Aut_T(M,\calh)}(e)$ in $\Aut(M_e)$ is virtually cyclic. The conclusion follows.

\medskip

We now prove the second assertion. Let $e \in EU_{T,\calh}$ and let $v \in VU_{T,\calh}$ be its endpoint in $V_1$. By Lemma~\ref{lem:factorisation-restriction}, the homomorphism $\Out_T^1(M,\calh^{(t)}) \to \Out(M_e)$ factors through the restriction homomorphism $\Out_T^1(M,\calh^{(t)}) \to \Out(M_v)$. By the first assertion, the image of $\Out_T^1(M,\calh^{(t)})$ in $\Out(M_v)$ is virtually polycyclic. As this property passes to quotients, the conclusion follows.
\end{proof}

\section{Dévissages of free-by-cyclic groups}\label{sec:dévissage}

Let $M$ be a free-by-cyclic group. The goal of this section is to build a canonical tree for $M$ that reduces the study of $\Out(M)$ to the case where $M$ has finite-order monodromy. This is the content of the following theorem, establishing Theorem~\ref{theo:tree} from the introduction.

\begin{thm}
\label{theo:dévissage}Let $M$ be a free-by-cyclic group, and $\calh$ be a finite set of conjugacy classes of finitely generated subgroups of $M$.

Then there exists an $\Aut(M,\calh)$-invariant $(M,\calh)$-tree $S$, equipped with a partition $VS=V_{\per}\sqcup V_{\vpc}$ of its vertex set, 
such that 
\begin{enumerate}
    \item\label{devissage1} no two vertices in $V_\per$ are adjacent;
    \item\label{devissage2} for every $v\in V_\per$, the vertex group $M_v$ is a free-by-cyclic group with finite-order monodromy;
    \item\label{devissage3} for every $v\in V_\vpc$, the image of the restriction homomorphism $\Out_S^0(M,\calh^{(t)})\to\Out(M_v)$ is virtually polycyclic;
    \item\label{devissage4} for every $e \in ET$, the stabiliser $M_e$ is finitely generated, and the image of the restriction homomorphism $\Out_S^0(M,\calh^{(t)})\to\Out(M_e)$ is virtually polycyclic; 
    \item\label{devissage5} the $M$-action on $S$ is acylindrical.
\end{enumerate}
\end{thm}

\begin{defn}[Dévissage]
Let $M$ be a free-by-cyclic group, and let $\calh$ be a finite set of conjugacy classes of finitely generated subgroups of $M$. 

An $\Aut(M,\calh)$-invariant $(M,\calh)$-tree that satisfies all conclusions from Theorem~\ref{theo:dévissage} will be called a \emph{dévissage} of $(M,\calh)$. 
\end{defn}

When $M$ is polynomially growing, Theorem~\ref{theo:dévissage} will be proved by induction on the polynomial degree. The case where $M$ is exponentially growing will then be deduced from the polynomially-growing case, using the canonical JSJ tree coming from relative hyperbolicity. The main tool for carrying this inductive strategy is the following lemma. If $U$ is an $(M,\calh)$-tree and $v \in VU$, recall from Section~\ref{sec:def_auto} the finite set $\Inc^\calh_v$ of $M_v$-conjugacy classes of subgroups of $M_v$ consisting of $\Inc_v$ together with the set of $M_v$-conjugacy classes of subgroups in $\calh$ having $v$ as their unique fixed point. Note that $\Inc^\calh_v$ is a finite set of conjugacy classes of subgroups of $M_v$, and these are finitely generated provided that edge stabilisers of $U$ are.
\begin{lemma}
\label{lemma:dévissage-induction}Let $M$ be a free-by-cyclic group, and let $\calh$ be a finite set of conjugacy classes of finitely generated subgroups of $M$. Assume that there exists an  $\Aut(M,\calh)$-invariant $(M,\calh)$-tree $U$ with finitely generated edge stabilisers, with a bipartition of its vertex set as $VU=W_\slow\sqcup W_\vpc$, such that 
\begin{enumerate}
\item\label{dev1} the $M$-action on $U$ is acylindrical;
    \item\label{dev2} for every vertex $v\in W_\vpc$, the image of the restriction homomorphism \[\Out_U^0(M,\calh^{(t)})\to\Out(M_v)\] is virtually polycyclic;
    \item\label{dev3} for every edge $e\in EU$, the image of the restriction homomorphism \[\Out^0_U(M,\calh^{(t)})\to\Out(M_e)\] is virtually polycyclic;
    \item\label{dev4} for every $v\in W_\slow$, the group $M_v$ is a free-by-cyclic group, 
    and there is a dévissage for $(M_v,\Inc^\calh_v)$.
\end{enumerate}
Then there is a dévissage for $(M,\calh)$.
\end{lemma}

\begin{proof}
If $M$ has finite-order monodromy, then the conclusion is clear because the trivial action of $M$ on a point is a dévissage for $M$ (with this unique point being in $V_\per$). We now assume otherwise. In particular, by Proposition~\ref{prop:periodic_implies_centre}, the group $M$ has trivial centre.

For every $v\in W_\slow$, let $S_v$ be a dévissage for $(M_v,\calh_v)$ -- reduced to a point if $M_v$ has finite-order monodromy. It comes equipped with a partition of its vertex set as $VS_v=V_\per^v\sqcup V_\vpc^v$. Since $M$ is one-ended (see e.g.\ Lemma~\ref{lemma:edge-stab}), every incident edge stabiliser $M_e$ is infinite. The fixed point set of $M_e$ in $S_v$ is non-empty as $[M_e] \in \calh_v$, and bounded as $M_e$ is infinite and the $M_v$-action on $S_v$ is acylindrical. 

Let $S$ be the tree obtained from $U$ by blowing up every vertex $v\in W_\slow$ into $S_v$, attaching any incident edge $e$ to the centre of the fixed point set of $M_e$. Our goal is to prove that $S$ is a dévissage for $(M,\calh)$.

Consider the partition $VS=V_\per\sqcup V_\vpc$, where vertices in $V_\per$ are the $V_\per^v$-vertices in the subtrees $S_v$, and vertices in $V_\vpc$ are the vertices of $U$ in $W_\vpc$ together with the $V_\vpc^v$-vertices in the subtrees $S_v$. 
Let $ES_U$ be the set of edges of $S$ with at least one endpoint in $W_\vpc$. Note that edges in $ES_U$ come from edges of $U$. 
By construction, no two vertices in $V_\per$ are adjacent, and stabilisers of vertices in $V_\per$ are free-by-cyclic groups with finite-order monodromy.

The tree $S$ is an $(M,\calh)$-tree: indeed, if $H\subseteq M$ is a subgroup such that $[H]\in\calh$, then $H$ fixes a vertex $v$ of $U$. If $v\in W_\slow$, then $H$ is elliptic in $S_v$, and therefore in $S$. If $v\in W_\vpc$, then $v$ does not get blown-up when constructing $S$, so $H$ remains elliptic in $S$. The tree $S$ is $\Aut(M,\calh)$-invariant, and the projection $\pi_U \colon S\to U$ (given by collapsing each tree $S_v$ to a point) is $\Aut(M,\calh)$-equivariant by e.g.\ \cite[Lemma~4.8(2)]{Fio} and its proof.

We next claim that the $M$-action on $S$ is acylindrical. Indeed, let $L$ be the acylindricity constant of the $M$-action on $U$ (Assumption~(\ref{dev1})), and for every $v\in W_\slow$, let $L_v$ be the acylindricity constant of the $M_v$-action on $S_v$ (Assumption~(\ref{devissage5}) in the definition of a dévissage). Let \[K:=(L/2+1)(\max_{v\in W_\slow}L_v+1)+L-1.\] Let $\gamma$ be an edge path of length $K+1$ in $S$. The path $\gamma$ decomposes as \[\gamma=\gamma_1 e_1e_2 \gamma_2 e_3e_4\dots \gamma_k e_{2k-1}e_{2k}\gamma_{k+1},\] where each $e_i$ is an edge from $U$ (the vertex incident on $e_{2i-1}$ and $e_{2i}$ belongs to $W_{\vpc}$), and each $\gamma_i$ is a (possibly trivial) path in some $S_v$ with $v\in W_\slow$. If $\pi_U(\gamma)$ has length at least $L+1$, then $\gamma$ has trivial $M$-stabiliser. So we can assume that $\gamma$ crosses at most $L$ edges coming from $U$, i.e.\ $k\le L/2$. Therefore, there are at most $L/2+1$ subpaths of the form $\gamma_i$, and their total length is at least $K+1-L$, so at least one of the subpaths $\gamma_i$ has length $(K+1-L)/(L/2+1)$, which is greater than $\max_{v\in W_\slow}L_v$. If follows that $\gamma_i$ has trivial $M_v$-stabiliser, and therefore it has trivial $M$-stabiliser. In particular $\gamma$ has trivial $M$-stabiliser.  

We finally prove that for every cell $c\in V_\vpc \cup ES$, the image of the restriction homomorphism $\Out^0_S(M,\calh^{(t)})\to\Out(M_c)$ is virtually polycyclic. For vertices in $W_\vpc$ and edges in $ES_U$, this directly follows from the properties of $U$ (Assumptions~(\ref{dev2}) and~(\ref{dev3})) and the $\Aut(M,\calh)$-equivariance of the projection $\pi_U:S\to U$. We now prove it when $c\in V_\vpc^v  \cup ES_v$ for some $v\in W_\slow$. Let $K$ be the kernel of the product of restriction homomorphisms \[\Out^0_S(M,\calh^{(t)})\to\prod_{v\in W_\vpc/M}\Out(M_v) \times \prod_{e \in ES_U/M} \Out(M_e),\] which as we just saw has virtually polycyclic image. By Lemma~\ref{lemma:vpc}, it is enough to prove that the image of the restriction homomorphism $K\to\Out(M_c)$ is virtually polycyclic for every $c\in V_\vpc^v \cup ES_v$ with $v\in W_\slow$. 

To prove this, we first apply Lemma~\ref{lemma:restriction} and derive that the restriction homomorphism $K\to\Out(M_v)$ takes its values in $\Out(M_v,(\Inc^\calh_v)^{(t)})$, where $\Inc^\calh_v$ is the peripheral structure introduced in Section~\ref{sec:def_auto}. Let $K^0\subseteq K$ be the finite-index subgroup that maps to $\Out^0_{S_v}(M_v,(\Inc^\calh_v)^{(t)})$. Since the map $\Out^0_{S_v}(M_v,(\Inc^\calh_v)^{(t)})\to\Out(M_c)$ has virtually polycyclic image by assumption, and since the map $K^0\to\Out(M_c)$ factors through the homomorphism $\Out^0_{S_v}(M_v,(\Inc^\calh_v)^{(t)})\to\Out(M_c)$ (by $\Aut(M,\calh)$-equivariance of the projection $S\to U$), it follows that the image of $K$ in $\Out(M_c)$ is virtually polycyclic. 
\end{proof}

\begin{proof}[Proof of Theorem~\ref{theo:dévissage}]
We first treat the case where $M$ is polynomially growing. In this case the proof is by induction on the polynomial degree $d$ of $M$. 

$\drsh$ The conclusion holds when $d=0$ by taking $S$ to be reduced to a single vertex (in $V_\per$).

$\drsh$ Assume now that $d\ge 1$. Let $T_c^*$ be the collapsed tree of cylinders of the limiting tree of $\Phi$, as in Definition~\ref{de:collapsed-tree-of-cyl}. Let $U_{T_c^*,\calh}$ be the skeleton of the $\calh$-transverse covering of $T_c^*$ (as in Section~\ref{sec:tree-of-cylinders}), and consider the bipartition $VU_{T_c^*,\calh}=V_0\sqcup V_1$ of its vertex set. We aim to apply Lemma~\ref{lemma:dévissage-induction} to the tree $U_{T_c^*,\calh}$, with $W_\slow=V_0$ and $W_\vpc=V_1$, so we now check its assumptions. Theorem~\ref{theo:collapsed-tree-of-cylinders} ensures that $T_c^*$ is an $\Aut(M,\calh)$-invariant $(M,\calh)$-tree, it is acylindrical as an $M$-tree, and its edge stabilisers are finitely generated. Theorem~\ref{thm:properties-tree-cyl} therefore implies that all these properties are also satisfied by $U_{T_c^*,\calh}$. 
Assumptions~(\ref{dev2}),~(\ref{dev3}) are satisfied in view of Proposition~\ref{prop:tree-cyl-F-by-Z-case}. For every $v\in V_0$, the stabiliser $M_v$ is a point stabiliser in $T_c^*$ (Theorem~\ref{thm:properties-tree-cyl}(\ref{prop1})), so it is free-by-cyclic and has growth smaller than $d$ (Theorem~\ref{theo:collapsed-tree-of-cylinders}(4)). Assumption~(\ref{dev4}) is therefore satisfied by induction. Theorem~\ref{theo:dévissage} therefore follows from Lemma~\ref{lemma:dévissage-induction} in this case.

\medskip

We now treat the case where $M$ is exponentially growing. Recall from Theorem~\ref{thm:relhyp} that the group $M$ is hyperbolic relative to the ($\Out(M)$-invariant) set $\calp(M)$ of conjugacy classes of maximal non-cyclic free-by-cyclic subgroups with polynomially-growing monodromy.

Let $\cala$ be the set of conjugacy classes of all elementary subgroups of $M$, i.e.\ cyclic or contained in a subgroup whose conjugacy class belongs to $\calp(M)$. Theorem~\ref{thm:tree_rel_aut} provides us with an $\Aut(M,\calh)$-invariant tree $U$ with finitely generated stabilisers (namely, the canonical JSJ  $(M ,\calp(M) \cup \calh;\cala)$-elementary splitting), coming with a bipartition $VU=V_0\sqcup V_1$ of its vertex set. We now check that $U$ satisfies the assumptions of Lemma~\ref{lemma:dévissage-induction}, with $W_\slow=V_1$ and $W_\vpc=V_0$. Assumption~(\ref{dev1}) follows from Theorem~\ref{thm:tree_rel_aut}(\ref{jsj-1}). Assumption~(\ref{dev2}) follows from Theorem~\ref{thm:tree_rel_aut}(\ref{jsj-3},\ref{jsj-4}), after noticing that there are no QH vertices in our case (see Remark~\ref{rk:no-qh}). Assumption~(\ref{dev3}) follows from Theorem~\ref{thm:tree_rel_aut}(\ref{jsj-5}). Finally, stabilisers of vertices in $W_\slow$ are either cyclic or polynomially-growing free-by-cyclic groups (Theorem~\ref{thm:tree_rel_aut}(\ref{jsj-3})), so Assumption~(\ref{dev4}) is satisfied in all cases (using the first half of this proof for the polynomially-growing case). Theorem~\ref{theo:dévissage} then follows from Lemma~\ref{lemma:dévissage-induction}. 
\end{proof}

\section{Twists}\label{sec:twist}

In this section, we review groups of twists associated to group actions on trees from work of Levitt \cite{Lev}. We discuss some properties of these in the presence of acylindricity, as well as for actions of free-by-cyclic groups.

\subsection{Twists and bitwists: general review}\label{sec:twists_def}

Let $G$ be a finitely generated group, and let $T$ be a splitting of $G$ which is not a line. Let $\calk_T$ be the kernel of the homomorphism \[\Out_T^0(G) \to \prod_{v \in VT/G} \Out(G_v) \] given by the action on vertex stabilisers. By~\cite[Proposition~2.2]{Lev}, the group $\calk_T$ is generated by bitwists around edges of $T$, defined as follows. 
 
Let $e\in ET$ be an edge, and let $v,w\in VT$ be its endpoints. Let $z_v \in G_v$, $z_w \in G_w$ be such that $z_v,z_w \in N_G(G_e)$ and that $z_vz_w^{-1} \in C_G(G_e)$. Let $\overline{T}$ be the splitting of $G$ obtained from $T$ by collapsing every $G$-orbit of edges except the one of $e$. Let $\overline{v}$ and $\overline{w}$ be the images of $v$ and $w$ in $\overline{T}$. The \emph{bitwist} of $T$ around $e$ associated with $z_v$ and $z_w$ is the outer class of the automorphism $D_{z_v,z_w}$ defined as follows. 
\begin{itemize}
\item If $\overline{v}$ and $\overline{w}$ are in distinct $G$-orbits, then $G=G_{\overline{v}} \ast_{G_e} G_{\overline{w}}$ and $D_{z_v,z_w}$ is the automorphism acting on $G_{\overline{v}}$ as the conjugation by $z_v$ and acting on $G_{\overline{w}}$ as the conjugation by $z_w$. 
\item If $\overline{v}$ and $\overline{w}$ are in the same $G$-orbit, then $G=G_{\overline{v}} \ast_{G_e}$. Let $u \in G$ be such that $u \overline{v}=\overline{w}$, so that $u$ is a stable letter for this HNN extension. Define $D_{z_v,z_w}$ as the automorphism acting on $G_{\overline{v}}$ as the global conjugation by $z_v$ and sending $u$ to $z_wuz_v^{-1}$. Two such stable letters differ by the right multiplication by an element of $G_v$, so that $D_{z_v,z_w}$ does not depend on $u$. 
\end{itemize}

When $z_v=1$ (which forces $z_w$ to belong to the centraliser $C_{G_w}(G_e)$), the outer class of $D_{1,z_w}$ is called a \emph{twist}. We denote by $\mathrm{Tw}$ the \emph{group of twists}, i.e.\ the subgroup of $\Out^0_T(G)$ generated by all twists. Note that every element of $\mathrm{Tw}$ induces the trivial outer automorphism on every edge group. We will be able to reduce from bitwists to twists thanks to the following observation:

\bigskip

\begin{lem}
\label{lemma:bitwist}Suppose that $e$ is an edge with endpoints $v,w$ and $N_{G_v}(G_e)=G_e \cdot C_{G_v}(G_e)$ (equivalently, the conjugation action of $N_{G_v}(G_e)$ on $G_e$ is by inner automorphisms of $G_e$). Then every bitwist around $e$ is a product of two twists.
\end{lem}

\bigskip

\begin{proof}
If the bitwist is $D_{z_v,z_w}$, write $z_v=z_ez'_v$ with $z_e\in G_e$ and $z'_v\in C_{G_v}(G_e)$. We can push $z_e$ across $e$ and deduce that $D_{z_v,z_w}$ is in the same outer class as the product $D_{z'_v,1}D_{1,z_e^{-1}z_w}$.
\end{proof}

\subsection{Twists and bitwists for acylindrical actions} We will now prove that when $G$ is torsion free and acts acylindrically on $T$, every bitwist belongs to the group of twists. 

\bigskip

\begin{lemma}
\label{lemma:acylindrical-self-normalising}Let $T$ be a tree equipped with an acylindrical action of a finitely generated torsion-free group $G$. Then every edge $e\in ET$ with non-trivial stabiliser has an incident vertex $u$ such that $N_{G_u}(G_e)=G_e$.
\end{lemma}

\begin{proof}
Since the $G$-action on $T$ is acylindrical and $G_e$ is infinite, the fixed subtree $\Fix_T(G_e)$ of $G_e$ has finite diameter. In particular, the normaliser $N_G(G_e)$ fixes the centre $c$ of $\Fix_T(G_e)$. Let $u$ be an endpoint of $e$ which is farthest from $c$ (this endpoint is unique unless $c$ is the midpoint of $e$). If $g\in G_u$ normalises $G_e$, then $g$ fixes both $u$ and $c$, so $g$ fixes $e$, i.e.\ $g\in G_e$, as desired.
\end{proof}

If $e\in ET$ is an edge such that either $G_e$ is trivial, or $G_e$ is self-normalising in one incident vertex stabiliser, then Lemma~\ref{lemma:bitwist} ensures that every bitwist around $e$ belongs to the group of twists. As a consequence, we reach the following statement.

\begin{cor}
\label{cor:bitwist_is_twist}Let $T$ be a tree equipped with an acylindrical action of a finitely generated torsion-free group $G$. Then the kernel of the restriction homomorphism \[\Out_T^0(G) \to \prod_{v \in VT/G} \Out(G_v)\] is equal to the group of twists. \qed
\end{cor}

\subsection{Twists for actions of free-by-cyclic groups} In this section $M$ is a free-by-cyclic group with fibre $F$.

\begin{lem}
\label{lem:centraliseredgegroups}Let $M$ be a free-by-cyclic group. Let $K\subseteq H$ be two subgroups of $M$, with $H$ finitely generated and $K$ non-trivial. Then $C_H(K)$ is either trivial or a free-by-cyclic group with finite-order monodromy.
\end{lem}

We insist that with our convention that free-by-cyclic groups have finitely generated fibre, an important part of the conclusion is that $C_{H}(K)\cap F$ is finitely generated.

\begin{proof}
If $K \cap F \neq \{1\}$, then $C_{H}(K) \cap F$ is cyclic (possibly trivial), so that $C_{H}(K)$ is cyclic (possibly trivial) or virtually $\ZZ^2$. 

If $K \cap F=\{1\}$, then $K$ is infinite cyclic and a generator $t$ of $K$ induces an automorphism $\psi$ of $H \cap F$. Since $H$ is finitely generated, by~\cite[Proposition~2.3]{FeiHan}, the group $H$ is an HNN extension $F_1\ast_{F_2}$ where $F_1$ and $F_2$ are finitely generated free subgroups of $H \cap F$.

Thus, the element $t$ acts hyperbolically on the tree associated with $F_1\ast_{F_2}$ and the group $C_{H}(K)$ preserves its oriented axis. Thus, we have a translation homomorphism $\tau:C_{H}(K) \to \ZZ$. The kernel of this homomorphism is contained in (a conjugate of) $F_1$ and corresponds to the fixed subgroup of the injective homomorphism $F_1\cap\psi(F_1) \to \psi(F_1)$ induced by $t$. Since $F_1\cap\psi(F_1)$ is finitely generated, it follows from \cite{GoldTur} that this fixed subgroup is finitely generated, in other words $\ker(\tau)$ is finitely generated. It follows that $C_{H}(K)$ is isomorphic to a free-by-cyclic group (i.e.\ with finitely generated fibre). Since $C_{H}(K)$ has infinite centre, Proposition~\ref{prop:periodic_implies_centre} ensures that the monodromy of this free-by-cyclic group has finite order.
\end{proof}

\begin{prop}
\label{prop:twist-of-type-F}Let $M$ be a free-by-cyclic group and let $N$ be a finitely generated subgroup of $M$ such that $Z(N)=\{1\}$. Let $T$ be a splitting of $N$
with finitely generated edge stabilisers, and let $\mathrm{Tw}$ be its group of twists.
\begin{enumerate}
    \item\label{twist1} Any finite-index torsion-free subgroup of $\mathrm{Tw}$ is of type $\mathrm{F}$. If $M=N$ then any torsion-free finite-index subgroup of $\mathrm{Tw} \cap \Out(M,F)$ is also of type $\mathrm{F}$. 
    \item\label{twist2} The group $\mathrm{Tw}$ contains only finitely many conjugacy classes of finite subgroups. If $M=N$, then  $\mathrm{Tw} \cap \Out(M,F)$ also contains only finitely many conjugacy classes of finite subgroups.
\end{enumerate}
\end{prop}

We encourage the reader to assume that $N=M$ as the case $N \neq M$ will only be used in the proof of Theorem~\ref{thm:VF_no_fg_fibre}.

\begin{proof}
We denote by $\vec{E}=\vec{E}(T/N)$ the set of oriented edges of $T/N$, and will write $o(e)$ for the origin of an oriented edge $e$.
Since $Z(N)=\{1\}$, by~\cite[Proposition~3.1]{Lev}, the group $\mathrm{Tw}$ is isomorphic to a quotient of 
\[G=\prod_{e \in \vec{E}}C_{N_{o(e)}}(N_e)\] 
by a subgroup $H$ of $Z(G)$ (Levitt gives a precise description of $H$ but for the purpose of this proof it is enough to know that $H$ is central in $G$). The quotient map $G \to \mathrm{Tw}$ is obtained by taking $z \in C_{N_{o(e)}}(N_e)$ and mapping it to the twist $D_{z,1}$ around the edge $e$.

Let $e \in \vec{E}$. Since edge stabilisers of $T$ are finitely generated, the vertex stabilisers $N_{o(e)}$ are also finitely generated. By Lemma~\ref{lem:centraliseredgegroups} (applied with $H=N_{o(e)}$ and $K=N_e$), the group $C_{N_{o(e)}}(N_e)$ is either trivial, equal to $N_{o(e)}$ when $N_e=\{1\}$ or isomorphic to a free-by-cyclic
group with finite-order monodromy. Let $\vec{E}^\per$ be the subset of $\vec{E}$ consisting of all edges such that $C_{N_{o(e)}}(N_e)$ is trivial or isomorphic to a free-by-cyclic group with finite order monodromy, and let $G^\per=\prod_{e \in \vec{E}^\per} C_{N_{o(e)}}(N_e)$. Notice that, if $N$ is one-ended (in particular, when $N=M$ by e.g.\ Lemma~\ref{lemma:edge-stab}), any edge group of $T$ is non-trivial. In that case, by Lemma~\ref{lem:centraliseredgegroups}, we have $\vec{E}=\vec{E}^\per$ and $G=G^\per$.

One can therefore decompose $G$ as a direct product
\[G=\left(\prod_{e \notin \vec{E}^\per}C_{N_{o(e)}}(N_e)\right) \times G^\per.\]
If $e \notin \vec{E}^\per$, then $N_e=\{1\}$ and $C_{N_{o(e)}}(N_e)=N_{o(e)}$ is not a free-by-cyclic group with finite order monodromy. By Lemma~\ref{lem:centraliseredgegroups} applied to $H=N_{o(e)}$ and $K=Z(N_{o(e)})$, the group $Z(N_{o(e)})$ is necessarily trivial. In particular, as $Z(G) \simeq \prod_{e \in \vec{E}}Z(C_{N_{o(e)}}(N_e))$, the groups $Z(G)$ and $H$ are contained in $G^\per$ and 
\[\mathrm{Tw} \simeq \left(\prod_{e \notin \vec{E}^\per} N_{o(e)}\right) \times G^\per/H.\]

The group $\prod_{e \notin \vec{E}^\per} N_{o(e)}$ is torsion free and of type $\F$ by the work of Feighn--Handel~\cite{FeiHan}. Thus, in order to prove the assertions of Proposition~\ref{prop:twist-of-type-F}, it suffices to prove the results for the group $G^\per/H$ in place of $\mathrm{Tw}$.

\medskip

For every $e \in \vec{E}^\per$, the group $C_{N_{o(e)}}(N_e)$ is isomorphic to a free-by-cyclic group with finite-order monodromy. One can therefore apply Proposition~\ref{prop:periodic_trees} to find a (possibly trivial) splitting $T_{o(e)}$ of each group $C_{N_{o(e)}}(N_e)$ such that

\begin{itemize}
 \item edge and vertex stabilizers of $T_{o(e)}$ are finitely generated free abelian groups;
 \item $Z(C_{N_{o(e)}}(N_e))$ acts on $T_{o(e)}$ trivially;
 \item any fibre $F_{o(e)}$ of $C_{N_{o(e)}}(N_e)$ acts on $T_{o(e)}$ cocompactly.

\end{itemize}

Consider the product action of $G^\per$ on $X=\prod_{e \in \vec{E}^\per} T_{o(e)}$. The space $X$ is contractible as it is a product of trees. We consider $X$ as a CW-complex with cells coming from products of the edges and vertices in the $T_{o(e)}$. For every $e \in \vec{E}^\per$, the group $C_{M_{o(e)}}(M_e)$ acts cocompactly on $T_{o(e)}$, so $G^\per$ acts cocompactly on $X$. As every edge and vertex stabiliser in $T_{o(e)}$ is finitely generated and free abelian and this property is closed under taking finite products, the same is true for the stabiliser of any cell in $X$.

For every $e \in \vec{E}^\per$, the group $Z(C_{N_{o(e)}}(N_e))$ acts trivially on $T_{o(e)}$, so that $Z(G^\per)=\prod_{e \in \vec{E}^\per} Z(C_{N_{o(e)}}(N_e))$ acts trivially on $X$. As $H \subseteq Z(G^\per)$, there is an induced action of $G^\per/H$ on $X$ which is cocompact (it has the same quotient) and with finitely generated abelian cell stabilisers (possibly now with torsion). Since the action of each factor $C_{N_{o(e)}}(N_e)$ on $T_{o(e)}$ is without edge inversions, any element of $G^\per/H$ that preserves a cell fixes it pointwise. As any torsion-free subgroup of a cell stabiliser in $G^\per/H$ is finitely generated free abelian and therefore of type $\F$, Proposition~\ref{prop:actiononcontractiblespacetypeVF} tells us that any torsion-free finite-index subgroup of $G^\per/H$ and hence of $\mathrm{Tw}$ is also of type $\F$. 

When $M=N$, the proof of the first assertion for $\mathrm{Tw} \cap \Out(M,F)$ is similar using the fact that the action of $\mathrm{Tw} \cap \Out(M,F)$ on $X$ is also cocompact. Indeed, since $M$ is one-ended, as explained above, we have $\vec{E}=\vec{E}^\per$ and $G=G^\per$. Thus, for every $e \in \vec{E}$, the fibre $F_{o(e)}:=C_{M_{o(e)}}(M_e) \cap F$ is finitely generated and acts cocompactly on $T_{o(e)}$. Moreover, any element $z \in F_{o(e)}$ maps to a twist $D_{z,1}$ under the quotient $G \to \mathrm{Tw}$ which is fibre and orientation preserving by definition (see Section~\ref{sec:twists_def} for the construction). Thus, $\mathrm{Tw} \cap \Out(M,F)$ acts cocompactly on each factor of $X$, so it also acts cocompactly on $X$. Hence Proposition~\ref{prop:actiononcontractiblespacetypeVF} also applies to $\mathrm{Tw} \cap \Out(M,F)$. This proves the first assertion. 

We now prove the second assertion. First, any finite subgroup of $\mathrm{Tw}$ (or $\mathrm{Tw} \cap \Out(M,F)$ when $N=M$) fixes a cell (in fact, even a vertex) in $X$. Since every cell stabiliser of $X$ is finitely generated abelian, any cell stabiliser has only finitely many finite subgroups. Since the actions of $\mathrm{Tw}$ and $\mathrm{Tw} \cap \Out(M,F)$ on $X$ are cocompact, there are only finitely many conjugacy classes of cell stabilisers, hence finitely many conjugacy classes of finite subgroups.
\end{proof}

\begin{cor}
\label{cor:twist-type-F}Let $M$ be a free-by-cyclic group and let $N$ be a finitely generated subgroup of $M$ such that $Z(N)=\{1\}$. Let $T$ be a splitting of $N$ with finitely generated edge stabilisers, and let $\mathrm{Tw}$ be its group of twists. 

There exists a finite-index subgroup $\Out^\mathrm{Tw}(N)$ of $\Out(N)$ such that 
\begin{itemize}
\item $\Out^\mathrm{Tw}(N) \cap \mathrm{Tw}$ is of type $\F$; 
\item if $M=N$, then $\Out^\mathrm{Tw}(M) \cap \mathrm{Tw} \cap \Out(M,F)$ is also of type $\F$.
\end{itemize}
\end{cor}

\begin{proof}
By Proposition~\ref{prop:twist-of-type-F}(\ref{twist1}), it suffices to prove that there exists a finite-index subgroup $\Out^\mathrm{Tw}(N)$ of $\Out(N)$ such that $\Out^\mathrm{Tw}(N) \cap \mathrm{Tw}$ (and $\Out^\mathrm{Tw}(M) \cap \mathrm{Tw} \cap \Out(M,F)$ when $M=N$) are torsion free. By Proposition~\ref{prop:twist-of-type-F}(\ref{twist2}) it suffices to construct a finite-index subgroup $\Out^\mathrm{Tw}(N)$ of $\Out(N)$ which avoids finitely many conjugacy classes of finite subgroups of $\Out(N)$. But this follows from the fact that $\Out(N)$ is residually finite \cite[Theorem~B]{DKHT25} (see also \cite[Theorem~1.3]{LevittMinasyan} in the specific case where $N$ is relatively hyperbolic).
\end{proof}

\begin{rk}\label{rk:twist-type-F}
If $\calh$ is a set of conjugacy classes of subgroups of $M$ that are elliptic in $T$, then $\mathrm{Tw}\subseteq\Out(M,\calh^{(t)})$. In particular, $\Out^{\mathrm{Tw}}(M,\calh^{(t)}):=\Out^{\mathrm{Tw}}(M)\cap\Out(M,\calh^{(t)})$ is a finite-index subgroup of $\Out(M,\calh^{(t)})$ whose intersection with $\mathrm{Tw}$ is of type $\F$. 
\end{rk}

\section{Extending automorphisms of vertex groups}\label{sec:extend}

Let $M$ be a free-by-cyclic group, let $T$ be an $M$-tree, let $v\in VT$ and let $\Psi_v \in \Out(M_v)$. The following section is dedicated to providing sufficient conditions in order to extend $\Psi_v$ to an outer automorphism $\Psi$ of $M$. Following \cite{Lev}, it is always possible to extend $\Psi_v$ when $\Psi_v \in \Out(M_v,\Inc_v^{(t)})$. When $T$ is acylindrical, one obtains additional information about $\Psi$. This relies on the following lemma.

\bigskip

\begin{lemma}
\label{lemma:acylindrical-out-e}Let $T$ be a tree equipped with an acylindrical action of a finitely generated torsion-free group $G$. Let $\Psi\in\Out^0_T(G)$. Assume that for every $v\in VT$, the image of $\Psi$ under the restriction homomorphism $\Out^0_T(G)\to\Out(G_v)$ belongs to $\Out(G_v,\Inc_v^{(t)})$.

Then for every edge $e\in ET$, the image of $\Psi$ under the restriction homomorphism $\Out_T^0(G)\to\Out(G_e)$ is trivial. 
\end{lemma}

\bigskip

\begin{proof}
Let $e\in ET$. The conclusion of the lemma is clear if $G_e$ is trivial, so we assume otherwise. Lemma~\ref{lemma:acylindrical-self-normalising} provides an endpoint $u$ of $e$ such that $N_{G_u}(G_e)=G_e$. Note that $\Psi$ has a representative $\varphi_u \in \Aut(G)$ such that $\varphi_u(u)=u$ and $\varphi_{u|G_e}=\id_{G_e}$. Let $\varphi_e$ be a representative of $\Psi$ fixing the edge $e$. Then $\varphi_u$ and $\varphi_e$ differ by an inner automorphism $\ad_g$ with $g \in N_{G_u}(G_e)=G_e$. Thus, $\varphi_e$ acts on $G_e$ as a global conjugation by an element of $G_e$ and the image of $\varphi_e$ (and $\Psi$) in $\Out(G_e)$ is trivial. 
\end{proof}

\bigskip

\begin{cor}
\label{cor:extension_trivial_edge}Let $T$ be a tree equipped with an acylindrical action of a finitely generated torsion-free group $G$. Let $v\in VT$ and $\Psi_v\in\Out(G_v,\Inc_v^{(t)})$. Then there exists $\Psi\in\Out^0_T(G)$ such that  
\begin{itemize}
\item the image of $\Psi$ under the restriction homomorphism $\Out^0_T(G)\to\Out(G_v)$ is $\Psi_v$;
\item for every vertex $w\in VT$ not in the orbit of $v$, the image of $\Psi$ under the restriction homomorphism $\Out^0_T(G)\to\Out(G_w)$ is trivial;
\item for every edge $e\in ET$, the image of $\Psi$ under the restriction homomorphism $\Out^0_T(G)\to\Out(G_e)$ is trivial.
\end{itemize}
\end{cor}

\bigskip

\begin{proof}
By~\cite[Proposition~2.1]{Lev}, we can find $\Psi\in\Out^0_T(G)$ satisfying the first two assertions. The third assertion then follows from Lemma~\ref{lemma:acylindrical-out-e}. 
\end{proof}

Suppose now that $M_v$ is a free-by-cyclic subgroup of $M$, with fibre $F_v=M_v \cap F$. Let $\Psi_v \in \Out(M_v,F_v,\Inc_v^{(t)})$. By  Corollary~\ref{cor:extension_trivial_edge}, one may always extend $\Psi_v$ to an outer automorphism $\Psi \in \Out(M)$ preserving $T$. However, the extended outer automorphism $\Psi$ is not necessarily fibre preserving. Our next statements ensure that, when $M_v$ has finite-order monodromy, one can find a finite-index subgroup of $\Out(M_v,F_v,\Inc_v^{(t)})$ all of whose elements extend to fibre-preserving outer automorphisms of $M$. 

\bigskip

\begin{lemma}
\label{lem:magic_lemma}Let $M$ be a free-by-cyclic group, let $T$ be an acylindrical $M$-tree, let $v\in VT$ be such that $M_v$ is a free-by-cyclic subgroup of $M$, with fibre $F_v=M_v \cap F$. Let $\Psi_v \in \Out^+(M_v,F_v,\Inc_v^{(t)})$. 

Suppose that $\Psi_v$ has a representative $\psi_v\in\Aut(M_v)$ such that for every subgroup $H\subseteq M_v$ whose $M_v$-conjugacy class belongs to $\Inc_v$, the restriction $\psi_{v|H}$ coincides with a conjugation by an element of $F_v$. Then there exists $\Psi\in\Out^{0}_T(M,F)\cap\Out^+(M,F)$ such that  
\begin{itemize}
\item the image of $\Psi$ under the restriction homomorphism $\Out^0_T(M)\to\Out(M_v)$ is $\Psi_v$;
\item for every vertex $w\in VT$ not in the orbit of $v$, the image of $\Psi$ under the restriction homomorphism $\Out^0_T(M)\to\Out(M_w)$ is trivial;
\item for every edge $e\in ET$, the image of $\Psi$ under the restriction homomorphism $\Out^0_T(M)\to\Out(M_e)$ is trivial.
\end{itemize}
\end{lemma}

\begin{proof}
Let $\calm$ be the quotient graph of groups for the action of $M$ on $T$, $\tau$ be a spanning tree of $T/M$, and identify $M$ with the fundamental group $\pi_1(\calm,S)$. A generating set for $M$ is then given by the vertex groups $M_w$, together with HNN generators $s_e$ corresponding to unoriented edges outside $\tau$.

We will define an automorphism of $M$ on these generators, check it is indeed well defined, and then observe that it has the claimed properties. To enable this, for each oriented edge $e$ at $v$, let $g_e \in F_v$ be an element such that $\ad_{g_e}$ agrees with $\psi_v$ on the inclusion $\alpha_e(M_e)$. Extend this definition to vertices $w\neq v$ by setting $g_w$ to be the element $g_e$ associated to the first edge in the path in $S$ joining $v$ to $w$, and to edges outside $S$, and with $o(e) \neq v$ by setting $g_e$ to be $g_{o(e)}$. We now define $\psi$ on generators as $\psi_v$ on $M_v$, $\ad_{g_w}$ on $M_w$ with $w \neq v$, and by declaring it to send $s_e$ to $g_{\overline{e}} s_e g_{e}^{-1}$. 

To see that $\psi$ is well defined, we need to check the various edge relations. Let $e$ be an edge of $\tau$. First suppose $v \notin \{ o(e), t(e) \}$. Then $g_{o(e)}$ and $g_{t(e)}$ are the same element of $F_v$, and so the two conjugations agree on $M_e$. If $e$ lies outside $S$, then since $g_e=g_{o(e)}$ and $g_{\overline{e}}=g_{t(e)}$ we see that for $g \in M_e$, \begin{align*} 
 \psi(s_e) \psi(\alpha_e(g)) \psi(s_e)^{-1} &= g_{\overline{e}} s_e g_{e}^{-1} g_{o(e)} \alpha_e(g) g_{o(e)}^{-1} g_e s_e^{-1} g_{\overline{e}}^{-1} \\
    &= g_{\overline{e}} s_e  \alpha_e(g) s_e^{-1} g_{\overline{e}}^{-1}  \\
    &= g_{\overline{e}}  \alpha_{\overline{e}} (g) g_{\overline{e}}^{-1} \\
    &= \psi (\alpha_{\overline{e}}(g)).
    \end{align*}

If $o(e)=v$, then instead we have $\psi(\alpha_e(g))=\psi_v(\alpha_e(g))=g_e\alpha_e(g) g_e^{-1}$, and again the relations hold; a similar argument is needed in the case $\tau(e)=v$.

Note that applying the construction to $\psi_v^{-1}$ (choosing the conjugating elements $g'_e=g_e^{-1}$) gives an inverse to $\psi$, and so this is indeed an automorphism. It is contained in $\Aut^+(M,F)$ since it preserves the exponent sum on the generating set. The outer automorphism $\Psi$ satisfies the first two assertions by construction. The last assertion follows from Lemma~\ref{lemma:acylindrical-out-e}, using the fact that $T$ is acylindrical and that $M$ is finitely generated and torsion free.
\end{proof}

\begin{rk}
For those familiar with the terminology of \cite{BassCoveringTheory}, we note that the construction above can be translated into an automorphism of the underlying graph of groups by setting the edge group and all vertex group morphisms away from $v$ to be trivial, the vertex group morphism at $v$ to be $\psi_v$, and the ``Bass data'' $\delta_e$ to be $1$ if $o(e) \neq v$, and $g_e \in F_v$ if $o(e)=v$.
\end{rk}

If $M$ has finite-order monodromy, recall the construction of the groups $\Out^{\mathrm{IA}}(M,F)$ and $\Out^\per(M,F)$ from Definitions~\ref{de:outia} and~\ref{dfn:out_0_finite_monodromy}.

\begin{cor}
\label{cor:extending_aut_finite_order}Let $M$ be a free-by-cyclic group, let $T$ be an acylindrical $M$-tree, let $v\in VT$ be such that $M_v$ is a free-by-cyclic subgroup of $M$ with fibre $F_v=M_v \cap F$ and finite-order monodromy. For every $\Psi_v \in \Out^\per(M_v,F_v,\Inc_v^{(t)})$, there exists  $\Psi\in\Out^0_T(M,F)$ such that  
\begin{itemize}
\item the image of $\Psi$ under the restriction homomorphism $\Out^0_T(M)\to\Out(M_v)$ is $\Psi_v$;
\item for every vertex $w\in VT$ not in the orbit of $v$, the image of $\Psi$ under the restriction homomorphism $\Out^0_T(M)\to\Out(M_w)$ is trivial;
\item for every edge $e\in ET$, the image of $\Psi$ under the restriction homomorphism $\Out^0_T(M)\to\Out(M_e)$ is trivial.
\end{itemize}
\end{cor}

\begin{proof}
Let $\Psi_v \in \Out^\per(M_v,F_v,\Inc_v^{(t)})$. By Lemma~\ref{lem:magic_lemma}, it suffices to prove that $\Psi_v$ has a representative $\psi_v$ such that, for every $[H] \in \Inc_v$, the automorphism $\psi_v$ acts on $H$ as a global conjugation by an element of $F_v$. 

Suppose first that $M_v$ is virtually abelian. If $M_v$ is abelian, then $\Out(M_v)=\Aut(M_v)$ and any element of $\Out^\per(M_v,F_v,\Inc_v^{(t)})$ fixes the subgroups in $\Inc_v$ pointwise. If $M_v$ is the fundamental group of the Klein bottle, then $\Out^\per(M_v,F_v)=\{1\}$ and the statement is immediate. So we assume from now on that $F_v$ is non-abelian.

By Lemma~\ref{lemma:ia}(\ref{ia-1}), it suffices to prove that, for every $[H] \in \Inc_v$, there exists $\psi_v^H \in \Aut^{\mathrm{IA}}(M,F)$ representing $\Psi_v$ such that $\psi_v^H$ fixes $H$ pointwise. If $H$ is virtually abelian, the result follows from Proposition~\ref{prop:trivial_action_on_abelians_finite_case}(\ref{item:1}). Suppose that $H$ is non-abelian. By Lemma~\ref{lemma:ia}(\ref{ia-3}), $\Psi_v$ has a representative $\psi_v^H$ in $\Aut^{\mathrm{IA}}(M,F)$ fixing $H \cap F$ pointwise. Since $H \cap F$ is non-abelian, the automorphism $\psi_v^H$ fixes $H$ pointwise by Lemma~\ref{lemma:non-abelian-fibre}. 
\end{proof}

\section{Conclusion}\label{sec:conclusion}

We are now in position to prove the main theorem of the paper.

\begin{thm}[restating Theorem~\ref{theo:mccool}]
\label{thm:McCool-free-by-Z}Let $M=F \rtimes \ZZ$ be a free-by-cyclic group and let $\calh$ be a finite set of conjugacy classes of finitely generated subgroups of $M$. Then $\Out(M,\calh^{(t)})$ and $\Out(M,F,\calh^{(t)})$ are of type $\VF$.
\end{thm}

\begin{proof}
Suppose first that $M$ is a free-by-cyclic group with finite-order monodromy.  Then $\Out(M,\calh^{(t)})$ and $\Out(M,F,\calh^{(t)})$ are of type $\VF$ \ by Theorem~\ref{th:section-4}. Hence we may suppose that $M$ has infinite-order monodromy, so that $Z(M)=\{1\}$ by Proposition~\ref{prop:periodic_implies_centre}. 

Let $S$ be a dévissage of $(M,\calh)$, as provided by Theorem~\ref{theo:dévissage}. Let $VS=V_\per \amalg V_\vpc$ be the associated partition of its vertex set. Let $V_\per^*\subseteq V_\per$ be a finite set of representatives of the $M$-orbits of vertices in $V_\per$. For every $v \in V_\per^*$, the group $M_v$ is a free-by-cyclic group with finite-order monodromy, so Theorem~\ref{th:section-4} provides us with a finite-index subgroup $\Out^\per(M_v)\subseteq\Out(M_v)$ such that, for every finite set $\calf$ of $M_v$-conjugacy classes of finitely generated subgroups of $M_v$, the groups $\Out^\per(M_v,\calf^{(t)})$ and $\Out^\per(M_v,F_v,\calf^{(t)})$ are of type $\F$. For every $v\in V_\per^*$, let \[r_v:\Out^0_S(M,\calh)\to\Out(M_v)\] be the restriction homomorphism, and let $\Out^v(M,\calh):=r_v^{-1}(\Out^{\per}(M_v))$. Let $\Out^\mathrm{Tw}(M)$ be the finite-index subgroup of $\Out(M)$ given by Corollary~\ref{cor:twist-type-F} (applied with $N=M$ and $T=S$), which we can assume without loss of generality to be contained in $\Out^0_S(M)$. We finally let  \begin{align*}
    \Out^*(M,\calh):=\Out^\mathrm{Tw}(M) \cap \bigcap_{v \in V_\per^*} \Out^v(M,\calh).
\end{align*}
This is a finite-index subgroup of $\Out(M,\calh)$. We also let $\Out^*(M,\calh^{(t)})$ be its intersection with $\Out(M,\calh^{(t)})$. 

Consider the restriction homomorphisms
\[r_\per:\Out^0_S(M,\calh^{(t)})\to\prod_{v\in V_\per/M}\Out(M_v);\]
\[r_\vpc:\Out^0_S(M,\calh^{(t)})\to\prod_{v\in V_\vpc/M}\Out(M_v);\]
\[r_\edge:\Out^0_S(M,\calh^{(t)})\to\prod_{e\in ES/M}\Out(M_e).\]
Theorem~\ref{theo:dévissage} ensures that the image of $r_\vpc\times r_\edge$ is virtually polycyclic, hence all of its subgroups are of type $\VF$. Let $\hat K$ be the kernel of $r_\vpc\times r_\edge$, and $\hat K_F:\hat K\cap\Out(M,F)$. Let $K:=\hat K\cap\Out^*(M,\calh^{(t)})$ and $K_F:=\hat K_F\cap\Out^*(M,\calh^{(t)})$. Using Proposition~\ref{prop:type_VF_ses}, it is enough to prove that $K$ and $K_F$ are of type $\F$. 

Let $v \in V_\per^*$. Recall from Section~\ref{sec:def_auto} the peripheral structure $\Inc^\calh_v$ of $M_v$ induced by $\calh$ and $\Inc_v$. It is a finite set of $M_v$-conjugacy classes of finitely generated subgroups. Since $K$ is in the kernel of $r_\edge$, Lemma~\ref{lemma:restriction} ensures that the image of $K$ in $\Out(M_v)$ under $r_v$ is contained in $\Out^\per(M_v,(\Inc^\calh_v)^{(t)})$. Similarly, the image of $K_F$ in $\Out(M_v)$ is contained in $\Out^{\per}(M_v,F_v,(\Inc^\calh_v)^{(t)})$.

By Theorem~\ref{theo:dévissage}, the $M$-action on $S$ is acylindrical. We can therefore apply Corollary~\ref{cor:extension_trivial_edge} (respectively Corollary~\ref{cor:extending_aut_finite_order}) and derive that for every $\Psi$ in $\Out^{\per}(M_v,(\Inc^\calh_v)^{(t)})$ (respectively in $\Out^{\per}(M_v,F_v,(\Inc^\calh_v)^{(t)})$), there exists an outer automorphism in $\Out^v(M,\calh^{(t)})$ (respectively in $\Out^v(M,F,\calh^{(t)})$) whose image under $r_v$ is $\Psi$, and whose image in \[\prod_{w \in V_\per^*\setminus\{v\}} \Out(M_w) \times \prod_{w \in V_\vpc/M} \Out(M_w) \times \prod_{e \in ES/M} \Out(M_e)\] is trivial. In other words, $\hat K$ contains an element whose image in $\Out(M_v)$ under $r_v$ is $\Psi$, and whose image in $\Out(M_w)$ with $w\in V_\per^*\setminus\{v\}$ is trivial. 
Therefore \[\prod_{v \in V_\per^*} \Out^{\per}(M_v,(\Inc^\calh_v)^{(t)})\subseteq r_\per(\hat K),\] and as a consequence, $r_\per(K)$ is a finite-index subgroup of $\prod_{v \in V_\per^*} \Out^{\per}(M_v,(\Inc^\calh_v)^{(t)})$. Similarly, $r_\per(K_F)$ is a finite-index subgroup of $\prod_{v \in V_\per^*} \Out^{\per}(M_v,F_v,(\Inc^\calh_v)^{(t)})$. By Theorem~\ref{th:section-4}, the images of $K$ and $K_F$ under $r_\per$ are of type $\F$. Hence, letting $r_\vertex:=r_\per\times r_\vpc$, the images of $K$ and $K_F$ under $r_\vertex$ are of type $\F$. Let $K'$ be the intersection of $K$ with kernel of $r_\vertex$, and let $K_F':=K_F \cap K'$.

In order to prove that $K$ and $K_F$ are of type $\F$, it remains to prove that $K'$ and $K_F'$ are of type $\F$. Since the $M$-action on $S$ is acylindrical, Corollary~\ref{cor:bitwist_is_twist} ensures that $K'$ is contained in the group of twists $\mathrm{Tw}$ of $S$. Note that the image of any element of $\mathrm{Tw}$ under $r_\vertex\times r_\edge$ is trivial, so $K'$ is a finite-index subgroup of $\mathrm{Tw}$ contained in the group $\Out^\mathrm{Tw}(M,\calh)$. By Corollary~\ref{cor:twist-type-F}, the groups $K'$ and $K_F'$ are of type $\F$. This shows that $K$ and $K_F$ are of type $\F$. By Proposition~\ref{prop:type_VF_ses}, it follows that $\Out(M,\calh^{(t)})$ and $\Out(M,F,\calh^{(t)})$ are of type $\VF$.
\end{proof}

We have already explained in Section~\ref{sec:autos} how Theorem~\ref{theo:centralisers-out} follows from Theorem~\ref{theo:mccool}. Hence our proofs of all theorems stated in the introduction are complete.

\section{An extension to infinitely generated fibres}\label{sec:linton}

Using the work of Linton~\cite{Linton}, one can drop the assumption that the fibre $F$ is finitely generated.

\begin{thm}
\label{thm:VF_no_fg_fibre}Let $G$ be a finitely generated group isomorphic to $F \rtimes \ZZ$, where $F$ is a (not necessarily finitely generated) free group. Let $\calh$ be a finite set of conjugacy classes of finitely generated subgroups of $G$. Then $\Out(G,\calh^{(t)})$ is of type $\VF$.
\end{thm}

\begin{proof}
By \cite[Theorem~1.2]{Linton}, there exists a free-by-cyclic group $M$ (i.e.\ with finitely generated fibre) such that $G$ embeds in $M$. We fix once and for all such an embedding $G \subseteq M$.

Suppose first that $Z(G) \neq \{1\}$. By Lemma~\ref{lem:centraliseredgegroups} applied inside $M$ with $H=G$ and $K=Z(G)$, the group $G=C_G(Z(G))$ is free-by-cyclic with finite-order monodromy. The conclusion then follows from Theorem~\ref{th:section-4}.

So we may suppose that $Z(G)=\{1\}$. By~\cite[Theorem~4.3]{Linton}, the group $G$ is hyperbolic relative to a finite set $\calp$ consisting of conjugacy classes of free-by-cyclic subgroups. Moreover, any free-by-cyclic subgroup of $G$ is contained in a subgroup whose conjugacy class belongs to $\calp$. Therefore, $\calp$ is $\Aut(G)$-invariant and $\Out(G,\calh^{(t)}) \subseteq \Out(G,\calp)$.

Recall that $G$ is \emph{one-ended relative to $\calp \cup \calh$} if there does not exist any non-degenerate $(G,\calp \cup \calh)$-free splitting. We divide the proof into two cases, depending on whether $G$ is relatively one-ended or not.

\bigskip

\textbf{Case 1:} $G$ is one-ended relative to $\calp \cup \calh$. 

By Theorem~\ref{thm:tree_rel_aut}, there exists an $\Out(G,\calh^{(t)})$-invariant JSJ $(G,\calp \cup \calh)$-elementary splitting $U$ with non-trivial edge stabilisers. Let $V_0^f$ be a finite set of representatives of the $G$-orbits of flexible vertices of $U$ and let $V_1^p$ be a finite set of representatives of the $G$-orbits of vertices of $U$ whose stabiliser belongs to $\calp$. For every $v \in V_0^f$, the group $G_v$ is QH with trivial fibre. Let $\Sigma_v$ be the associated compact, connected, hyperbolic surface.

By~\cite[Theorem~4.6]{GLsplittings}, there exists a finite-index subgroup $\Out^1(G,\calh^{(t)})$ of $\Out(G,\calh^{(t)})$ which fits into a short exact sequence 
\[1 \to \mathrm{Tw} \to \Out^1(G,\calh^{(t)}) \to \prod_{v \in V_0^f} \mathrm{MCG}^1(\Sigma_v) \times \prod_{v \in V_1^p} \Out(G_v,(\Inc^\calh_v)^{(t)}) \to 1,\] where $\mathrm{MCG}^1(\Sigma_v)$ is a finite-index subgroup of the mapping class group $\mathrm{MCG}(\Sigma_v)$ of $\Sigma_v$, $\Inc^\calh_v=\Inc_v \cup \calh_v^0$ is the finite set of conjugacy classes of subgroups of $G_v$ defined in Section~\ref{sec:def_auto} and $\mathrm{Tw}$ is the group of twists of $U$. By Theorem~\ref{thm:McCool-free-by-Z}, for every $v \in V_1^p$, the group $\Out(G_v,(\Inc^\calh_v)^{(t)})$ is of type $\VF$. For every $v \in V_0^f$, the group $\mathrm{MCG}(\Sigma_v)$ is of type $\VF$ \cite[Section~4]{Ivanov}. Moreover, since $G \subseteq M$ and $Z(G)=\{1\}$, by Corollary~\ref{cor:twist-type-F} (see also Remark~\ref{rk:twist-type-F}), the group $\Out(G,\calh^{(t)})$ has a finite-index subgroup $\Out^{\mathrm{Tw}}(G,\calh^{(t)})$ such that $\Out^{\mathrm{Tw}}(G,\calh^{(t)}) \cap \mathrm{Tw}$ is of type $\F$. By Proposition~\ref{prop:type_VF_ses}, the group $\Out(G,\calh^{(t)})$ is of type $\VF$. This proves Theorem~\ref{thm:VF_no_fg_fibre} when $G$ is relatively one-ended.

\bigskip

\textbf{Case 2:} $G$ is infinitely ended relative to $\calp \cup \calh$. 

Let $\calf$ be the \emph{free factor support} of $\calp \cup \calh$: it is the set of non-trivial conjugacy classes of vertex stabilisers of some (in fact any) maximal $(G,\calp \cup \calh)$-free splitting. This set is independent of the choice of such splitting, so that $\Out(G,\calh^{(t)}) \subseteq \Out(G,\calf)$. A tree is a JSJ free splitting for $G$ relative to $\calp \cup \calh$ (usually called a Grushko splitting of $(G,\calp\cup\calh)$) if and only if the conjugacy class of every non-trivial vertex stabiliser belongs to $\calf$.

The group $\Out(G,\calf)$ acts on the \emph{spine $\calk$ of Outer space of $G$ relative to $\calf$}~\cite{GL-OuterSpace}. This is the geometric realisation of the poset of $(G,\calp \cup \calh)$-free splittings relative to $\calf$, where $T \preccurlyeq T'$ if $T$ is a collapse of $T'$. This space is a contractible simplicial complex \cite[Corollary~4.4]{GL-OuterSpace}. 

The action of $\Out(G,\calf^{(t)})$ on $\calk$ is cocompact as there are only finitely many possibilities for the quotient graph $T/G$ when $T$ is a JSJ $(G,\calp \cup \calh)$-free splitting. Since $\Out(G,\calf^{(t)}) \subseteq \Out(G,\calh^{(t)})$, the group $\Out(G,\calh^{(t)})$ also acts cocompactly on $\calk$. Therefore, by Proposition~\ref{prop:actiononcontractiblespacetypeVF}, in order to prove that $\Out(G,\calh^{(t)})$ is of type $\VF$, it suffices to prove that $\Out(G,\calh^{(t)})$ has a finite-index subgroup $\Out^1(G,\calh^{(t)})$ such that the intersection of $\Out^1(G,\calh^{(t)})$ with every cell stabiliser of $\calk$ is of type $\F$. Since the action of $\Out(G,\calh^{(t)})$ on $\calk$ is cocompact, it suffices to prove that, for every simplex $c$ of $\calk$, there exists a finite-index subgroup $\Out^c(G,\calh^{(t)})$ of $\Out(G,\calh^{(t)})$ such that $\Stab(c) \cap \Out^c(G, \calh^{(t)})$ is of type $\F$.

Let $c$ be a simplex of $\calk$ and let $T_0, \ldots,T_k$ be the vertices of $c$. Then, up to reordering, for every $i \in \{0,\ldots,k-1\}$, the splitting $T_i$ is a collapse of $T_{i+1}$. Since $T_k$ has only finitely many collapses, the group $\Stab(c)$ is a finite-index subgroup of $\Out_{T_k}(G,\calh^{(t)})$. Set $T:=T_k$. Then it suffices to prove that there exists a finite-index subgroup $\Out^c(G,\calh^{(t)})$ of $\Out(G,\calh^{(t)})$ such that $\Out_T^c(G,\calh^{(t)}):=\Out^c(G,\calh^{(t)}) \cap \Out_T(G,\calh^{(t)})$ is of type $\F$. We construct the finite-index subgroup $\Out^c(G,\calh^{(t)})$ in three steps:

\begin{itemize}
    \item  Let $\Out^{c,\star}(G,\calh^{(t)})$ be the finite-index subgroup obtained as the kernel of the homomorphism \[\Out(G,\calh^{(t)}) \to \mathrm{Sym}(\calf) \times \Out(H_1(G,\ZZ/3\ZZ)),\] 
where the map to $\mathrm{Sym}(\calf)$ is given by the permutation action of $\Out(G,\mathcal{H}^{(t)})$ on the finite set of conjugacy classes $\calf$.
    \item Let $\mathrm{Tw}$ be the group of twists of $T$. By Corollary~\ref{cor:twist-type-F} (see also Remark~\ref{rk:twist-type-F}), since $G \subseteq M$ and $Z(G)=\{1\}$, there exists a finite-index subgroup $\Out^{\mathrm{Tw}}(G,\calh^{(t)})$ of $\Out(G,\calh^{(t)})$ such that $\Out^{\mathrm{Tw}}(G,\calh^{(t)}) \cap \mathrm{Tw}$ is of type $\F$. 

    \item Let $[A] \in \calf$. Let $\calp_A \cup \calh_A$ be the finite set of $A$-conjugacy classes of subgroups of $A$ whose $G$-conjugacy class belongs to $\calp \cup \calh$. Then $A$ is one-ended relative to $\calp_A \cup \calh_A$ since $T$ dominates every $(G,\calp \cup \calh)$-free splitting. By the one-ended case, there exists a finite-index subgroup $\Out^A(A,\calh_A^{(t)})$ of $\Out(A,\calh_A^{(t)})$ such that $\Out^A(A,\calh_A^{(t)})$ is of type $\F$. Let $\Out^\calf(G,\calh^{(t)})$ be the preimage of $\prod_{[A] \in \calf} \Out^A(A,\calh_A^{(t)})$ in $\Out(G,\calh^{(t)})$.
 \end{itemize}

Let \[\Out^c(G,\calh^{(t)}):=\Out^{\mathrm{Tw}}(G,\calh^{(t)}) \cap \Out^{c,\star}(G,\calh^{(t)}) \cap \Out^\calf(G,\calh^{(t)}).\] We claim that $\Out_T^c(G,\calh^{(t)})$ is of type $\F$, which will conclude the proof by Proposition~\ref{prop:actiononcontractiblespacetypeVF}. By for instance~\cite[Lemma~3.7(1)]{guerch2026}, since $\Out_T^c(G,\calh^{(t)}) \subseteq \Out_T^{c,\star}(G,\calh^{(t)})$, the group $\Out_T^c(G,\calh^{(t)})$ is contained in $\Out_T^0(G,\calh^{(t)})$, that is, it acts trivially on the quotient graph $T/G$. Thus, $\Out_T^c(G,\calh^{(t)})$ fits in an exact sequence 
\[1 \to \mathrm{Tw} \cap \Out_T^c(G,\calh^{(t)}) \to \Out_T^c(G,\calh^{(t)}) \to \prod_{[A] \in \calf}\Out^A(A,\calh_A^{(t)}). \] By Corollary~\ref{cor:extension_trivial_edge}, since edge groups in $T$ are trivial, the image of $\Out_T^c(G,\calh^{(t)})$ in \[\prod_{[A] \in \calf}\Out^A(A,\calh_A^{(t)})\] is a finite-index subgroup, hence is of type $\F$. Since $\Out_T^c(G,\calh^{(t)}) \subseteq \Out^{\mathrm{Tw}}(G,\calh^{(t)})$, the group $\mathrm{Tw} \cap \Out_T^c(G,\calh^{(t)})$ is of type $\F$. By Proposition~\ref{prop:type_VF_ses}, the group $\Out_T^c(G,\calh^{(t)})$ is of type $\VF$. The conclusion follows.
\end{proof}

\small

\bibliographystyle{alpha}
\bibliography{autofbyz}

@article {AKP,
    AUTHOR = {Algom-Kfir, Y. and Pfaff, C.},
     TITLE = {Normalizers and centralizers of cyclic subgroups generated by
              lone axis fully irreducible outer automorphisms},
   JOURNAL = {New York J. Math.},
    VOLUME = {23},
      YEAR = {2017},
     PAGES = {365--381},
}

@article {AGH,
    AUTHOR = {Andrew, N. and Guerch, Y. and Hughes, S.},
     TITLE = {Automorphisms of relatively hyperbolic groups and the
              {F}arrell-{J}ones conjecture},
   JOURNAL = {Math. Ann.},
  FJOURNAL = {Mathematische Annalen},
    VOLUME = {395},
      YEAR = {2026},
    NUMBER = {4},
     PAGES = {Paper No. 89, 39},
      ISSN = {0025-5831,1432-1807},
}

@article {AHKtorsionhomology,
    AUTHOR = {Andrew, N. and Hughes, S. and Kudlinska, M.},
     TITLE = {Torsion homology growth of polynomially growing free-by-cyclic
              groups},
   JOURNAL = {Rocky Mountain J. Math.},
  FJOURNAL = {The Rocky Mountain Journal of Mathematics},
    VOLUME = {54},
      YEAR = {2024},
    NUMBER = {4},
     PAGES = {933--941},
}

@misc{aimpl,
    title = {Aim{PL}: Rigidity properties of free-by-cyclic groups},
    key = {AimPL},
    note = {Available at \url{http://aimpl.org/freebycyclic}}
}

@article {AMFbyZ,
    AUTHOR = {Andrew, N. and Martino, A.},
     TITLE = {Free-by-cyclic groups, automorphisms and actions on nearly canonical trees},
   JOURNAL = {J. Algebra},
    VOLUME = {604},
      YEAR = {2022},
     PAGES = {451--495},
}

@article {AMCentralisers,
    AUTHOR = {Andrew, N. and Martino, A.},
     TITLE = {Centralisers of linear growth automorphisms of free groups},
   JOURNAL = {Math. Proc. Cambridge Philos. Soc.},
  FJOURNAL = {Mathematical Proceedings of the Cambridge Philosophical
              Society},
    VOLUME = {177},
      YEAR = {2024},
    NUMBER = {2},
     PAGES = {219--240},
}

@misc{bartlett2026,
      title={The conjugacy problem in $\mathrm{Out}({F}_m)$ when the polynomial restrictions are non-growing}, 
      author={G. Bartlett},
      note = {arXiv:2507.13050},
      year={2025},
      eprint={2507.13050},
      archivePrefix={arXiv},
      primaryClass={math.GR},
      url={https://arxiv.org/abs/2507.13050}, 
}

@article {BassCoveringTheory,
    AUTHOR = {Bass, H.},
     TITLE = {Covering theory for graphs of groups},
   JOURNAL = {J. Pure Appl. Algebra},
  FJOURNAL = {Journal of Pure and Applied Algebra},
    VOLUME = {89},
      YEAR = {1993},
    NUMBER = {1-2},
     PAGES = {3--47},
}

@article{Bestvina-Brady,
 author = {Bestvina, M. and Brady, N.},
 title = {Morse theory and finiteness properties of groups},
 journal = {Invent. Math.},
 volume = {129},
 number = {3},
 pages = {445--470},
 year = {1997},
}

@article {BT,
    AUTHOR = {Baumslag, G. and Taylor, T.},
     TITLE = {The centre of groups with one defining relator},
   JOURNAL = {Math. Ann.},
    VOLUME = {175},
      YEAR = {1968},
     PAGES = {315--319},
}

@article {BFcombination,
    AUTHOR = {Bestvina, M. and Feighn, M.},
     TITLE = {A combination theorem for negatively curved groups},
   JOURNAL = {J. Differential Geom.},
    VOLUME = {35},
      YEAR = {1992},
    NUMBER = {1},
     PAGES = {85--101},
}

@article {BFHlaminations,
    AUTHOR = {Bestvina, M. and Feighn, M. and Handel, M.},
     TITLE = {Laminations, trees, and irreducible automorphisms of free
              groups},
   JOURNAL = {Geom. Funct. Anal.},
  FJOURNAL = {Geometric and Functional Analysis},
    VOLUME = {7},
      YEAR = {1997},
    NUMBER = {2},
     PAGES = {215--244},
}

@article{BFHTits1,
 author = {Bestvina, M. and Feighn, M. and Handel, M.},
 title = {The {Tits} alternative for {{\(\text{Out}(F_n)\)}}. {I}: {Dynamics} of exponentially-growing automorphisms},
 journal = {Ann. Math. (2)},
 volume = {151},
 number = {2},
 pages = {517--623},
 year = {2000},
}

@article {BFHTits2,
    AUTHOR = {Bestvina, M. and Feighn, M. and Handel, M.},
     TITLE = {The {T}its alternative for {${\rm Out}(F_n)$}. {II}. {A}
              {K}olchin type theorem},
   JOURNAL = {Ann. of Math. (2)},
  FJOURNAL = {Annals of Mathematics. Second Series},
    VOLUME = {161},
      YEAR = {2005},
    NUMBER = {1},
     PAGES = {1--59},
}

@article{BFHTits3,
 author = {Bestvina, M. and Feighn, M. and Handel, M.},
 title = {Solvable subgroups of {{\(\text{Out}(F_n)\)}} are virtually {Abelian}.},
 journal = {Geom. Dedicata},
 volume = {104},
 pages = {71--96},
 year = {2004}
}

@article {BFH2023,
    AUTHOR = {Bestvina, M. and Feighn, M. and Handel, M.},
     TITLE = {A {M}c{C}ool {W}hitehead type theorem for finitely generated subgroups of {$\Out(F_n)$}},
   JOURNAL = {Ann. H. Lebesgue},
    VOLUME = {6},
      YEAR = {2023},
     PAGES = {65--94},
}

@article {BH,
    AUTHOR = {Bestvina, M. and Handel, M.},
     TITLE = {Train tracks and automorphisms of free groups},
   JOURNAL = {Ann. of Math. (2)},
    VOLUME = {135},
      YEAR = {1992},
    NUMBER = {1},
     PAGES = {1--51},
}

@article {BMV,
    AUTHOR = {Bogopolski, O. and Martino, A. and Ventura, E.},
     TITLE = {The automorphism group of a free-by-cyclic group in rank 2},
   JOURNAL = {Comm. Algebra},
    VOLUME = {35},
      YEAR = {2007},
    NUMBER = {5},
     PAGES = {1675--1690},
}

@article {Bri,
    AUTHOR = {Brinkmann, P.},
     TITLE = {Hyperbolic automorphisms of free groups},
   JOURNAL = {Geom. Funct. Anal.},
    VOLUME = {10},
      YEAR = {2000},
    NUMBER = {5},
     PAGES = {1071--1089},
}

@article {CharneyVogtmann,
    AUTHOR = {Charney, R. and Vogtmann, K.},
     TITLE = {Finiteness properties of automorphism groups of right-angled
              {A}rtin groups},
   JOURNAL = {Bull. Lond. Math. Soc.},
  FJOURNAL = {Bulletin of the London Mathematical Society},
    VOLUME = {41},
      YEAR = {2009},
    NUMBER = {1},
     PAGES = {94--102},
      ISSN = {0024-6093,1469-2120},
   MRCLASS = {20F36},
  MRNUMBER = {2481994},
MRREVIEWER = {Volker\ Gebhardt},
       DOI = {10.1112/blms/bdn108},
       URL = {https://doi.org/10.1112/blms/bdn108},
}

@article {CollinsLevin1983,
    AUTHOR = {Collins, D. J. and Levin, F.},
     TITLE = {Automorphisms and {H}opficity of certain {B}aumslag-{S}olitar
              groups},
   JOURNAL = {Arch. Math. (Basel)},
    VOLUME = {40},
      YEAR = {1983},
    NUMBER = {5},
     PAGES = {385--400},
}

@misc{CHHL,
 author = {Coulon, R. and Hilion, A. and Horbez, C. and Levitt, G.},
 title = {Pol{E}xp growth for automorphisms of toral relatively hyperbolic groups},
 year = {2025},
 note = {arXiv:2512.05569},
}

@misc{Cul,
 author = {Culler, M.},
 title = {Finite groups of outer automorphisms of a free group},
 year = {1984},
 language = {English},
 howpublished = {Contributions to group theory, {Contemp}. {Math}. 33, 197-207 (1984).}
}

@article {CM,
    AUTHOR = {Culler, M. and Morgan, J.W.},
     TITLE = {Group actions on $\mathbb{R}$-trees},
   JOURNAL = {Proc. London Math. Soc. (3)},
    VOLUME = {55},
      YEAR = {1987},
    NUMBER = {3},
     PAGES = {571--604},
}

@article {CL95,
    AUTHOR = {Cohen, M.M. and Lustig, M.},
     TITLE = {Very small group actions on {$\mathbb{R}$}-trees and {D}ehn twist
              automorphisms},
   JOURNAL = {Topology},
  FJOURNAL = {Topology. An International Journal of Mathematics},
    VOLUME = {34},
      YEAR = {1995},
    NUMBER = {3},
     PAGES = {575--617},
}

@article {CL99,
    AUTHOR = {Cohen, M.M. and Lustig, M.},
     TITLE = {The conjugacy problem for {D}ehn twist automorphisms of free
              groups},
   JOURNAL = {Comment. Math. Helv.},
  FJOURNAL = {Commentarii Mathematici Helvetici},
    VOLUME = {74},
      YEAR = {1999},
    NUMBER = {2},
     PAGES = {179--200},
}

@article {Cooper,
    AUTHOR = {Cooper, D.},
     TITLE = {Automorphisms of free groups have finitely generated fixed
              point sets},
   JOURNAL = {J. Algebra},
  FJOURNAL = {Journal of Algebra},
    VOLUME = {111},
      YEAR = {1987},
    NUMBER = {2},
     PAGES = {453--456},
}

@article {CV86,
    AUTHOR = {Culler, M. and Vogtmann, K.},
     TITLE = {Moduli of graphs and automorphisms of free groups},
   JOURNAL = {Invent. Math.},
  FJOURNAL = {Inventiones Mathematicae},
    VOLUME = {84},
      YEAR = {1986},
    NUMBER = {1},
     PAGES = {91--119},
}

@article{dahmani2016suspensions,
  AUTHOR = {Dahmani, F.},
  TITLE = {On suspensions and conjugacy of hyperbolic automorphisms},
  JOURNAL = {Trans. Amer. Math. Soc.},
  VOLUME = {368},
  YEAR = {2016},
  NUMBER = {8},
  PAGES = {5565--5577},
}

@article {DFMT,
    AUTHOR = {Dahmani, F. and Francaviglia, S. and Martino,
              A. and Touikan, N.},
     TITLE = {The conjugacy problem for {${\rm Out}(F_3)$}},
   JOURNAL = {Forum Math. Sigma},
    VOLUME = {13},
      YEAR = {2025},
     PAGES = {Paper No. e41, 21},
}

@misc {DKHT25,
    AUTHOR = {Dahmani, F. and Hughes, S. and Kudlinska, M. and Touikan, N.},
     TITLE = {Free-by-cyclic groups are conjugacy separable},
      NOTE = {arXiv:2509.22346},
      YEAR = {2025},
}

@article {DahLi,
    AUTHOR = {Dahmani, F. and Li, R.},
     TITLE = {Relative hyperbolicity for automorphisms of free products and
              free groups},
   JOURNAL = {J. Topol. Anal.},
    VOLUME = {14},
      YEAR = {2022},
    NUMBER = {1},
     PAGES = {55--92},
}

@misc{DT24,
AUTHOR = {Dahmani, F. and Touikan, N.},
TITLE = {Unipotent linear suspensions of free groups},
NOTE = {arXiv:2305.11274v2},
YEAR = {2024},
}

@article{FH,
 author = {Feighn, M. and Handel, M.},
 title = {Abelian subgroups of {{\(\text{Out}(F_n)\)}}.},
 journal = {Geom. Topol.},
 issn = {1465-3060},
 volume = {13},
 number = {3},
 pages = {1657--1727},
 year = {2009},
}

@article {FeiHan,
    AUTHOR = {Feighn, M. and Handel, M.},
     TITLE = {Mapping tori of free group automorphisms are coherent},
   JOURNAL = {Ann. of Math. (2)},
    VOLUME = {149},
      YEAR = {1999},
    NUMBER = {3},
     PAGES = {1061--1077},
}

@article {FeiHan25,
    AUTHOR = {Feighn, M. and Handel, M.},
     TITLE = {The conjugacy problem for {UPG} elements of {${\rm
              Out}(F_n)$}},
   JOURNAL = {Geom. Topol.},
  FJOURNAL = {Geometry \& Topology},
    VOLUME = {29},
      YEAR = {2025},
    NUMBER = {4},
     PAGES = {1693--1817},
}

@misc{Fio,
author = {Fioravanti, E.},
title = {Growth of automorphisms of virtually special groups},
year = {2025},
note = {arXiv:2501.12321},
}

@article{FMScocompact,
     author = {Francaviglia, S. and Martino, A. and Syrigos, D.},
     title = {On the action of relatively irreducible automorphisms on their train tracks},
     journal = {Annales de l'Institut Fourier},
     year = {2026},
}

@article {FMSmindisplaced,
    AUTHOR = {Francaviglia, S. and Martino, A. and Syrigos, D.},
     TITLE = {The minimally displaced set of an irreducible automorphism of
              {$F_N$} is co-compact},
   JOURNAL = {Arch. Math. (Basel)},
  FJOURNAL = {Archiv der Mathematik},
    VOLUME = {116},
      YEAR = {2021},
    NUMBER = {4},
     PAGES = {369--383},
}

@misc{gaulus,
  title={The mapping-torus of a free group automorphism is hyperbolic relative to the canonical subgroups of polynomial growth},
  author={Gautero, F. and Lustig, M.},
  note={arXiv:0707.0822},
  year={2007}
}

@book {Geo,
    AUTHOR = {Geoghegan, R.},
     TITLE = {Topological methods in group theory},
    SERIES = {Graduate Texts in Mathematics},
    VOLUME = {243},
 PUBLISHER = {Springer, New York},
      YEAR = {2008},
     PAGES = {xiv+473},
}

@article {Ger,
    AUTHOR = {Gersten, S. M.},
     TITLE = {Fixed points of automorphisms of free groups},
   JOURNAL = {Adv. in Math.},
    VOLUME = {64},
      YEAR = {1987},
    NUMBER = {1},
     PAGES = {51--85},
}

@article {Ghosh,
    AUTHOR = {Ghosh, P.},
     TITLE = {Relative hyperbolicity of free-by-cyclic extensions},
   JOURNAL = {Compos. Math.},
    VOLUME = {159},
      YEAR = {2023},
    NUMBER = {1},
     PAGES = {153--183},
}

@misc{GHSS,
AUTHOR = {Guerch, Y. and Hughes, S. and S{\'a}nchez Salda{\~n}a, L.J.},
TITLE = {Centralisers and the virtually cyclic dimension of {$\Out(F_N)$}},
NOTE = {arXiv:2308.01590v1},
YEAR = {2023},
}

@article {GoldTur,
    AUTHOR = {Goldstein, R.Z. and Turner, E.C.},
     TITLE = {Fixed subgroups of homomorphisms of free groups},
   JOURNAL = {Bull. London Math. Soc.},
    VOLUME = {18},
      YEAR = {1986},
    NUMBER = {5},
     PAGES = {468--470},
}

@article {GW,
    AUTHOR = {Gonçalves, D. and Wong, P.},
     TITLE = {Twisted conjugacy classes in nilpotent groups},
   JOURNAL = {J. Reine Angew. Math.},
    VOLUME = {633},
      YEAR = {2009},
     PAGES = {11--27},
}

@article {GuerchRoots,
    AUTHOR = {Guerch, Y.},
     TITLE = {Roots of outer automorphisms of free groups and centralizers
              of abelian subgroups of {${\rm Out}(F_N)$}},
   JOURNAL = {Michigan Math. J.},
  FJOURNAL = {Michigan Mathematical Journal},
    VOLUME = {75},
      YEAR = {2025},
    NUMBER = {5},
     PAGES = {1003--1020},
     }

@misc{guerch2026,
    AUTHOR={Guerch, Y.},
    TITLE={Aperiodicity properties of automorphism groups of free products}, 
    NOTE = {arXiv:2601.03947},
    YEAR = {2026},
}

@article {Guirardellimit,
    AUTHOR = {Guirardel, V.},
     TITLE = {Limit groups and groups acting freely on {$\mathbb{R}^n$}-trees},
   JOURNAL = {Geom. Topol.},
  FJOURNAL = {Geometry and Topology},
    VOLUME = {8},
      YEAR = {2004},
     PAGES = {1427--1470},
}

@article {GL-defo,
    AUTHOR = {Guirardel, V. and Levitt, G.},
     TITLE = {Deformation spaces of trees},
   JOURNAL = {Groups Geom. Dyn.},
    VOLUME = {1},
      YEAR = {2007},
    NUMBER = {2},
     PAGES = {135--181},
}

@article {GLjsj,
    AUTHOR = {Guirardel, V. and Levitt, G.},
     TITLE = {J{SJ} decompositions of groups},
   JOURNAL = {Ast\'erisque},
    volume = {395},
      YEAR = {2017},
     PAGES = {vii+165},
}

@article {GLMcCool,
    AUTHOR = {Guirardel, V. and Levitt, G.},
     TITLE = {Mc{C}ool groups of toral relatively hyperbolic groups},
   JOURNAL = {Algebr. Geom. Topol.},
    VOLUME = {15},
      YEAR = {2016},
    NUMBER = {6},
     PAGES = {3485--3534},
}

@article {GLsplittings,
    AUTHOR = {Guirardel, V. and Levitt, G.},
     TITLE = {Splittings and automorphisms of relatively hyperbolic groups},
   JOURNAL = {Groups Geom. Dyn.},
    VOLUME = {9},
      YEAR = {2015},
    NUMBER = {2},
     PAGES = {599--663},
}

@article {GL-cylinders,
    AUTHOR = {Guirardel, V. and Levitt, G.},
     TITLE = {Trees of cylinders and canonical splittings},
   JOURNAL = {Geom. Topol.},
    VOLUME = {15},
      YEAR = {2011},
    NUMBER = {2},
     PAGES = {977--1012},
}

@article {GL-OuterSpace,
    AUTHOR = {Guirardel, V. and Levitt, G.},
     TITLE = {The outer space of a free product},
   JOURNAL = {Proc. Lond. Math. Soc. (3)},
  FJOURNAL = {Proceedings of the London Mathematical Society. Third Series},
    VOLUME = {94},
      YEAR = {2007},
    NUMBER = {3},
     PAGES = {695--714},
}

@article {Hagen,
    AUTHOR = {Hagen, M.},
     TITLE = {A remark on thickness of free-by-cyclic groups},
   JOURNAL = {Illinois J. Math.},
    VOLUME = {63},
      YEAR = {2019},
    NUMBER = {4},
     PAGES = {633--643},
}

@article {HaglundWise,
    AUTHOR = {Haglund, F. and Wise, D.T.},
     TITLE = {A note on finiteness properties of graphs of groups},
   JOURNAL = {Proc. Amer. Math. Soc. Ser. B},
  FJOURNAL = {Proceedings of the American Mathematical Society. Series B},
    VOLUME = {8},
      YEAR = {2021},
     PAGES = {121--128},
      ISSN = {2330-1511},
   MRCLASS = {20J06 (20E08)},
  MRNUMBER = {4234060},
MRREVIEWER = {Olympia\ Talelli},
       DOI = {10.1090/bproc/81},
       URL = {https://doi.org/10.1090/bproc/81},
}

@article {HM,
    AUTHOR = {Handel, M. and Mosher, L.},
     TITLE = {Subgroup decomposition in {${\rm Out}(F_n)$}},
   JOURNAL = {Mem. Amer. Math. Soc.},
    VOLUME = {264},
      YEAR = {2020},
    NUMBER = {1280},
     PAGES = {vii+276},
}

@article {Ivanov,
    AUTHOR = {Ivanov, N. V.},
     TITLE = {Complexes of curves and {T}eichm\"uller modular groups},
   JOURNAL = {Uspekhi Mat. Nauk},
  FJOURNAL = {Akademiya Nauk SSSR i Moskovskoe Matematicheskoe Obshchestvo.
              Uspekhi Matematicheskikh Nauk},
    VOLUME = {42},
      YEAR = {1987},
    NUMBER = {3(255)},
     PAGES = {49--91, 255},
}

@article{Khr,
 author = {Khramtsov, D. G.},
 title = {Finite groups of automorphisms of free groups},
 journal = {Math. Notes},
 volume = {38},
 pages = {721--724},
 year = {1985}
}

@article {Khr95,
    AUTHOR = {Khramtsov, D. G.},
     TITLE = {Decidability of the conjugacy problem for finite subgroups of
              automorphism groups of free groups},
   JOURNAL = {Algebra i Logika},
  FJOURNAL = {Sibirski\u{\i} Fond Algebry i Logiki. Algebra i Logika},
    VOLUME = {34},
      YEAR = {1995},
    NUMBER = {5},
     PAGES = {558--606, 609},
}

@article {KLV01,
    AUTHOR = {Krsti\'c, S. and Lustig, M. and Vogtmann, K.},
     TITLE = {An equivariant {W}hitehead algorithm and conjugacy for roots
              of {D}ehn twist automorphisms},
   JOURNAL = {Proc. Edinb. Math. Soc. (2)},
  FJOURNAL = {Proceedings of the Edinburgh Mathematical Society. Series II},
    VOLUME = {44},
      YEAR = {2001},
    NUMBER = {1},
     PAGES = {117--141},
}

@article {KV93,
    AUTHOR = {Krsti\'{c}, S. and Vogtmann, K.},
     TITLE = {Equivariant outer space and automorphisms of free-by-finite
              groups},
   JOURNAL = {Comment. Math. Helv.},
  FJOURNAL = {Commentarii Mathematici Helvetici},
    VOLUME = {68},
      YEAR = {1993},
    NUMBER = {2},
     PAGES = {216--262},
}

@article {KV24,
    AUTHOR = {Kudlinska, M. and Valiunas, M.},
     TITLE = {Free-by-cyclic groups are equationally {N}oetherian},
   JOURNAL = {Math. Ann.},
  FJOURNAL = {Mathematische Annalen},
    VOLUME = {392},
      YEAR = {2025},
    NUMBER = {1},
     PAGES = {525--546},
}

@article {LN,
    AUTHOR = {Leary, I. J. and Nucinkis, B. E. A.},
     TITLE = {Some groups of type {$VF$}},
   JOURNAL = {Invent. Math.},
  FJOURNAL = {Inventiones Mathematicae},
    VOLUME = {151},
      YEAR = {2003},
    NUMBER = {1},
     PAGES = {135--165},
}

@article {Lev,
    AUTHOR = {Levitt, G.},
     TITLE = {Automorphisms of hyperbolic groups and graphs of groups},
   JOURNAL = {Geom. Dedicata},
    VOLUME = {114},
      YEAR = {2005},
     PAGES = {49--70},
}

@article {LevittGBS,
    AUTHOR = {Levitt, G.},
     TITLE = {On the automorphism group of generalized {B}aumslag--{S}olitar groups},
   JOURNAL = {Geom. Topol.},
    VOLUME = {11},
      YEAR = {2007},
    NUMBER = {1},
     PAGES = {473--515},
}

@article {Levittgrowth,
    AUTHOR = {Levitt, G.},
     TITLE = {Counting growth types of automorphisms of free groups},
   JOURNAL = {Geom. Funct. Anal.},
    VOLUME = {19},
      YEAR = {2009},
    NUMBER = {4},
     PAGES = {1119--1146},
}

@article {LevittRankOfGBS,
    AUTHOR = {Levitt, G.},
     TITLE = {Generalized {B}aumslag--{S}olitar groups: rank and finite index subgroups},
   JOURNAL = {Ann. Inst. Fourier (Grenoble)},
    VOLUME = {65},
      YEAR = {2015},
    NUMBER = {2},
     PAGES = {725--762},
}

@article {LevittMinasyan,
    AUTHOR = {Levitt, G. and Minasyan, A.},
     TITLE = {Residual properties of automorphism groups of (relatively)
              hyperbolic groups},
   JOURNAL = {Geom. Topol.},
  FJOURNAL = {Geometry \& Topology},
    VOLUME = {18},
      YEAR = {2014},
    NUMBER = {5},
     PAGES = {2985--3023},
}

@article {Linton,
    AUTHOR = {Linton, M.},
     TITLE = {Embedding finitely generated free-by-cyclic groups in
              \{finitely generated free\}-by-cyclic groups},
   JOURNAL = {Int. Math. Res. Not. IMRN},
  FJOURNAL = {International Mathematics Research Notices. IMRN},
      YEAR = {2026},
      volume={2016},
    NUMBER = {4},
     PAGES = {Paper No. rnag020, 15},
}

@article {Los,
    AUTHOR = {Los, J.E.},
     TITLE = {On the conjugacy problem for automorphisms of free groups},
   JOURNAL = {Topology},
    VOLUME = {35},
      YEAR = {1996},
    NUMBER = {3},
     PAGES = {779--808},
}

@incollection {Lustig,
    AUTHOR = {Lustig, M.},
     TITLE = {Conjugacy and centralizers for iwip automorphisms of free
              groups},
 BOOKTITLE = {Geometric group theory},
    SERIES = {Trends Math.},
     PAGES = {197--224},
 PUBLISHER = {Birkh\"auser, Basel},
      YEAR = {2007},
}

@article {Macura,
    AUTHOR = {Macura, N.},
     TITLE = {Detour functions and quasi-isometries},
   JOURNAL = {Q. J. Math.},
    VOLUME = {53},
      YEAR = {2002},
    NUMBER = {2},
     PAGES = {207--239},
}

@article {MutanguhaPretreesUniv,
    AUTHOR = {Mutanguha, J.P.},
     TITLE = {Limit trees for free group automorphisms: universality},
   JOURNAL = {Forum Math. Sigma},
  FJOURNAL = {Forum of Mathematics. Sigma},
    VOLUME = {12},
      YEAR = {2024},
     PAGES = {Paper No. e122, 45},
}

@misc {MutanguhaPG,
    AUTHOR = {Mutanguha, J.P.},
     TITLE = {On polynomial free-by-cyclic groups},
      NOTE = {arXiv:2412.16150},
      YEAR = {2024},
}

@article{Nie,
 author = {Nielsen, J.},
 title = {Die {Isomorphismengruppe} der freien {Gruppen}.},
 fjournal = {Mathematische Annalen},
 journal = {Math. Ann.},
 issn = {0025-5831},
 volume = {91},
 pages = {169--209},
 year = {1924},
}

@article {RS,
    AUTHOR = {Rips, E. and Sela, Z.},
     TITLE = {Cyclic splittings of finitely presented groups and the
              canonical {JSJ} decomposition},
   JOURNAL = {Ann. of Math. (2)},
    VOLUME = {146},
      YEAR = {1997},
    NUMBER = {1},
     PAGES = {53--109},
}

@PhDThesis{Rod,
author = {Rodenhausen, M.},
title = {Centralisers of polynomially growing automorphisms of free groups},
school = {Universität Bonn}, 
year ={2013},
}

@article {RW,
    AUTHOR = {Rodenhausen, M. and Wade, R.D.},
     TITLE = {Centralisers of {D}ehn twist automorphisms of free groups},
   JOURNAL = {Math. Proc. Cambridge Philos. Soc.},
    VOLUME = {159},
      YEAR = {2015},
    NUMBER = {1},
     PAGES = {89--114},
}

@article {Sel95,
    AUTHOR = {Sela, Z.},
     TITLE = {The isomorphism problem for hyperbolic groups. {I}},
   JOURNAL = {Ann. of Math. (2)},
  FJOURNAL = {Annals of Mathematics. Second Series},
    VOLUME = {141},
      YEAR = {1995},
    NUMBER = {2},
     PAGES = {217--283},
}

@article{Zim,
 author = {Zimmermann, B.},
 title = {{\"U}ber {Homoeomorphismen} {{\(n\)}}-dimensionaler {Henkelkoerper} und endliche {Erweiterungen} von {Schottky}-{Gruppen}},
 journal = {Comment. Math. Helv.},
 volume = {56},
 pages = {474--486},
 year = {1981}
}

\vfill

	\begin{flushleft}
		\textcolor{blue}{Naomi Andrew}\\
		Université Paris-Saclay, CNRS,  Laboratoire de mathématiques d'Orsay, 91405, Orsay, France \\
		\emph{e-mail:~}\texttt{naomi.andrew@universite-paris-saclay.fr}\\[4mm]
	\end{flushleft}

	\begin{flushleft}
		\textcolor{blue}{Yassine Guerch}\\
	Normandie Univ., UNICAEN, CNRS, Laboratoire de mathématiques Nicolas Oresme, 14000 Caen, France \\
		\emph{e-mail:~}\texttt{yassine.guerch@cnrs.fr}\\[4mm]
	\end{flushleft}

	\begin{flushleft}
		\textcolor{blue}{Camille Horbez}\\
	ENS de Lyon,  Unité de Mathématiques Pures et Appliquées, 69007, Lyon, France \\
		\emph{e-mail:~}\texttt{camille.horbez@ens-lyon.fr}\\[4mm]
	\end{flushleft}
	
		\begin{flushleft}
		\textcolor{blue}{Richard D. Wade}\\
        Mathematics Institute, Zeeman Building, University of Warwick, Coventry, CV4 7AL, UK
~\\
		\emph{e-mail:~}\texttt{ric.wade@warwick.ac.uk}\\[4mm]
	\end{flushleft}
	
	\noindent\textcolor{black!70}{\small Naomi Andrew and Camille Horbez acknowledge support from the European Research Council, through Grants~101040507 Artin-Out-ME-OA and~850930 Fibring. Views and opinions expressed are however those of the authors only and do not necessarily reflect those of the European Union or the European Research Council; neither the European Union nor the granting authority can be held responsible for them. Naomi Andrew and Richard D. Wade acknowledge support from the Royal Society of Great Britain through grant URF\textbackslash  R1\textbackslash 201297. Some of the research was carried out at the Isaac Newton Institute for Mathematical Sciences during the programme Operators, Graphs, Groups (EPSRC grant EP/Z000580/1s).}

\end{document}